\documentclass[12pt]{article}

\usepackage{amsmath}
\usepackage{amssymb}
\usepackage{amsfonts}
\usepackage{latexsym}
\usepackage{graphicx}

\newtheorem{proposition}{Proposition}[section]

\newtheorem{remark}{Remark}[section]

\date{ }

\begin{document}

\title{Penalization of stationary
Navier-Stokes equations and applications in topology optimization}

\author{Cornel Marius Murea$^1$, Dan Tiba$^2$\\
{\normalsize $^1$ D\'epartement de Math\'ematiques, IRIMAS,}\\
{\normalsize Universit\'e de Haute Alsace, France,}\\
{\normalsize cornel.murea@uha.fr}\\
{\normalsize $^2$ Institute of Mathematics (Romanian Academy) and}\\ 
{\normalsize Academy of Romanian Scientists, Bucharest, Romania,}\\ 
{\normalsize dan.tiba@imar.ro}
}

\maketitle

\begin{abstract}
We consider the steady Navier-Stokes system with mixed boundary conditions,
in subdomains of a holdall domain. We study, via the penalization method,
its approximation properties. Error estimates, obtained using the extension operator, other evaluations 
and the uniqueness of the solution, when the viscosity may be arbitrarily small
in certain subdomains, are also discussed.
Numerical tests, including topological optimization applications, are presented.
A general convergence result for the approximation of this type of geometric inverse problems and
of the associated optimal control problems, is investigated in the last part of the paper.

  
\end{abstract}

\section{Introduction and preliminaries\label{sec:1}}
\setcounter{equation}{0}

We consider a bounded domain $D\subset \mathbb{R}^d$, $d\in\{2,3\}$,
with Lipschitz boundary $\partial D=\overline{\Gamma}_D \cup \overline{\Gamma}_N$,
where $\Gamma _D$ and $\Gamma _N$ are disjoint, relatively open subsets.
Let $\mathbf{f} \in L^2(D)^d$ be given 
and $\Omega \subset \subset D$ (the obstacle) be an open set, not necessarily connected,
with Lipschitz boundary $\partial\Omega$. We assume that $\omega=D\setminus\overline{\Omega}$
is connected.
A fluid satisfying the Navier-Stokes equations occupies the domain $\omega$:
\begin{eqnarray}
-\nu \Delta \mathbf{y}_\omega + (\mathbf{y}_\omega\cdot\nabla)\mathbf{y}_\omega + \nabla p_\omega
& = & \mathbf{f}_\omega,\hbox{ in }\omega, \label{1.1}\\
\nabla \cdot \mathbf{y}_\omega & = & 0,\hbox{ in }\omega, \label{1.2}\\
\mathbf{y}_\omega & = & 0,\hbox{ on }\Gamma_D, \label{1.3}\\
\nu \frac{\partial\mathbf{y}_\omega}{\partial \mathbf{n}}-p_\omega\mathbf{n}
& = & \boldsymbol{\psi},\hbox{ on }\Gamma_N, \label{1.4}\\
\nu \frac{\partial\mathbf{y}_\omega}{\partial \mathbf{n}}-p_\omega\mathbf{n}
& = & 0,\hbox{ on }\partial\Omega, \label{1.5}
\end{eqnarray}
where $\mathbf{y}_\omega : \overline{\omega}\rightarrow \mathbb{R}^d$
is the velocity and $p_\omega : \overline{\omega}\rightarrow \mathbb{R}$ denotes the pressure
of the fluid.
Here $\nu>0$ is the kinematic viscosity, $\mathbf{n}$ is the unit
outward normal to the boundary $\partial\omega$, 
$\mathbf{f}_\omega : \omega\rightarrow \mathbb{R}^d$, $\mathbf{f}_\omega = \mathbf{f}\vert_{\omega}$ are the body forces 
and $\boldsymbol{\psi}:\Gamma_N\rightarrow \mathbb{R}^d$ is the given Neumann boundary
condition.
We have denoted by $\nabla \cdot \mathbf{y}=\sum_{i=1}^d \frac{\partial y_i}{\partial x_i}$
the divergence operator and
$(\mathbf{y}\cdot\nabla)\mathbf{y}=\sum_{j=1}^d y_j\frac{\partial \mathbf{y}}{\partial x_j}$.
We quote some classical textbooks on Navier-Stokes equations \cite{Temam1984},
\cite{Girault1986}, \cite{Galdi2011}, that we use in this paper.
The articles \cite{Zhou2014}, \cite{Zhou2016}, \cite{Zhou2017}  discuss related subjects and
methods, while \cite{bf} applies a different type of penalization in a general situation.
The extension operator technique is used in \cite{Zhou2017} for the Stokes equation. We also
mention the well known immersed boundary approach \cite{pes} that also uses penalization.
Such methods enter the class of fictitious domain methodology  and their applications in shape
optimization have been already considered in \cite{N_Tiba2012}, for linear state systems
(with Dirichlet boundary conditions). 
Here, we consider mixed boundary conditions for the stationary Navier-Stokes equation.
Our results (uniqueness, error estimates) and the
applications to the difficult topology optimization question, are certainly new.

We set
$$
W_\omega=\left\{ \mathbf{w}_\omega\in H^1(\omega)^d;\ \mathbf{w}_\omega=0\hbox{ on }\Gamma_D\right\},
\quad
Q_\omega=L^2(\omega)
$$
and
$V_\omega=\left\{ \mathbf{w}_\omega\in W_\omega; \nabla \cdot \mathbf{w}_\omega=0 \hbox{ in }\omega\right\}$.
To shorten notation, we write $\| \mathbf{w}_\omega \|_{1,\omega} =\| \mathbf{w}_\omega \|_{H^1(\omega)^d}$,
$\vert \mathbf{w}_\omega \vert_{1,\omega} =\| \nabla \mathbf{w}_\omega \|_{L^2(\omega)^{d\times d}}$ and
$\| q_\omega \|_{0,\omega} =\| q_\omega \|_{L^2(\omega)}$.

From the generalized Poincar\'e inequality, \cite{Boyer2013}, p. 179, Proposition III.3.28, we obtain
\begin{equation}\label{1.poincare}
\| \mathbf{w}_\omega \|_{L^2(\omega)^d} \leq C_P \| \nabla \mathbf{w}_\omega \|_{L^2(\omega)^{d\times d}},\quad \forall
\mathbf{w}_\omega\in W_\omega
\end{equation}
where $C_P >0$ depends on $\omega$ and $d$.
Consequently,
$\vert \mathbf{w}_\omega \vert_{1,\omega}\leq \| \mathbf{w}_\omega \|_{1,\omega}
\leq  \sqrt{C_P^2+1}\, \vert\mathbf{w}_\omega \vert_{1,\omega}$, i.e.
the norms $\vert \cdot \vert_{1,\omega}$ and $\| \cdot \|_{1,\omega}$ are
equivalent on $W_\omega$.

We introduce
\begin{eqnarray*}
&& 
a_\omega : W_\omega \times W_\omega \rightarrow \mathbb{R},\quad
a_\omega(\mathbf{v}_\omega,\mathbf{w}_\omega)=
\nu \int_{\omega} \nabla \mathbf{v}_\omega : \nabla\mathbf{w}_\omega \,d\mathbf{x}, \\
&& 
b_\omega : W_\omega \times Q_\omega \rightarrow \mathbb{R},\quad
b_\omega(\mathbf{v}_\omega,q_\omega) = 
- \int_{\omega}\left(\nabla\cdot\mathbf{v}_\omega\right)q_\omega \,d\mathbf{x}, \\
&& 
c_\omega : W_\omega \times W_\omega \times W_\omega \rightarrow \mathbb{R},\quad
c_\omega(\mathbf{u}_\omega,\mathbf{v}_\omega,\mathbf{w}_\omega)  = 
\int_{\omega} \left[ 
(\mathbf{u}_\omega\cdot\nabla ) \mathbf{v}_\omega
\right] \cdot \mathbf{w}_\omega \,d\mathbf{x}.
\end{eqnarray*}
We have denoted by ``$\cdot$'' the scalar product in
$\mathbb{R}^d$ with the associated norm $\vert\cdot\vert$
and by ``:'' the Frobenius matrix product $(a_{ij}):(b_{ij})=\sum_{i,j=1}^d a_{ij}b_{ij}$.

For $\mathbf{u}_\omega\in V_\omega$ and $\mathbf{w}_\omega\in W_\omega$, integrating by parts we get
\begin{eqnarray*}
  c_\omega(\mathbf{u}_\omega,\mathbf{w}_\omega,\mathbf{w}_\omega)  &=&
  \int_\omega \sum_{i=1}^d\sum_{j=1}^d u_j\frac{\partial w_i}{\partial x_j}w_i \,d\mathbf{x}
  = \frac{1}{2}\int_\omega \sum_{i=1}^d\sum_{j=1}^d u_j \frac{\partial (w_i)^2}{\partial x_j}
  \,d\mathbf{x}\\
  &=& - \int_\omega (\nabla \cdot \mathbf{u}_\omega) \vert\mathbf{w}_\omega\vert^2\,d\mathbf{x}
  +\int_{\partial\omega} \mathbf{u}_\omega\cdot \mathbf{n} \vert\mathbf{w}_\omega\vert^2\,ds\\
  &=&\int_{\Gamma_N\cup \partial\Omega} \mathbf{u}_\omega\cdot \mathbf{n} \vert\mathbf{w}_\omega\vert^2\,ds.
\end{eqnarray*}
The last expression does not necessarily vanish, due to the boundary
conditions on $\partial\omega$.
We introduce $\tilde{c}_\omega : W_\omega \times W_\omega \times W_\omega \rightarrow\mathbb{R}$,
$$
\tilde{c}_\omega(\mathbf{u}_\omega,\mathbf{v}_\omega,\mathbf{w}_\omega)=
\frac{1}{2}c_\omega(\mathbf{u}_\omega,\mathbf{v}_\omega,\mathbf{w}_\omega)
-\frac{1}{2}c_\omega(\mathbf{u}_\omega,\mathbf{w}_\omega,\mathbf{v}_\omega)
$$
and we have
$$
\tilde{c}_\omega(\mathbf{u}_\omega,\mathbf{w}_\omega,\mathbf{w}_\omega)=0,
\quad\forall \mathbf{u}_\omega,\mathbf{w}_\omega\in W_\omega.
$$
Consequently, for all $\mathbf{u}_\omega,\mathbf{v}_\omega,\mathbf{w}_\omega\in W_\omega$, we have
\begin{eqnarray*}
0&=&\tilde{c}_\omega(\mathbf{u}_\omega,\mathbf{v}_\omega+\mathbf{w}_\omega,\mathbf{v}_\omega+\mathbf{w}_\omega)
\\
&=&\tilde{c}_\omega(\mathbf{u}_\omega,\mathbf{v}_\omega,\mathbf{v}_\omega)
+\tilde{c}_\omega(\mathbf{u}_\omega,\mathbf{v}_\omega,\mathbf{w}_\omega)
+\tilde{c}_\omega(\mathbf{u}_\omega,\mathbf{w}_\omega,\mathbf{v}_\omega)
+\tilde{c}_\omega(\mathbf{u}_\omega,\mathbf{w}_\omega,\mathbf{w}_\omega)
\\
&=&\tilde{c}_\omega(\mathbf{u}_\omega,\mathbf{v}_\omega,\mathbf{w}_\omega)
+\tilde{c}_\omega(\mathbf{u}_\omega,\mathbf{w}_\omega,\mathbf{v}_\omega)
\end{eqnarray*}
and finally
$$
\tilde{c}_\omega(\mathbf{u}_\omega,\mathbf{v}_\omega,\mathbf{w}_\omega)
=-\tilde{c}_\omega(\mathbf{u}_\omega,\mathbf{w}_\omega,\mathbf{v}_\omega).
$$
Reciprocally, if the above equality holds for all
$\mathbf{u}_\omega,\mathbf{v}_\omega,\mathbf{w}_\omega\in W_\omega$, in particular
for $\mathbf{v}_\omega=\mathbf{w}_\omega$, we get
$$
\tilde{c}_\omega(\mathbf{u}_\omega,\mathbf{w}_\omega,\mathbf{w}_\omega)
=-\tilde{c}_\omega(\mathbf{u}_\omega,\mathbf{w}_\omega,\mathbf{w}_\omega)
\Leftrightarrow
\tilde{c}_\omega(\mathbf{u}_\omega,\mathbf{w}_\omega,\mathbf{w}_\omega)=0
$$
for all $\mathbf{u}_\omega,\mathbf{w}_\omega\in W_\omega$.

Clearly $\mathbf{f}_\omega\in L^2(\omega)^d$ and we assume that $\boldsymbol{\psi}\in H^{1/2}(\Gamma_N)^d$,
we define $F_\omega\in W_\omega^\prime$, where $ W_\omega^\prime$ is the dual space of $W_\omega$, by
$$
\langle F_\omega, \mathbf{w}_\omega\rangle_{*,\omega} =
\int_\omega \mathbf{f}_\omega\cdot \mathbf{w}_\omega\,d\mathbf{x}
+\int_{\Gamma_N} \boldsymbol{\psi}\cdot \mathbf{w}_\omega\,ds
$$
and we set
$$
\vert F_\omega \vert_{*,\omega} = \sup_{\mathbf{w}_\omega\in W_\omega \setminus \{0\} }
\frac{\langle F_\omega, \mathbf{w}_\omega\rangle_{*,\omega}}{\vert\mathbf{w}_\omega \vert_{1,\omega}},
$$
where $\langle F_\omega , \mathbf{w}_\omega \rangle_{*,\omega}=F_\omega(\mathbf{w}_\omega)$
is the value of $F_\omega$ in $\mathbf{w}_\omega$.

From \cite{Girault1986}, Chapter IV, Lemma 2.1, p. 284, 
the application $c_\omega$ is continuous on $H^1(\omega)^d \times H^1(\omega)^d \times H^1(\omega)^d$, then
$\tilde{c}_\omega$ is continuous too. We set
$$
\vert \tilde{c}_\omega \vert_{1,\omega} = \sup_{\mathbf{u}_\omega,\mathbf{v}_\omega,\mathbf{w}_\omega \in W_\omega \setminus \{0\} }
\frac{\left\vert \tilde{c}_\omega(\mathbf{u}_\omega,\mathbf{v}_\omega,\mathbf{w}_\omega) \right\vert}
{\vert\mathbf{u}_\omega \vert_{1,\omega} \vert\mathbf{v}_\omega \vert_{1,\omega} \vert\mathbf{w}_\omega \vert_{1,\omega}}.
$$
Estimations of $\vert c_\omega\vert_{1,\omega}$ in the case of homogeneous Dirichlet boundary conditions
all over the boundary $\partial\omega$ can be found in \cite{Galdi2011}, Lemma IX.1.1,
p. 588.

Using \cite{MR_2007}, Lemma 5.1, the application
$
\mathbf{w}_\omega \in W_\omega \rightarrow \nabla\cdot\mathbf{w}_\omega \in Q_\omega
$
is onto. It follows, see
\cite{Girault1986}, Chapter I, Lemma 4.1, p. 58, 
that the below \textit{inf-sup} condition holds
$$
\exists \beta_\omega >0,\quad
\beta_\omega \leq \inf_{q_\omega\in Q_\omega \setminus \{0\} }\sup_{\mathbf{w}_\omega\in W_\omega \setminus \{0\} }
\frac{b_\omega(\mathbf{w}_\omega,q_\omega)}{\| q_\omega\|_{0,\omega}\vert\mathbf{w}_\omega \vert_{1,\omega}}.
$$

The weak variational formulation of (\ref{1.1})-(\ref{1.5}) is:
find $\mathbf{y}_\omega\in W_\omega$ and $p_\omega\in Q_\omega$ such that
\begin{eqnarray}
a_\omega(\mathbf{y}_\omega,\mathbf{w}_\omega)
+b_\omega(\mathbf{w}_\omega,p_\omega)
+\tilde{c}_\omega(\mathbf{y}_\omega,\mathbf{y}_\omega,\mathbf{w}_\omega)
&=&\langle F_\omega, \mathbf{w}_\omega\rangle_{*,\omega},\forall \mathbf{w}_\omega\in W_\omega
\label{1.ns1}\\
b_\omega(\mathbf{y}_\omega,q_\omega)&=&0,\forall q_\omega\in Q_\omega.
\label{1.ns2}
\end{eqnarray}
The condition (\ref{1.ns2}) is equivalent to $\mathbf{y}_\omega\in V_\omega$: if we put
$q_\omega=\nabla \cdot \mathbf{y}_\omega$ in (\ref{1.ns2}), we get $\nabla \cdot \mathbf{y}_\omega=0$
in $\omega$ and reciprocally, if $\mathbf{y}_\omega\in V_\omega$, then (\ref{1.ns2}) holds.

We also consider the problem: find $\mathbf{y}_\omega\in V_\omega$ such that
\begin{equation}\label{1.ns3}
a_\omega(\mathbf{y}_\omega,\mathbf{v}_\omega)
+\tilde{c}_\omega(\mathbf{y}_\omega,\mathbf{y}_\omega,\mathbf{v}_\omega)
=\langle F_\omega, \mathbf{v}_\omega\rangle_{*,\omega},\quad \forall \mathbf{v}_\omega\in V_\omega.
\end{equation}
If $(\mathbf{y}_\omega,p_\omega)\in W_\omega\times Q_\omega$ is a solution of (\ref{1.ns1})-(\ref{1.ns2}),
then $\mathbf{y}_\omega\in V_\omega$ is a solution of (\ref{1.ns3}).
Reciprocally, if $\mathbf{y}_\omega\in V_\omega$ is a solution of (\ref{1.ns3}),
there exists a unique $p_\omega\in Q_\omega$ such that $(\mathbf{y}_\omega,p_\omega)$ is a solution
of (\ref{1.ns1})-(\ref{1.ns2}), see \cite{Girault1986}, Chapter IV, Theorem 1.4, p. 283.

Notice that multiplying (\ref{1.1}) by $\mathbf{w}_\omega\in W_\omega$ and integrating by parts, we get
(\ref{1.ns1}) with $c_\omega$ in place of $\tilde{c}_\omega$.
We point out that $c_\omega(\mathbf{u}_\omega,\mathbf{v}_\omega,\mathbf{w}_\omega)
\neq -c_\omega(\mathbf{u}_\omega,\mathbf{w}_\omega,\mathbf{v}_\omega)$ for 
$\mathbf{u}_\omega,\mathbf{v}_\omega,\mathbf{w}_\omega\in V_\omega$, then 
$c_\omega(\mathbf{y}_\omega,\mathbf{y}_\omega,\mathbf{v}_\omega)\neq
\tilde{c}_\omega(\mathbf{y}_\omega,\mathbf{y}_\omega,\mathbf{v}_\omega)$, for $\mathbf{v}_\omega\in V_\omega$.
Then, if we replace $\tilde{c}_\omega$ by $c_\omega$ in (\ref{1.ns3}), we do not obtain an
equivalent formulation. This replacement technique is used and accepted for the finite element approximation $\mathbf{u}_h$ when 
$\nabla\cdot\mathbf{u}_h\neq 0$ in $\omega$, see \cite{Temam1984}, Chap. 2, Sect. 3
or for certain boundary conditions, see \cite{Zhou2016}, and we also use it in this paper.

Multiplying (\ref{1.2}) by $q_\omega\in Q_\omega$ and integrating over $\omega$, we get
(\ref{1.ns2}).

\begin{remark}\label{rem:1.1}
Let $\mathbf{u}^n=(u_i^n)_{1\leq i \leq d}$, $n\in\mathbb{N}$, be a weakly convergent sequence to
$\mathbf{u}=(u_i)_{1\leq i \leq d}$ in $V_\omega$. 
From the Sobolev embedding theorem, for $d\in\{2,3\}$, the injection
$H^1(\omega)\subset L^4(\omega)$ is compact,
then $u_i^n$ converges strongly to $u_i$ in $L^{4}(\omega)$ and 
for $\mathbf{v}\in V_\omega\subset L^4(\omega)^d$, we get
$u_i^n v_j$ converges strongly to $u_i v_j$ in $L^{2}(\omega)$, $1\leq i \leq d$, $1\leq j \leq d$.
We also have that $\frac{\partial u_j^n}{\partial x_i}$ converges weakly to
$\frac{\partial u_j}{\partial x_i}$ in $L^{2}(\omega)$.
It follows
$$
\lim_{n\rightarrow\infty}
c_\omega(\mathbf{u}^n,\mathbf{u}^n,\mathbf{v})
=\lim_{n\rightarrow\infty}{\sum_{i,j=1}^d \int_\omega u_i^n \frac{\partial u_j^n}{\partial x_i} v_j\, d\mathbf{x}}
=c_\omega(\mathbf{u},\mathbf{u},\mathbf{v}).
$$
Similarly, $u_i^n u_j^n$ converges strongly to $u_i u_j$ in $L^{2}(\omega)$. For
$\mathbf{v}\in V_\omega\subset L^2(\omega)^d$ we get that
$$
\lim_{n\rightarrow\infty}
c_\omega(\mathbf{u}^n,\mathbf{v},\mathbf{u}^n)
=\lim_{n\rightarrow\infty}{\sum_{i,j=1}^d \int_\omega u_i^n \frac{\partial v_j}{\partial x_i} u_j^n\, d\mathbf{x}}
=c_\omega(\mathbf{u},\mathbf{v},\mathbf{u}).
$$
The free divergence condition or the boundary conditions of $\mathbf{u}^n$, $\mathbf{u}$, $\mathbf{v}$
are not mandatory to get this property. Similar convergences are valid for $\tilde{c}_\omega$.
\end{remark}

We can apply \cite{Girault1986}, Chapter IV, Theorem 1.2, p. 280 and get that
the problem (\ref{1.ns3})
has at least one solution in $V_\omega$. Moreover, we have the same estimations as in the case
of Dirichlet homogeneous boundary conditions
\begin{eqnarray}
\vert \mathbf{y}_\omega \vert_{1,\omega} &\leq & \frac{\vert F_\omega \vert_{*,\omega}}{\nu}, \label{1.vit}\\
\| p_\omega \|_{0,\omega} &\leq &  \frac{1}{\beta_\omega}
\left( 2\vert F_\omega \vert_{*,\omega} + \frac{ \vert \tilde{c}_\omega \vert_{1,\omega} \vert F_\omega \vert_{*,\omega}^2 }{\nu^2}
\right)\label{1.pres}
\end{eqnarray}
where $p_\omega$ is defined in (\ref{1.ns1}), see \cite{Boyer2013}, Theorem V.3.1, p. 392.

The uniqueness results for homogeneous Dirichlet boundary condition from
\cite{Quarteroni}, Theorem 10.1.1, p. 341 can be adapted for the mixed boundary
conditions (\ref{1.ns1})-(\ref{1.ns2}).
The idea is to prove that the application
$$
\mathbf{w}\in V_\omega \rightarrow \mathbf{u}_\mathbf{w}\in V_\omega
$$
is a contraction from
$$
K=\left\{\mathbf{w}\in V_\omega ;\  \vert \mathbf{w} \vert_{1,\omega} \leq  \frac{\vert F_\omega \vert_{*,\omega}}{\nu}
\right\}
$$
to $K$, where $\mathbf{u}_\mathbf{w}\in V_\omega$ is the solution of
$$
a_\omega(\mathbf{u}_\mathbf{w},\mathbf{v}_\omega)
+\tilde{c}_\omega(\mathbf{w},\mathbf{u}_\mathbf{w},\mathbf{v}_\omega)
=\langle F_\omega, \mathbf{v}_\omega\rangle_{*,\omega},\quad \forall \mathbf{v}_\omega\in V_\omega.
$$
Then, we apply the Banach fixed point theorem and obtain:

\begin{proposition}\label{prop:1.1}
If
\begin{equation}\label{1.H1}
  \frac{ \vert \tilde{c}_\omega \vert_{1,\omega} \vert F_\omega \vert_{*,\omega} }{\nu^2} < 1
\end{equation}
then the problem (\ref{1.ns1})-(\ref{1.ns2}) has a unique solution
$\mathbf{y}_\omega\in V_\omega$ and $p_\omega\in Q_\omega$.
\end{proposition}

Let us introduce 
\begin{equation}\label{1.13}
V_1=\left\{
\mathbf{v}^1\in V_\omega;\ 
\int_{\partial \Omega_i} \mathbf{v}^1\cdot \mathbf{n}\, ds=0,
\hbox{ for each } \Omega_i \hbox{ connected comp. of } \Omega
\right\}, 
\end{equation}
a Hilbert subspace of $V_\omega$, and the problem: find $\widehat{\mathbf{y}}_\omega\in V_1$
such that
\begin{equation}\label{1.ns3_v1}
a_\omega(\widehat{\mathbf{y}}_\omega,\mathbf{v}^1)
+\tilde{c}_\omega(\widehat{\mathbf{y}}_\omega,\widehat{\mathbf{y}}_\omega,\mathbf{v}^1)
=\langle F_\omega, \mathbf{v}^1\rangle_{*,\omega},\quad \forall \mathbf{v}^1\in V_1.
\end{equation}

If $\mathbf{v}^1\in V_1$, then $\int_{\Gamma_N} \mathbf{v}^1\cdot \mathbf{n}\, ds=0$, too. In fact,
$$
0=\int_{\omega} \nabla \cdot \mathbf{v}^1\, d\mathbf{x}
=\int_{\partial D} \mathbf{v}^1\cdot \mathbf{n}\, ds 
+\sum_i \int_{\partial \Omega_i} \mathbf{v}^1\cdot \mathbf{n}\, ds
=\int_{\partial D} \mathbf{v}^1\cdot \mathbf{n}\, ds
=\int_{\Gamma_N} \mathbf{v}^1\cdot \mathbf{n}\, ds.
$$

Putting $\mathbf{v}^1=\widehat{\mathbf{y}}_\omega$  in (\ref{1.ns3_v1}), we get
\begin{equation}\label{1.15}
\vert \widehat{\mathbf{y}}_\omega \vert_{1,\omega} \leq  \frac{\vert F_\omega \vert_{*,\omega}}{\nu}.
\end{equation}

\begin{proposition}\label{prop:1.2}
Under the hypotheses (\ref{1.H1}),
the problem (\ref{1.ns3_v1}) has a unique solution
$\widehat{\mathbf{y}}_\omega\in V_1$.
\end{proposition}

The proof is as before, we just replace $V_\omega$ by $V_1$.
We 
apply the Banach 
fixed point theorem in the closed set below, with the norm $\vert \cdot \vert_{1,\omega}$:
$$
K_1=\left\{\mathbf{v}\in V_1 ;\  \vert \mathbf{v} \vert_{1,\omega} \leq  \frac{\vert F_\omega \vert_{*,\omega}}{\nu}
\right\}.
$$

We can introduce 
$$
W_1=\left\{
\mathbf{w}^1\in W_\omega;\ 
\int_{\partial \Omega_i} \mathbf{w}^1\cdot \mathbf{n}\, ds=0,
\hbox{ for each } \Omega_i \hbox{ connected comp. of } \Omega
\right\} 
$$
and $Q_1$ the image of the divergence operator 
$\mathbf{w}^1 \in W_1 \rightarrow \nabla\cdot\mathbf{w}_\omega \in Q_\omega$. 
From \cite{Girault1986}, Chap. 1, Sect. 4,
for each $\widehat{\mathbf{y}}_\omega\in V_1$ solution of (\ref{1.ns3_v1}), there exists only one
$\widehat{p}_\omega\in Q_1$, such that
\begin{equation}\label{1.ns1_v1}
a_\omega(\widehat{\mathbf{y}}_\omega,\mathbf{w}^1)
+b_\omega(\mathbf{w}^1,\widehat{p}_\omega)
+\tilde{c}_\omega(\widehat{\mathbf{y}}_\omega,\widehat{\mathbf{y}}_\omega,\mathbf{w}^1)
=\langle F_\omega, \mathbf{w}^1\rangle_{*,\omega},\quad \forall \mathbf{w}^1\in W_1.
\end{equation}

We give an interpretation of (\ref{1.ns1_v1}).
Suppose $\widehat{\mathbf{y}}_\omega\in V_1\cap H^2(\omega)^d$ and
$\widehat{p}_\omega\in Q_1\cap H^1(\omega)$ satisfy (\ref{1.ns1_v1}),
with $c_\omega$ in place of $\tilde{c}_\omega$.
As before, using $\mathbf{w}_1\in \mathcal{D}(\omega)$, we get that (\ref{1.1}) holds in $L^2(\omega)^d$.
From $\widehat{\mathbf{y}}_\omega\in V_1$, (\ref{1.2}) and (\ref{1.3}) hold.
If $\mathbf{w}_1\in W_1$, we also get
$$
\int_{\partial \Omega_i}
\left( \nu \frac{\partial \widehat{\mathbf{y}}_\omega}{\partial \mathbf{n}}
-\widehat{p}_\omega\mathbf{n}\right)
\cdot \mathbf{w}_1\, ds=0
$$
and
$$
\int_{\Gamma_N}
\left( \nu \frac{\partial \widehat{\mathbf{y}}_\omega}{\partial \mathbf{n}}
-\widehat{p}_\omega\mathbf{n}\right)
\cdot \mathbf{w}_1\, ds
=\int_{\Gamma_N}\boldsymbol{\psi}\cdot \mathbf{w}_1\, ds.
$$

The boundary conditions (\ref{1.4}), (\ref{1.5}) hold 
in the dual space of the traces on $\partial \Omega_i$ and $\Gamma_N$ of the test functions from $W_1$.
In addition to (\ref{1.1})-(\ref{1.5}), the solution $\widehat{\mathbf{y}}_\omega$ satisfies the relations
$\int_{\partial \Omega_i} \widehat{\mathbf{y}}_\omega\cdot \mathbf{n}\, ds=0,
\hbox{ for each } \Omega_i \hbox{ connected component of } \Omega$.

Reciprocally, let $\mathbf{y}_\omega\in H^2(\omega)^d$ and
$p_\omega\in H^1(\omega)$ be a solution of (\ref{1.1})-(\ref{1.5}) such that
$\int_{\partial \Omega_i} \mathbf{y}_\omega\cdot \mathbf{n}\, ds=0,
\hbox{ for each } \Omega_i \hbox{ connected component of } \Omega$. Integrating
by parts, we obtain (\ref{1.ns1_v1}) with $c_\omega$ in the place of $\tilde{c}_\omega$ (that are different in general).

In this article, due to the methodology of fixed domain type that we use, the problem 
(\ref{1.ns3_v1}) will play a central role. This is adapted to the shape and topology optimization
problems that we associate to the Navier-Stokes system and we solve by penalizing the state equation in $D$ and
working just in $D$ at the approximating level.

The plan of the paper is as follows.
In the next section we study the penalized Navier-Stokes system
in $D$, its approximation properties (with error estimates) and the uniqueness of
the solution, other evaluations, that are obtained in a non standard setting.
For the discrete and the penalized/regularized Navier-Stokes system in $D$, a convergence
result for its solution is obtained. In the last section,
several numerical examples are indicated, stressing the topological optimization question and discussing a general convergence property.

Our results provide new insight in the uniqueness theory of Navier-Stokes, in their approximation and a
new effective way to solve the  difficult topological optimization problems associated to them. 

The scientific literature contains works devoted mainly to shape
optimization problems
for the Navier-Stokes system and we quote just the books \cite{mp} (that also includes brief presentations of some topology optimization approaches), \cite{PS} and the
articles \cite{F}, \cite{DFOP}, \cite{SN2022}. Concerning technical applications of topology optimization for
fluids, we mention the very recent paper \cite{chen} and its references.
In the case of second order linear elliptic equations, topology optimization via
penalization methods were discussed in \cite{Tiba2009} and in the survey \cite{N_Tiba2012}.
For higher order equations, we indicate the recent article \cite{MT_2019}. The classical topology
optimization methodology is discussed in \cite{NST}, Ch. 5.3 and in the dedicated monograph \cite{NOS}.

\section{Fixed domain approach}
\setcounter{equation}{0}

\subsection{Penalization}

We set $Q=L^2(D)$,
$$
W=\left\{ \mathbf{w}\in H^1(D)^d;\ \mathbf{w}=0\hbox{ on }\Gamma_D\right\},\quad
V=\left\{ \mathbf{w}\in W; \nabla \cdot \mathbf{w}=0 \hbox{ in }D\right\}
$$
and for $ \epsilon > 0$, we introduce (here, $a_\omega(\mathbf{v},\mathbf{w})$ means
$a_\omega(\mathbf{v}_{\vert\omega},\mathbf{w}_{\vert\omega})$, etc.)

\begin{eqnarray*}
&& 
a^\epsilon : W \times W \rightarrow \mathbb{R},\quad
a^\epsilon(\mathbf{v},\mathbf{w})=
a_\omega(\mathbf{v},\mathbf{w})+\epsilon a_\Omega(\mathbf{v},\mathbf{w}), \\
&&
b^\epsilon : W \times Q \rightarrow \mathbb{R},\quad
b^\epsilon(\mathbf{w},q)=
b_\omega(\mathbf{w},q)+\epsilon b_\Omega(\mathbf{w},q), \\
&& 
\tilde{c}^\epsilon : W \times W \times W \rightarrow \mathbb{R},\quad
\tilde{c}^\epsilon(\mathbf{u},\mathbf{v},\mathbf{w})  = 
\tilde{c}_\omega(\mathbf{u},\mathbf{v},\mathbf{w})
+\epsilon\tilde{c}_\Omega(\mathbf{u},\mathbf{v},\mathbf{w})
\end{eqnarray*}
where
$\tilde{c}_\Omega: H^1(\Omega)^d\times H^1(\Omega)^d \times H^1(\Omega)^d
\rightarrow \mathbb{R}$ and $b_\Omega: H^1(\Omega)^d\times L^2(\Omega)\rightarrow \mathbb{R}$
are defined
similarly to $\tilde{c}_\omega$ and $b_\omega$ but integrating over $\Omega$ and
$a_\Omega:H^1(\Omega)^d \times H^1(\Omega)^d
\rightarrow \mathbb{R}$ is defined by
$$
a_\Omega(\mathbf{v}_\Omega,\mathbf{w}_\Omega)
=\int_\Omega \nabla \mathbf{v}_\Omega : \nabla \mathbf{w}_\Omega\, d\mathbf{x}
+\int_\Omega \mathbf{v}_\Omega \cdot \mathbf{w}_\Omega\, d\mathbf{x},
$$
where $\mathbf{v}_\Omega, \mathbf{w}_\Omega \in H^1(\Omega)^d$.  

We set
$$
\| \tilde{c}_\Omega \|_{1,\Omega} = \sup_{\mathbf{u}_\Omega,\mathbf{v}_\Omega,\mathbf{w}_\Omega \in H^1(\Omega)^d\setminus \{0\} }
\frac{\left\vert \tilde{c}_\Omega(\mathbf{u}_\Omega,\mathbf{v}_\Omega,\mathbf{w}_\Omega) \right\vert}
{\|\mathbf{u}_\Omega \|_{1,\Omega} \|\mathbf{v}_\Omega \|_{1,\Omega} \|\mathbf{w}_\Omega \|_{1,\Omega}}.
$$
We point out that for $H^1(\Omega)^d$, we use the norm $\|\cdot \|_{1,\Omega}$, while
for $W_\omega$, we use $\vert\cdot \vert_{1,\omega}$.

The approximating problem written in the fixed domain $D$ is: find $\mathbf{y}^\epsilon\in V$ such that
\begin{equation}\label{2.ns3}
a^\epsilon(\mathbf{y}^\epsilon,\mathbf{v})
+\tilde{c}^\epsilon(\mathbf{y}^\epsilon,\mathbf{y}^\epsilon,\mathbf{v})
=\langle F_\omega, \mathbf{v}\rangle_{*,\omega},\quad \forall \mathbf{v}\in V.
\end{equation}
For $\mathbf{v}\in V$, $\langle F_\omega, \mathbf{v}\rangle_{*,\omega}$ means
$\langle F_\omega, \mathbf{v}_{\vert\omega}\rangle_{*,\omega}$.

\begin{proposition}\label{prop:2.1}
For each $ \epsilon > 0$, the problem (\ref{2.ns3}) has at least one solution and
\begin{eqnarray}
\vert \mathbf{y}^\epsilon \vert_{1,\omega} &\leq & \frac{\vert F_\omega \vert_{*,\omega}}{\nu}, \label{2.2}\\
\sqrt{\epsilon} \| \mathbf{y}^\epsilon \|_{1,\Omega} &\leq &\frac{\vert F_\omega \vert_{*,\omega}}{\sqrt{\nu}}.
\label{2.3}
\end{eqnarray}
\end{proposition}

\noindent
\textit{Proof.} 
For $\mathbf{w}\in W$, we have
\begin{eqnarray*}
a^\epsilon(\mathbf{w},\mathbf{w})&=&\nu \vert \mathbf{w} \vert_{1,\omega}^2
+\epsilon \| \mathbf{w} \|_{1,\Omega}^2 \geq
\frac{\nu}{C_P^2+1} \| \mathbf{w} \|_{1,\omega}^2 + \epsilon  \| \mathbf{w} \|_{1,\Omega}^2\\
&\geq &\min\left(\frac{\nu}{C_P^2+1},\epsilon\right) \| \mathbf{w} \|_{1,D}^2
\end{eqnarray*}
and for $\mathbf{v},\mathbf{w}\in W$, we have
$$
\tilde{c}^\epsilon(\mathbf{v},\mathbf{w},\mathbf{w})=
\tilde{c}_\omega(\mathbf{v},\mathbf{w},\mathbf{w})
+\epsilon\tilde{c}_\Omega(\mathbf{v},\mathbf{w},\mathbf{w})=0.
$$ We prove that $\tilde{c}$ is continuous.
\begin{eqnarray*}
  \vert\tilde{c}^\epsilon(\mathbf{u},\mathbf{v},\mathbf{w}) \vert
  &\leq &
  \vert \tilde{c}_\omega(\mathbf{u},\mathbf{v},\mathbf{w}) \vert +
  \epsilon  \vert \tilde{c}_\Omega(\mathbf{u},\mathbf{v},\mathbf{w}) \vert
  \\
  &\leq & 
  \vert \tilde{c}_\omega \vert_{1,\omega}
  \vert\mathbf{u} \vert_{1,\omega} \vert\mathbf{v} \vert_{1,\omega} \vert\mathbf{w} \vert_{1,\omega}
  + \epsilon \| \tilde{c}_\Omega \|_{1,\Omega}
  \|\mathbf{u} \|_{1,\Omega} \|\mathbf{v} \|_{1,\Omega} \|\mathbf{w} \|_{1,\Omega}
  \\ &\leq &
  ( \vert\tilde{c}_\omega \vert_{1,\omega} +\epsilon \| \tilde{c}_\Omega \|_{1,\Omega} )
\|\mathbf{u} \|_{1,D} \|\mathbf{v} \vert \|_{1,D} \|\mathbf{w} \|_{1,D}.
\end{eqnarray*}
We can prove that $a$ is continuous, too.  Using Rem. \ref{1.1}, for all $\mathbf{v}\in V$,
we have
$$
\lim_{n\rightarrow\infty}
\tilde{c}^\epsilon(\mathbf{u}^n,\mathbf{u}^n,\mathbf{v})
=\tilde{c}^\epsilon(\mathbf{u},\mathbf{u},\mathbf{v})
$$
if $\mathbf{u}^n$ converges weakly to $\mathbf{u}$ in $V$. Once again,
we can apply \cite{Girault1986}, Theorem 1.2, p. 280 
and we obtain that the problem (\ref{2.ns3}) has at least one solution.

Putting $\mathbf{v}=\mathbf{y}^\epsilon$ in (\ref{2.ns3}), we obtain
$$
\nu \vert \mathbf{y}^\epsilon \vert_{1,\omega}^2
+\epsilon \| \mathbf{y}^\epsilon \|_{1,\Omega}^2
=a^\epsilon(\mathbf{y}^\epsilon,\mathbf{y}^\epsilon)
\leq \vert F \vert_{*,\omega} \vert \mathbf{y}^\epsilon \vert_{1,\omega}.
$$
It follows $\nu \vert \mathbf{y}^\epsilon \vert_{1,\omega}^2\leq
\vert F_\omega \vert_{*,\omega} \vert \mathbf{y}^\epsilon \vert_{1,\omega}$
which is equivalent to (\ref{2.2}) and
$$
\epsilon \| \mathbf{y}^\epsilon \|_{1,\Omega}^2 \leq  \vert F_\omega \vert_{*,\omega} \vert \mathbf{y}^\epsilon \vert_{1,\omega}
\leq \vert F_\omega \vert_{*,\omega} \frac{\vert F_\omega \vert_{*,\omega}}{\nu}
$$
which is equivalent to (\ref{2.3}).
 \quad $\Box$

\begin{remark}\label{rem:2.1bis}
An element $\mathbf{v}_\omega\in V_\omega$ can be extended to $\mathbf{v}\in V$ (i.e. preserving the divergence free property) if and only
if
$\int_{\partial \Omega_i} \mathbf{v}_\omega\cdot \mathbf{n}_{\partial \Omega_i}\, ds=0$
for each $\Omega_i$ connected component of $\Omega$, see \cite{Galdi2011}, Corollary III.3.1, p. 181. 
This condition is obviously satisfied by any $\mathbf{v}^1\in V_1$.  
Then, there exists an extension operator 
$E:V_1\rightarrow V$, such that
\begin{eqnarray}
E(\mathbf{v}^1) & = & \mathbf{v}^1,\hbox{ in }\omega \label{2.4}\\
\| E(\mathbf{v}^1) \|_{1,D} &\leq & C_E \| \mathbf{v}^1 \|_{1,\omega}
\label{2.5}
\end{eqnarray}
where $C_E>0$ depends on $\omega$, $d$. Conversely, for any $\mathbf{v}\in V$, we have 
$0=\int_{\Omega_i} \nabla\cdot\mathbf{v}\,d\mathbf{x}
= \int_{\partial \Omega_i} \mathbf{v}\cdot \mathbf{n}_{\partial \Omega_i}\, ds$, then
$\mathbf{v}_{\vert\omega}\in V_1$. 
\end{remark}

\begin{proposition}\label{prop:2.2}
Let $0 <\theta <1$ be fixed. For each $0 < \epsilon$, if
\begin{equation}
\| \tilde{c}_\Omega \|_{1,\Omega}
C_E\sqrt{C_P^2+1}\frac{\vert F_\omega \vert_{*,\omega}}{\nu} \leq  \theta
\label{2.6}
\end{equation}  
then, $\mathbf{y}^\epsilon $, solution of (\ref{2.ns3}), satisfies
\begin{equation}
  \| \mathbf{y}^\epsilon \|_{1,\Omega} \leq
  \frac{1}{1-\theta}C_E\sqrt{C_P^2+1}\frac{\vert F_\omega \vert_{*,\omega}}{\nu}.
\label{2.7}
\end{equation}  
\end{proposition}

\noindent
\textit{Proof.} 
Let $\mathbf{y}^\epsilon $ be a solution of (\ref{2.ns3}). We put
$\mathbf{v}=\mathbf{y}^\epsilon- E(\mathbf{y}^\epsilon_{\vert\omega})$ in (\ref{2.ns3}).
Since $\mathbf{v}=0$ in $\omega$, we get
$$
\epsilon a_\Omega\left(\mathbf{y}^\epsilon,\mathbf{y}^\epsilon- E(\mathbf{y}^\epsilon_{\vert\omega})\right)
+\epsilon \tilde{c}_\Omega\left(\mathbf{y}^\epsilon,\mathbf{y}^\epsilon,
\mathbf{y}^\epsilon- E(\mathbf{y}^\epsilon_{\vert\omega})\right) =0.
$$
After simplification by $\epsilon $ and using that
$\tilde{c}_\Omega\left(\mathbf{y}^\epsilon,\mathbf{y}^\epsilon,\mathbf{y}^\epsilon\right) =0$, it follows
\begin{eqnarray*}
\| \mathbf{y}^\epsilon \|_{1,\Omega}^2 &=&
a_\Omega\left(\mathbf{y}^\epsilon,E(\mathbf{y}^\epsilon_{\vert\omega})\right)
+\tilde{c}_\Omega\left(\mathbf{y}^\epsilon,\mathbf{y}^\epsilon,E(\mathbf{y}^\epsilon_{\vert\omega})\right)
\\
&\leq & \| \mathbf{y}^\epsilon \|_{1,\Omega} \| E(\mathbf{y}^\epsilon_{\vert\omega}) \|_{1,\Omega}
+\| \tilde{c}_\Omega \|_{1,\Omega} \| \mathbf{y}^\epsilon \|_{1,\Omega}^2
\| E(\mathbf{y}^\epsilon_{\vert\omega}) \|_{1,\Omega}.
\end{eqnarray*}
Then
\begin{eqnarray*}
\| \mathbf{y}^\epsilon \|_{1,\Omega}
&\leq & \| E(\mathbf{y}^\epsilon_{\vert\omega}) \|_{1,\Omega}
+\| \tilde{c}_\Omega \|_{1,\Omega} \| \mathbf{y}^\epsilon \|_{1,\Omega}
\| E(\mathbf{y}^\epsilon_{\vert\omega}) \|_{1,\Omega}\\
&\leq & C_E \| \mathbf{y}^\epsilon \|_{1,\omega}
+\| \tilde{c}_\Omega \|_{1,\Omega} C_E \| \mathbf{y}^\epsilon \|_{1,\omega}\| \mathbf{y}^\epsilon \|_{1,\Omega}
\\
&\leq &  C_E \sqrt{C_P^2+1}  \vert \mathbf{y}^\epsilon \vert_{1,\omega}
+\| \tilde{c}_\Omega \|_{1,\Omega} C_E \sqrt{C_P^2+1}
\vert \mathbf{y}^\epsilon \vert_{1,\omega}
\| \mathbf{y}^\epsilon \|_{1,\Omega}
\end{eqnarray*}
and using (\ref{2.2}) and (\ref{2.6}), we get
\begin{eqnarray*}
\| \mathbf{y}^\epsilon \|_{1,\Omega} &\leq&
C_E \sqrt{C_P^2+1} \frac{\vert F_\omega \vert_{*,\omega}}{\nu}
+\| \tilde{c}_\Omega \|_{1,\Omega} C_E \sqrt{C_P^2+1}\frac{\vert F_\omega \vert_{*,\omega}}{\nu}
\| \mathbf{y}^\epsilon \|_{1,\Omega}\\
&\leq& C_E \sqrt{C_P^2+1} \frac{\vert F_\omega \vert_{*,\omega}}{\nu}
+\theta \| \mathbf{y}^\epsilon \|_{1,\Omega}
\end{eqnarray*}
which gives (\ref{2.7}).
 \quad $\Box$

\bigskip

\begin{proposition}\label{prop:2.3}
i) Let $0 <\theta <1$ be fixed. Under the hypothesis (\ref{2.6}), we have
\begin{equation}
\label{2.9}
\lim_{\epsilon\rightarrow 0} \mathbf{y}^\epsilon _{\vert\omega} = \widehat{\mathbf{y}}_\omega
\end{equation}
in $H^1(\omega)^d$ weakly on a subsequence, where $\widehat{\mathbf{y}}_\omega\in V_1$
is a solution of (\ref{1.ns3_v1}).

ii) Let $0 <\theta <1$ be fixed. Under the hypothesis 
\begin{equation}
\label{2.H3}
\vert \tilde{c}_\omega \vert_{1,\omega} \vert F_\omega \vert_{*,\omega}  \leq \theta\nu^2,
\end{equation}
together with (\ref{2.6}), then
\begin{equation}
\label{2.error}
\vert \mathbf{y}^\epsilon -\widehat{\mathbf{y}}_\omega \vert_{1,\omega} \leq \epsilon C,
\end{equation} 
where $C$ is some constant independent of $\epsilon$.
\end{proposition}

\noindent
\textit{Proof.} 
i) From (\ref{2.2}) and (\ref{2.7}), there exists $\widehat{\mathbf{y}}\in V$
such that $\lim_{\epsilon\rightarrow 0} \mathbf{y}^\epsilon= \widehat{\mathbf{y}}$
in $V$ weakly on a subsequence.
From (\ref{2.ns3}), we have
\begin{eqnarray*}
a_\omega(\mathbf{y}^\epsilon,\mathbf{v}) + \epsilon a_\Omega(\mathbf{y}^\epsilon,\mathbf{v})
+
\tilde{c}_\omega(\mathbf{y}^\epsilon,\mathbf{y}^\epsilon,\mathbf{v})
+\epsilon\tilde{c}_\Omega(\mathbf{y}^\epsilon,\mathbf{y}^\epsilon,\mathbf{v})
&=&\langle F_\omega, \mathbf{v}\rangle_{*,\omega},\quad\forall \mathbf{v}\in V.
\end{eqnarray*}
From (\ref{2.2}) and (\ref{2.7}) and the continuity of $ a_\Omega$, 
$\tilde{c}_\Omega$, $ a_\omega$ and
$\tilde{c}_\omega$ we  have that
$$
\lim_{\epsilon\rightarrow 0} a_\omega(\mathbf{y}^\epsilon,\mathbf{v})
=a_\omega(\widehat{\mathbf{y}},\mathbf{v}),\quad
\lim_{\epsilon\rightarrow 0} \tilde{c}_\omega(\mathbf{y}^\epsilon,\mathbf{y}^\epsilon,\mathbf{v})
=\tilde{c}_\omega(\widehat{\mathbf{y}},\widehat{\mathbf{y}},\mathbf{v})
$$
on a subsequence.
Passing to the limit, we get
\begin{equation}
\label{2.ns3_hat}
a_\omega(\widehat{\mathbf{y}},\mathbf{v})
+\tilde{c}_\omega(\widehat{\mathbf{y}},\widehat{\mathbf{y}},\mathbf{v})
=\langle F_\omega, \mathbf{v}\rangle_{*,\omega},\quad\forall \mathbf{v}\in V.
\end{equation}
In view of Rem. \ref{rem:2.1bis}, we get that 
$\widehat{\mathbf{y}}_\omega=\widehat{\mathbf{y}}_{\vert\omega}\in V_1$ is a solution of (\ref{1.ns3_v1}).

ii)  Assuming (\ref{2.H3}), then (\ref{1.H1}) holds and (\ref{1.ns3_v1}) has a unique solution
$\widehat{\mathbf{y}}_\omega\in V_1$.
Let $\mathbf{v}$ be in $V$. Then $\mathbf{v}_{\vert\omega}\in V_1$ and subtracting (\ref{1.ns3_v1})
from (\ref{2.ns3}), we obtain
\begin{eqnarray*}
a_\omega(\mathbf{y}^\epsilon -\widehat{\mathbf{y}}_\omega,\mathbf{v} ) & = &
-\epsilon a_\Omega(\mathbf{y}^\epsilon ,\mathbf{v} )
-\epsilon \tilde{c}_\Omega (\mathbf{y}^\epsilon ,\mathbf{y}^\epsilon , \mathbf{v} )
\\
&&  -\tilde{c}_\omega (\mathbf{y}^\epsilon ,\mathbf{y}^\epsilon , \mathbf{v} )
+\tilde{c}_\omega (\widehat{\mathbf{y}}_\omega ,\widehat{\mathbf{y}}_\omega , \mathbf{v} ).
\end{eqnarray*}
We have
$$
-\tilde{c}_\omega (\mathbf{y}^\epsilon ,\mathbf{y}^\epsilon , \mathbf{v} )
+\tilde{c}_\omega (\widehat{\mathbf{y}}_\omega ,\widehat{\mathbf{y}}_\omega , \mathbf{v} )
=-\tilde{c}_\omega (\mathbf{y}^\epsilon -\widehat{\mathbf{y}}_\omega,\mathbf{y}^\epsilon , \mathbf{v} )
-\tilde{c}_\omega (\widehat{\mathbf{y}}_\omega ,\mathbf{y}^\epsilon -\widehat{\mathbf{y}}_\omega , \mathbf{v} ).
$$
From Rem. \ref{rem:2.1bis}, we can extend $(\mathbf{y}^\epsilon -\widehat{\mathbf{y}}_\omega)_{\vert\omega}$ and 
we put $\mathbf{v}=E\left( \mathbf{y}^\epsilon -\widehat{\mathbf{y}}_\omega\right)\in V$. Using
that $E\left( \mathbf{y}^\epsilon -\widehat{\mathbf{y}}_\omega\right)=\mathbf{y}^\epsilon -\widehat{\mathbf{y}}_\omega$
in $\omega$, we get
\begin{eqnarray*}
  \nu \vert\mathbf{y}^\epsilon -\widehat{\mathbf{y}}_\omega\vert_{1,\omega}^2 & = &
-\epsilon a_\Omega(\mathbf{y}^\epsilon , E\left( \mathbf{y}^\epsilon -\widehat{\mathbf{y}}_\omega\right))
-\epsilon \tilde{c}_\Omega (\mathbf{y}^\epsilon ,\mathbf{y}^\epsilon ,
E\left( \mathbf{y}^\epsilon -\widehat{\mathbf{y}}_\omega\right) )  \\
&&
-\tilde{c}_\omega (\mathbf{y}^\epsilon -\widehat{\mathbf{y}}_\omega,\mathbf{y}^\epsilon ,
\mathbf{y}^\epsilon -\widehat{\mathbf{y}}_\omega )
-\tilde{c}_\omega (\widehat{\mathbf{y}}_\omega ,\mathbf{y}^\epsilon -\widehat{\mathbf{y}}_\omega ,
\mathbf{y}^\epsilon -\widehat{\mathbf{y}}_\omega ).
\end{eqnarray*}
Since the second and the third arguments are the same, the last term is zero.

Using (\ref{2.5}) and $\|\cdot\|_{1,\omega}\leq \sqrt{C_P^2+1}\vert\cdot\vert_{1,\omega}$, it follows
\begin{eqnarray*}
\nu \vert\mathbf{y}^\epsilon -\widehat{\mathbf{y}}_\omega\vert_{1,\omega}^2 & \leq &
\epsilon \|\mathbf{y}^\epsilon\|_{1,\Omega}  \|E\left( \mathbf{y}^\epsilon -\widehat{\mathbf{y}}_\omega\right)\|_{1,\Omega}
+\epsilon  \| \tilde{c}_\Omega \|_{1,\Omega} \| \mathbf{y}^\epsilon \|_{1,\Omega}^2
\|E\left( \mathbf{y}^\epsilon -\widehat{\mathbf{y}}_\omega\right)\|_{1,\Omega}
\\
&& + \vert \tilde{c}_\omega \vert_{1,\omega} \vert \mathbf{y}^\epsilon -\widehat{\mathbf{y}}_\omega \vert_{1,\omega}^2
\vert \mathbf{y}^\epsilon \vert_{1,\omega}
\\
& \leq &
\epsilon \|\mathbf{y}^\epsilon\|_{1,\Omega}  C_E\|\mathbf{y}^\epsilon -\widehat{\mathbf{y}}_\omega\|_{1,\omega}
+\epsilon  \| \tilde{c}_\Omega \|_{1,\Omega} \| \mathbf{y}^\epsilon \|_{1,\Omega}^2
C_E\|\mathbf{y}^\epsilon -\widehat{\mathbf{y}}_\omega\|_{1,\omega}
\\
&& + \vert \tilde{c}_\omega \vert_{1,\omega} \vert \mathbf{y}^\epsilon -\widehat{\mathbf{y}}_\omega \vert_{1,\omega}^2
\vert \mathbf{y}^\epsilon \vert_{1,\omega}
\\
& \leq &
\epsilon \|\mathbf{y}^\epsilon\|_{1,\Omega}
C_E\sqrt{C_P^2+1} \vert\mathbf{y}^\epsilon -\widehat{\mathbf{y}}_\omega \vert_{1,\omega}
\\
&&
+\epsilon  \| \tilde{c}_\Omega \|_{1,\Omega} \| \mathbf{y}^\epsilon \|_{1,\Omega}^2
C_E\sqrt{C_P^2+1}\vert\mathbf{y}^\epsilon -\widehat{\mathbf{y}}_\omega\vert_{1,\omega}
\\
&&
+ \vert \tilde{c}_\omega \vert_{1,\omega} \vert \mathbf{y}^\epsilon -\widehat{\mathbf{y}}_\omega \vert_{1,\omega}^2
\vert \mathbf{y}^\epsilon \vert_{1,\omega}.
\end{eqnarray*}
Assuming $\vert\mathbf{y}^\epsilon -\widehat{\mathbf{y}}_\omega\vert_{1,\omega}\neq 0$,
after simplification, we obtain
\begin{eqnarray*}
\nu \vert\mathbf{y}^\epsilon -\widehat{\mathbf{y}}_\omega\vert_{1,\omega}
&\leq&
\epsilon \|\mathbf{y}^\epsilon\|_{1,\Omega}
C_E\sqrt{C_P^2+1}
+\epsilon  \| \tilde{c}_\Omega \|_{1,\Omega} \| \mathbf{y}^\epsilon \|_{1,\Omega}^2
C_E\sqrt{C_P^2+1}\\
&&+\vert \tilde{c}_\omega \vert_{1,\omega} \vert \mathbf{y}^\epsilon -\widehat{\mathbf{y}}_\omega \vert_{1,\omega}
\vert \mathbf{y}^\epsilon \vert_{1,\omega}
\end{eqnarray*}
and using (\ref{2.7}) and (\ref{2.2}), we get
$$
\nu \vert\mathbf{y}^\epsilon -\widehat{\mathbf{y}}_\omega\vert_{1,\omega}\leq
\epsilon\, \widehat{C}
+\vert \tilde{c}_\omega \vert_{1,\omega}\frac{\vert F_\omega \vert_{*,\omega}}{\nu}
\vert \mathbf{y}^\epsilon -\widehat{\mathbf{y}}_\omega \vert_{1,\omega}.
$$
Under the hypotheses (\ref{2.H3}), we get the conclusion of ii).
 \quad $\Box$

\begin{remark}\label{rem:2.1}
We underline that the equality (\ref{2.ns3_hat}) is related to (\ref{1.ns3_v1}) that is
not equivalent to (\ref{1.ns3}). See Rem. \ref{rem:2.1bis}.
In fact $\{\mathbf{v}_{\vert\omega};\ \mathbf{v}\in V \}=V_1$ is a strict subspace of $V_\omega$.
They coincide when, for example,
$\Gamma_D=\partial D$, $\Gamma_N=\emptyset$ and $\Omega$ connected.
In this case, for $\mathbf{v}_{\omega}\in V_\omega$,
$$
0=\int_{\omega} \nabla \cdot \mathbf{v}_{\omega}\, d\mathbf{x}
=\int_{\partial D} \mathbf{v}_{\omega}\cdot \mathbf{n}\, ds 
+\int_{\partial \Omega} \mathbf{v}_{\omega}\cdot \mathbf{n}\, ds
=\int_{\partial \Omega} \mathbf{v}_{\omega}\cdot \mathbf{n}\, ds,
$$
then $\mathbf{v}_{\omega}\in V_1$.
When homogeneous Dirichlet conditions 
are imposed on $\partial \omega$, then all the above arguments remain valid for $\Omega$
not necessarily connected.
\end{remark}

\begin{proposition}\label{prop:2.3iii}
Let $0 <\theta <1$ be fixed.
Under the hypotheses 
\begin{equation}
\label{2.H4}
\| \tilde{c}_\Omega \|_{1,\Omega}
C_E\sqrt{C_P^2+1}\frac{\vert F_\omega \vert_{*,\omega}}{\nu} \leq  \theta(1-\theta)
\end{equation}
\begin{equation}
\label{2.H5}  
0 < \epsilon < \epsilon_0=\frac{1-\theta}{\theta \nu C_1^2},\quad
C_1=\frac{(1+2\theta)\,C_E\sqrt{C_P^2+1}}{(1-\theta)\nu }
\end{equation}
and (\ref{2.H3}), if $\mathbf{y}^\epsilon,\mathbf{z}^\epsilon\in V$
are two solutions of (\ref{2.ns3}),
then $\mathbf{y}^\epsilon=\mathbf{z}^\epsilon$ in $D$.
\end{proposition}

\noindent
\textit{Proof.} 
Subtracting (\ref{2.ns3}) written for $\mathbf{z}^\epsilon$ from (\ref{2.ns3}), we get
\begin{eqnarray*}
a_\omega(\mathbf{y}^\epsilon -\mathbf{z}^\epsilon,\mathbf{v} ) & = &
-\epsilon a_\Omega(\mathbf{y}^\epsilon -\mathbf{z}^\epsilon,\mathbf{v} )\\
&& 
-\epsilon \tilde{c}_\Omega (\mathbf{y}^\epsilon ,\mathbf{y}^\epsilon , \mathbf{v} )
+\epsilon \tilde{c}_\Omega (\mathbf{z}^\epsilon ,\mathbf{z}^\epsilon , \mathbf{v} )
\\
&&  -\tilde{c}_\omega (\mathbf{y}^\epsilon ,\mathbf{y}^\epsilon , \mathbf{v} )
+\tilde{c}_\omega (\mathbf{z}^\epsilon ,\mathbf{z}^\epsilon , \mathbf{v} )
\\
&=&
-\epsilon a_\Omega(\mathbf{y}^\epsilon -\mathbf{z}^\epsilon,\mathbf{v} )\\
&& 
-\epsilon \tilde{c}_\Omega (\mathbf{y}^\epsilon -\mathbf{z}^\epsilon,\mathbf{y}^\epsilon , \mathbf{v} )
-\epsilon \tilde{c}_\Omega (\mathbf{z}^\epsilon ,\mathbf{y}^\epsilon -\mathbf{z}^\epsilon , \mathbf{v} )
\\
&&  -\tilde{c}_\omega (\mathbf{y}^\epsilon  -\mathbf{z}^\epsilon,\mathbf{y}^\epsilon , \mathbf{v} )
-\tilde{c}_\omega (\mathbf{z}^\epsilon ,\mathbf{y}^\epsilon - \mathbf{z}^\epsilon , \mathbf{v} ).
\end{eqnarray*}

We can extend $(\mathbf{y}^\epsilon -\mathbf{z}^\epsilon)_{\vert\omega}$ and 
we put $\mathbf{v}=E\left( \mathbf{y}^\epsilon -\mathbf{z}^\epsilon\right)\in V$. We obtain
\begin{eqnarray*}
  \nu \vert\mathbf{y}^\epsilon -\mathbf{z}^\epsilon\vert_{1,\omega}^2 & = &
-\epsilon a_\Omega(\mathbf{y}^\epsilon -\mathbf{z}^\epsilon, E\left( \mathbf{y}^\epsilon -\mathbf{z}^\epsilon\right))
\\
&&
-\epsilon \tilde{c}_\Omega (\mathbf{y}^\epsilon -\mathbf{z}^\epsilon,\mathbf{y}^\epsilon ,
E\left( \mathbf{y}^\epsilon -\mathbf{z}^\epsilon\right) )  
-\epsilon \tilde{c}_\Omega (\mathbf{z}^\epsilon ,\mathbf{y}^\epsilon -\mathbf{z}^\epsilon , 
E\left( \mathbf{y}^\epsilon -\mathbf{z}^\epsilon\right) )
\\
&&  -\tilde{c}_\omega (\mathbf{y}^\epsilon  -\mathbf{z}^\epsilon,\mathbf{y}^\epsilon , 
E\left( \mathbf{y}^\epsilon -\mathbf{z}^\epsilon\right) )
-\tilde{c}_\omega (\mathbf{z}^\epsilon ,\mathbf{y}^\epsilon - \mathbf{z}^\epsilon , 
E\left( \mathbf{y}^\epsilon -\mathbf{z}^\epsilon\right) ).
\end{eqnarray*}
But $E\left( \mathbf{y}^\epsilon -\mathbf{z}^\epsilon\right)=\mathbf{y}^\epsilon -\mathbf{z}^\epsilon$ in $\omega$,
then the last term is zero.

It follows, using (\ref{2.5}) and $\|\cdot\|_{1,\omega}\leq \sqrt{C_P^2+1}\vert\cdot\vert_{1,\omega}$
\begin{eqnarray*}
&&\nu \vert\mathbf{y}^\epsilon -\mathbf{z}^\epsilon\vert_{1,\omega}^2  \leq 
\epsilon \|\mathbf{y}^\epsilon-\mathbf{z}^\epsilon\|_{1,\Omega}  
\|E\left( \mathbf{y}^\epsilon -\mathbf{z}^\epsilon\right)\|_{1,\Omega}
\\
&& 
+ \epsilon \|\tilde{c}_\Omega\|_{1,\Omega} 
\|\mathbf{y}^\epsilon-\mathbf{z}^\epsilon\|_{1,\Omega}
\|E\left( \mathbf{y}^\epsilon -\mathbf{z}^\epsilon\right)\|_{1,\Omega}
\left( \|\mathbf{y}^\epsilon\|_{1,\Omega}  + \|\mathbf{z}^\epsilon\|_{1,\Omega} \right)
\\
&&
+\vert \tilde{c}_\omega \vert_{1,\omega} 
\vert \mathbf{y}^\epsilon -\mathbf{z}^\epsilon \vert_{1,\omega} ^2
\vert \mathbf{y}^\epsilon \vert_{1,\omega}
\\
&& \leq 
\epsilon \|\mathbf{y}^\epsilon-\mathbf{z}^\epsilon\|_{1,\Omega}  
\,C_E\sqrt{C_P^2+1}\vert \mathbf{y}^\epsilon -\mathbf{z}^\epsilon \vert_{1,\omega}
\\
&& 
+ \epsilon \|\tilde{c}_\Omega\|_{1,\Omega} 
\|\mathbf{y}^\epsilon-\mathbf{z}^\epsilon\|_{1,\Omega}
\left( \|\mathbf{y}^\epsilon\|_{1,\Omega}  + \|\mathbf{z}^\epsilon\|_{1,\Omega} \right)
\,C_E\sqrt{C_P^2+1}\vert \mathbf{y}^\epsilon -\mathbf{z}^\epsilon \vert_{1,\omega}
\\
&&
+\vert \tilde{c}_\omega \vert_{1,\omega} 
\vert \mathbf{y}^\epsilon -\mathbf{z}^\epsilon \vert_{1,\omega} ^2
\vert \mathbf{y}^\epsilon \vert_{1,\omega}.
\end{eqnarray*}

If $\vert \mathbf{y}^\epsilon -\mathbf{z}^\epsilon \vert_{1,\omega}\neq 0$, we can simplify and we obtain
\begin{eqnarray*}
\nu \vert\mathbf{y}^\epsilon -\mathbf{z}^\epsilon\vert_{1,\omega} & \leq &
\epsilon \|\mathbf{y}^\epsilon-\mathbf{z}^\epsilon\|_{1,\Omega}  
\,C_E\sqrt{C_P^2+1}
\\
&& 
+ \epsilon \|\tilde{c}_\Omega\|_{1,\Omega} 
\|\mathbf{y}^\epsilon-\mathbf{z}^\epsilon\|_{1,\Omega}
\left( \|\mathbf{y}^\epsilon\|_{1,\Omega}  + \|\mathbf{z}^\epsilon\|_{1,\Omega} \right)
\,C_E\sqrt{C_P^2+1}
\\
&&
+\vert \tilde{c}_\omega \vert_{1,\omega} 
\vert \mathbf{y}^\epsilon -\mathbf{z}^\epsilon \vert_{1,\omega} 
\vert \mathbf{y}^\epsilon \vert_{1,\omega}.
\end{eqnarray*}
The inequality (\ref{2.H4}) implies (\ref{2.6}), then from Prop. \ref{prop:2.2}
we have (\ref{2.7}). Using (\ref{2.7}), (\ref{2.2}), it follows
\begin{eqnarray*}
\nu \vert\mathbf{y}^\epsilon -\mathbf{z}^\epsilon\vert_{1,\omega} & \leq &
\epsilon \|\mathbf{y}^\epsilon-\mathbf{z}^\epsilon\|_{1,\Omega}  
\,C_E\sqrt{C_P^2+1}
\\
&& 
+ \epsilon \|\tilde{c}_\Omega\|_{1,\Omega} 
\|\mathbf{y}^\epsilon-\mathbf{z}^\epsilon\|_{1,\Omega}
\,\frac{2}{1-\theta}C_E\sqrt{C_P^2+1}\frac{\vert F_\omega \vert_{*,\omega}}{\nu}
\,C_E\sqrt{C_P^2+1}
\\
&&
+\vert \tilde{c}_\omega \vert_{1,\omega} 
\vert \mathbf{y}^\epsilon -\mathbf{z}^\epsilon \vert_{1,\omega} 
\frac{\vert F_\omega \vert_{*,\omega}}{\nu}.
\end{eqnarray*}

Above, in the right-hand side, we use (\ref{2.H4}) in the second term and (\ref{2.H3})
in the third term, we get
\begin{eqnarray*}
\nu \vert\mathbf{y}^\epsilon -\mathbf{z}^\epsilon\vert_{1,\omega} & \leq &
\epsilon \|\mathbf{y}^\epsilon-\mathbf{z}^\epsilon\|_{1,\Omega}  
\,C_E\sqrt{C_P^2+1}
\\
&&
+ \epsilon 
\|\mathbf{y}^\epsilon-\mathbf{z}^\epsilon\|_{1,\Omega}
\, 2\theta C_E\sqrt{C_P^2+1}
\\
&&
+\theta\nu
\vert \mathbf{y}^\epsilon -\mathbf{z}^\epsilon \vert_{1,\omega} 
\end{eqnarray*}
or, by passing the last term to the left-hand side,
$$
(1-\theta)\nu \vert\mathbf{y}^\epsilon -\mathbf{z}^\epsilon\vert_{1,\omega}
\leq
\epsilon 
\|\mathbf{y}^\epsilon-\mathbf{z}^\epsilon\|_{1,\Omega}
\, (1+2\theta) C_E\sqrt{C_P^2+1}
$$
and finally
\begin{equation}
\label{2.error2}
\vert \mathbf{y}^\epsilon - \mathbf{z}^\epsilon\vert_{1,\omega} \leq \epsilon C_1
\| \mathbf{y}^\epsilon - \mathbf{z}^\epsilon\|_{1,\Omega},
\quad C_1=\frac{(1+2\theta)\,C_E\sqrt{C_P^2+1}}{(1-\theta)\nu } .
\end{equation}

Putting $\mathbf{v}=\mathbf{y}^\epsilon -\mathbf{z}^\epsilon$ in the first formula of
the proof of this Proposition, we get
\begin{eqnarray*}
&&  a_\omega(\mathbf{y}^\epsilon -\mathbf{z}^\epsilon,\mathbf{y}^\epsilon -\mathbf{z}^\epsilon )
+\epsilon a_\Omega(\mathbf{y}^\epsilon -\mathbf{z}^\epsilon,\mathbf{y}^\epsilon -\mathbf{z}^\epsilon )
\\
&=&
-\epsilon \tilde{c}_\Omega (\mathbf{y}^\epsilon -\mathbf{z}^\epsilon,\mathbf{y}^\epsilon ,
\mathbf{y}^\epsilon -\mathbf{z}^\epsilon )
-\epsilon \tilde{c}_\Omega (\mathbf{z}^\epsilon ,\mathbf{y}^\epsilon -\mathbf{z}^\epsilon ,
\mathbf{y}^\epsilon -\mathbf{z}^\epsilon )
\\
&&  -\tilde{c}_\omega (\mathbf{y}^\epsilon  -\mathbf{z}^\epsilon,\mathbf{y}^\epsilon ,
\mathbf{y}^\epsilon -\mathbf{z}^\epsilon )
-\tilde{c}_\omega (\mathbf{z}^\epsilon ,\mathbf{y}^\epsilon - \mathbf{z}^\epsilon ,
\mathbf{y}^\epsilon -\mathbf{z}^\epsilon )
\\
&=& -\epsilon \tilde{c}_\Omega (\mathbf{y}^\epsilon -\mathbf{z}^\epsilon,\mathbf{y}^\epsilon ,
\mathbf{y}^\epsilon -\mathbf{z}^\epsilon )
-\tilde{c}_\omega (\mathbf{y}^\epsilon  -\mathbf{z}^\epsilon,\mathbf{y}^\epsilon ,
\mathbf{y}^\epsilon -\mathbf{z}^\epsilon ).
\end{eqnarray*}
We obtain
\begin{eqnarray*}
  \epsilon \|\mathbf{y}^\epsilon -\mathbf{z}^\epsilon\|_{1,\Omega}^2 & \leq &
  \epsilon  \| \tilde{c}_\Omega\|_{1,\Omega} \| \mathbf{y}^\epsilon -\mathbf{z}^\epsilon\|_{1,\Omega}^2
  \|\mathbf{y}^\epsilon \|_{1,\Omega}
  +\vert \tilde{c}_\omega\vert_{1,\omega} \vert \mathbf{y}^\epsilon -\mathbf{z}^\epsilon\vert_{1,\omega}^2
  \vert\mathbf{y}^\epsilon \vert_{1,\omega} .
\end{eqnarray*}
We have seen that the inequality (\ref{2.H4}) implies (\ref{2.6}), then
from Prop. \ref{prop:2.2} we have (\ref{2.7}).
It follows 
\begin{eqnarray*}
  \epsilon \|\mathbf{y}^\epsilon -\mathbf{z}^\epsilon\|_{1,\Omega}^2 & \leq &
  \epsilon  \| \tilde{c}_\Omega\|_{1,\Omega} \| \mathbf{y}^\epsilon -\mathbf{z}^\epsilon\|_{1,\Omega}^2
  \frac{1}{1-\theta}C_E\sqrt{C_P^1+1}\frac{\vert F_\omega \vert_{*,\omega}}{\nu}
  \\
  &&
  +\vert \tilde{c}_\omega\vert_{1,\omega} \vert \mathbf{y}^\epsilon -\mathbf{z}^\epsilon\vert_{1,\omega}^2
  \vert\mathbf{y}^\epsilon \vert_{1,\omega}
\\
&\leq &  \epsilon \theta \| \mathbf{y}^\epsilon -\mathbf{z}^\epsilon\|_{1,\Omega}^2
+\vert \tilde{c}_\omega\vert_{1,\omega} \vert \mathbf{y}^\epsilon -\mathbf{z}^\epsilon\vert_{1,\omega}^2
  \vert\mathbf{y}^\epsilon \vert_{1,\omega}.
\end{eqnarray*}
For the last inequality, we used (\ref{2.H4}).
We obtain from (\ref{2.2}) and (\ref{2.H3})
\begin{eqnarray*}
&&(1-\theta)  \epsilon \|\mathbf{y}^\epsilon -\mathbf{z}^\epsilon\|_{1,\Omega}^2  \leq 
\vert \tilde{c}_\omega\vert_{1,\omega} \vert \mathbf{y}^\epsilon -\mathbf{z}^\epsilon\vert_{1,\omega}^2
\vert\mathbf{y}^\epsilon \vert_{1,\omega}
\\
& \leq &
\vert \tilde{c}_\omega\vert_{1,\omega} \frac{\vert F_\omega \vert_{*,\omega}}{\nu}
\vert \mathbf{y}^\epsilon -\mathbf{z}^\epsilon\vert_{1,\omega}^2
\leq \theta \nu \vert \mathbf{y}^\epsilon -\mathbf{z}^\epsilon\vert_{1,\omega}^2
\end{eqnarray*}
or
\begin{equation}
  \label{2.error3}
  \|\mathbf{y}^\epsilon -\mathbf{z}^\epsilon\|_{1,\Omega}
  \leq
  \sqrt{\frac{\theta \nu}{(1-\theta)  \epsilon}} \, \vert \mathbf{y}^\epsilon -\mathbf{z}^\epsilon\vert_{1,\omega}.
\end{equation}
From (\ref{2.error2}) and (\ref{2.error3}), we obtain
\begin{eqnarray*}
  \vert \mathbf{y}^\epsilon -\mathbf{z}^\epsilon\vert_{1,\omega}\leq \epsilon C_1
  \sqrt{\frac{\theta \nu}{(1-\theta)  \epsilon}} \, \vert \mathbf{y}^\epsilon -\mathbf{z}^\epsilon\vert_{1,\omega}
  \leq \sqrt{\epsilon} C_1 \sqrt{\frac{\theta \nu}{(1-\theta)} } \,
    \vert \mathbf{y}^\epsilon -\mathbf{z}^\epsilon\vert_{1,\omega}
\end{eqnarray*}

By (\ref{2.H5}), the coefficient of the last term is less than $1$ and we get
 $\vert \mathbf{y}^\epsilon -\mathbf{z}^\epsilon\vert_{1,\omega} = 0$.

\noindent
Again by (\ref{2.error3}), it yields $\| \mathbf{y}^\epsilon -\mathbf{z}^\epsilon\|_{1,\Omega}= 0$ and finally $\mathbf{y}^\epsilon =\mathbf{z}^\epsilon$
in $D$.
 \quad $\Box$

\begin{remark}\label{rem:2.3}
Assuming that $D$, $\Omega$ and $\nu$ are fixed, in Prop. \ref{prop:2.3iii}, we can
choose $F_\omega$ ``small'', satisfying (\ref{2.H3}) and (\ref{2.H4}), independently
of the choice of $\epsilon$ given by (\ref{2.H5}). The uniqueness result from
Prop. \ref{prop:2.3iii}  has
a non standard character due to the very weak coercivity properties of the
bilinear form (\ref{2.ns3}), defined in $D$.
For steady Navier-Stokes equations in $\omega$, if the body forces and the Neumann boundary condition are
given, from (\ref{1.H1}) we get the uniqueness of the solution, under the condition that the viscosity
is larger than a threshold.
For the penalized problem, the penalization parameter $\epsilon$ can be interpreted as a viscosity in $\Omega$.
The Prop. \ref{prop:2.3iii} proves the uniqueness of the solution in $D$ even for $\epsilon\rightarrow 0$.
\end{remark}

\begin{proposition}\label{prop:2.4}
Let $0 <\theta <1$ be fixed.  
For each $\mathbf{y}^\epsilon\in V$ solution of (\ref{2.ns3}), there exists a unique $p^\epsilon\in Q$
such that
\begin{equation}
\label{2.ns1}
a^\epsilon(\mathbf{y}^\epsilon,\mathbf{w}) + b^\epsilon(\mathbf{w},p^\epsilon)
+\tilde{c}^\epsilon(\mathbf{y}^\epsilon,\mathbf{y}^\epsilon,\mathbf{w})
=\langle F_\omega, \mathbf{w}\rangle_{*,\omega},\quad \forall \mathbf{w}\in W.
\end{equation}  
Under the hypotheses (\ref{2.6}) and $0 < \epsilon \leq 1$, we have
\begin{eqnarray}
\| p^\epsilon \|_{0,\omega} &\leq & C(\vert F_\omega \vert_{*,\omega} + \vert F_\omega \vert_{*,\omega}^2), \label{2.10}\\
\epsilon\| p^\epsilon \|_{0,\Omega}  &\leq & C (\vert F_\omega \vert_{*,\omega} + \vert F_\omega \vert_{*,\omega}^2).
\label{2.11}
\end{eqnarray}  
\end{proposition}

\noindent
\textit{Proof.} 
Let $b:W\times Q\rightarrow \mathbb{R}$, be defined by
\begin{equation}\label{b_D}
b(\mathbf{w},q)=-\int_D (\nabla \cdot \mathbf{w}) q\, d\mathbf{x}.
\end{equation}
Since the divergence operator
$
\mathbf{w} \in W \rightarrow \nabla\cdot\mathbf{w} \in Q
$
is onto, see \cite{MR_2007}, Lemma 5.1,
it follows, see \cite{Girault1986}, Chapter I, Lemma 4.1, p. 58, 
that the following \textit{inf-sup} condition holds
$$
\exists \beta >0,\quad
\beta \leq \inf_{q\in Q \setminus \{0\} }\sup_{\mathbf{w}\in W \setminus \{0\} }
\frac{b(\mathbf{w},q)}{\| q\|_{0,D} \|\mathbf{w} \|_{1,D} }
$$
and for each $\mathbf{y}^\epsilon\in V$ solution of (\ref{2.ns3}), there exists a unique $P^\epsilon\in Q$
such that
$$
b^\epsilon(\mathbf{w},P^\epsilon)
=\langle F_\omega, \mathbf{w}\rangle_{*,\omega}
-a^\epsilon(\mathbf{y}^\epsilon,\mathbf{w})
-\tilde{c}^\epsilon(\mathbf{y}^\epsilon,\mathbf{y}^\epsilon,\mathbf{w})
,\quad \forall \mathbf{w}\in W.
$$
Using (\ref{2.2}), (\ref{2.6}), (\ref{2.7}), the \textit{inf-sup}
condition and $0 < \epsilon \leq 1$, we get
$$
\| P^\epsilon \|_{0,D} \leq  C(\vert F_\omega \vert_{*,\omega} + \vert F_\omega \vert_{*,\omega}^2).
$$

If we set $p^\epsilon\in Q$ by $p^\epsilon=P^\epsilon$ in $\omega$ and
$p^\epsilon=\frac{1}{\epsilon}P^\epsilon$ in $\Omega$, it follows
$$
b^\epsilon(\mathbf{w},p^\epsilon)=b(\mathbf{w},P^\epsilon), \quad \forall \mathbf{w}\in W
$$
and we get (\ref{2.ns1}). The uniqueness of $p^\epsilon$ is a consequence of the
uniqueness of $P^\epsilon$.

From
$$
\| P^\epsilon \|_{0,D} ^2 = \| p^\epsilon \|_{0,\omega}^2 + \epsilon^2 \| p^\epsilon \|_{0,\Omega}^2
$$
we get (\ref{2.10}) and (\ref{2.11}).
 \quad $\Box$

\bigskip

\begin{proposition}\label{prop:2.6}
Let $0 <\theta <1$ be fixed. Under the hypotheses (\ref{2.H3}) and (\ref{2.6}), then
\begin{equation}
\label{2.errorP}
\| p^\epsilon -\widehat{p}_\omega \|_{0,\omega} \leq \epsilon C,
\end{equation} 
where $\widehat{p}_\omega$ is defined in  (\ref{1.ns1_v1}) and $C$ is some constant independent of $\epsilon$.
\end{proposition}

\noindent
\textit{Proof.} 
We set
$$
W_2=\left\{ \mathbf{w}\in H^1(\omega)^d;\ \mathbf{w}=0\hbox{ on }\Gamma_D\cup\partial\Omega\right\}
$$
and we have $W_2\subset W_1$.

Each $\mathbf{w}^2\in W_2$ can be extended by zero to $\widetilde{\mathbf{w}}^2 \in W$.
Subtracting from (\ref{2.ns1}) with $\mathbf{w}=\widetilde{\mathbf{w}}^2$ the equation
(\ref{1.ns1_v1}) with $\mathbf{w}^1=\mathbf{w}^2$, we get
\begin{eqnarray*}
  0 &=&
  a_\omega(\mathbf{y}^\epsilon -\widehat{\mathbf{y}}_\omega,\mathbf{w}^2 ) 
+b_\omega(\mathbf{w}^2,p^\epsilon -\widehat{p}_\omega )
 +\tilde{c}_\omega (\mathbf{y}^\epsilon ,\mathbf{y}^\epsilon , \mathbf{v} )
 -\tilde{c}_\omega (\widehat{\mathbf{y}}_\omega ,\widehat{\mathbf{y}}_\omega , \mathbf{v} )\\
 &=&
 a_\omega(\mathbf{y}^\epsilon -\widehat{\mathbf{y}}_\omega,\mathbf{w}^2 ) 
 +b_\omega(\mathbf{w}^2,p^\epsilon -\widehat{p}_\omega )\\
 &&
 +\tilde{c}_\omega (\mathbf{y}^\epsilon -\widehat{\mathbf{y}}_\omega,\mathbf{y}^\epsilon , \mathbf{w}^2 )
 +\tilde{c}_\omega (\widehat{\mathbf{y}}_\omega ,\mathbf{y}^\epsilon -\widehat{\mathbf{y}}_\omega ,
 \mathbf{w}^2 ).
\end{eqnarray*}
Then
\begin{eqnarray*}
  \vert b_\omega(\mathbf{w}^2,p^\epsilon -\widehat{p}_\omega ) \vert
  & \leq & \vert  a_\omega(\mathbf{y}^\epsilon -\widehat{\mathbf{y}}_\omega,\mathbf{w}^2 )  \vert
  + \vert \tilde{c}_\omega (\mathbf{y}^\epsilon -\widehat{\mathbf{y}}_\omega,\mathbf{y}^\epsilon , \mathbf{w}^2 )\vert\\
  &&+ \vert\tilde{c}_\omega (\widehat{\mathbf{y}}_\omega ,\mathbf{y}^\epsilon -\widehat{\mathbf{y}}_\omega ,
  \mathbf{w}^2 )\vert\\
  & \leq & \vert\mathbf{y}^\epsilon -\widehat{\mathbf{y}}_\omega\vert_{1,\omega} \vert\mathbf{w}^2\vert_{1,\omega}
  + 2 \vert \tilde{c}_\omega \vert_{1,\omega} \vert  \vert\mathbf{y}^\epsilon -\widehat{\mathbf{y}}_\omega\vert_{1,\omega}
  \vert\mathbf{y}^\epsilon \vert_{1,\omega}  \vert\mathbf{w}^2\vert_{1,\omega}\\
  & \leq &
  \vert\mathbf{y}^\epsilon -\widehat{\mathbf{y}}_\omega\vert_{1,\omega}
  (1+2 \vert \tilde{c}_\omega \vert_{1,\omega} \vert \vert\mathbf{y}^\epsilon \vert_{1,\omega} ) \|\mathbf{w}^2 \|_{1,\omega}
\end{eqnarray*}
and using (\ref{2.2}) and (\ref{2.error}) we get
$$
  \sup_{\mathbf{w}^2\in W_2\setminus \{ 0 \}}
  \frac{\vert b_\omega(\mathbf{w}^2,p^\epsilon -\widehat{p}_\omega ) \vert}{\|\mathbf{w}^2 \|_{1,\omega}}
  \leq \epsilon C_1.
$$

Using \cite{MR_2007}, Lemma 5.1, the application
$
\mathbf{w}^2 \in W_2 \rightarrow \nabla\cdot\mathbf{w}^2 \in Q_\omega=L^2(\omega)
$
is onto. It follows, see \cite{Girault1986}, Chapter I, Lemma 4.1, p. 58, that
$$
\exists \beta_\omega^\prime >0,\quad
\beta_\omega^\prime \leq \inf_{q_\omega\in Q_\omega \setminus \{0\} }\sup_{\mathbf{w}^2\in W_2 \setminus \{0\} }
\frac{ \vert b_\omega(\mathbf{w}^2,q_\omega)\vert }{\| q_\omega\|_{0,\omega} \|\mathbf{w}^2 \|_{1,\omega}}
$$
then
$$
\beta_\omega^\prime \| p^\epsilon -\widehat{p}_\omega \|_{0,\omega}
\leq \sup_{\mathbf{w}^2\in W_2\setminus \{ 0 \}}
\frac{\vert b_\omega(\mathbf{w}^2,p^\epsilon -\widehat{p}_\omega ) \vert}{\|\mathbf{w}^2 \|_{1,\omega}}
$$
and finally, we get
$ \| p^\epsilon -\widehat{p}_\omega \|_{0,\omega} \leq \epsilon C_1 / \beta_\omega^\prime$.
 \quad $\Box$

\bigskip

To summarize, under the hypotheses (\ref{2.H3}), (\ref{2.H4}) and (\ref{2.H5}),
the problem (\ref{1.ns3_v1}) has a unique solution $\widehat{\mathbf{y}}_\omega\in V_1$,
the problem (\ref{2.ns3}) has a unique solution $\mathbf{y}^\epsilon\in V$ and
(\ref{2.error}) holds.

As in the precedent section, we can introduce the mixed version of (\ref{2.ns3}): find
$\mathbf{y}^\epsilon\in W$ and $p^\epsilon\in Q$ such that
\begin{eqnarray}
a^\epsilon(\mathbf{y}^\epsilon,\mathbf{w}) + b^\epsilon(\mathbf{w},p^\epsilon)
+\tilde{c}^\epsilon(\mathbf{y}^\epsilon,\mathbf{y}^\epsilon,\mathbf{w})
&=&\langle F_\omega, \mathbf{w}\rangle_{*,\omega},\quad \forall \mathbf{w}\in W \label{2.20}\\
b^\epsilon(\mathbf{y}^\epsilon,q)&=&0, \quad \forall q\in Q.\label{2.21}
\end{eqnarray} 
If we replace in (\ref{2.20})-(\ref{2.21}), $b^\epsilon$ by $b$ given by (\ref{b_D}),
we get an equivalent system.
We point out that, 
if $\mathbf{y}^\epsilon\in W$ and $p^\epsilon\in Q$
is a solution of (\ref{2.20})-(\ref{2.21}), then $\mathbf{y}^\epsilon$ is a solution
of (\ref{2.ns3}). Reciprocally, from Prop. \ref{prop:2.4}, if
$\mathbf{y}^\epsilon\in V$ is a solution of (\ref{2.ns3}), there exists a unique
$p^\epsilon\in Q$ such that (\ref{2.20})-(\ref{2.21}) hold.

\subsection{Discretization and Regularization\label{sec:2.2}}

All the previous results hold for $d\in\{2,3\}$. In the following, we fix
 $d=2$.
We assume that $D$ is a polygonal domain in $\mathbb{R}^2$, and let 
$\{\mathcal{T}_h\}_{h>0}$ be a family of regular triangulations of $D$ with $h$ the mesh size.
Based on the free divergence finite element introduced in \cite{Temam1984}, Chap. 1, Sec. 4.4, 
one can construct  a finite dimensional subspace of $V$ giving an internal approximation $V_h$
of $V$. In \cite{Temam1984}, only the homogeneous Dirichlet boundary condition for
$\partial D=\Gamma_D$ is discussed.
In \cite{Scott1985} and in \cite{Brenner2008}, Chap. 12, using continuous $\mathbb{P}_k$
for $W$ and discontinuous $\mathbb{P}_{k-1}^{dc}$ for $Q$, $k\geq 4$, 
other cases when $V_h\subset V$ are discussed.

Alternatively, based on mixed finite elements, see \cite{Girault1986} and \cite{Brezzi2013}, we can
construct $W_h\subset W$ and $Q_h\subset Q$, finite dimensional subspaces. We define
\begin{equation}\label{V_h}
V_h=\left\{ \mathbf{v}_h\in W_h;\ b(\mathbf{v}_h,q_h)=0,\ \forall q_h\in Q_h\right\}.
\end{equation}
We assume that for all $\mathbf{w}\in W$, $q\in Q$, $\mathbf{v}\in V$ there exist
$\{\mathbf{w}_h\}_{h>0}$ in $W_h$, $\{q_h\}_{h>0}$ in $Q_h$ and $\{\mathbf{v}_h\}_{h>0}$ in $V_h$, such that
\begin{equation}\label{density}
\lim_{h\rightarrow 0} \|\mathbf{w}- \mathbf{w}_h\|_{1,D}=0,\quad
\lim_{h\rightarrow 0} \| q- q_h\|_{0,D}=0,\quad
\lim_{h\rightarrow 0} \|\mathbf{v}- \mathbf{v}_h\|_{1,D}=0.
\end{equation}
There are several pairs of mixed finite elements satisfying the inf-sup condition, see 
\cite{Brezzi2013} Chap. 8, \cite{Girault1986} Chap. II, in order to construct $W_h$;$Q_h$: 
$\mathbb{P}_1+bubble$;$\mathbb{P}_1$, $\mathbb{P}_2$;$\mathbb{P}_1$, $\mathbb{P}_2$;$\mathbb{P}_0$, etc.
The mixed system can be solved without knowing an explicit basis for $V_h$.
In fact, as for the Stokes equations, the global system is solved using the finite
element basis of $W_h\times Q_h$, see \cite{Girault1986}, Chap. II. 

In the following, we use the external approximation of $V_h$ based on the mixed finite elements
$\mathbb{P}_1+bubble$;$\mathbb{P}_1$.
Using density property, there exits $\mathbf{u}\in W\cap \left(H^2(D)\right)^2$ such that
$\|\mathbf{w}- \mathbf{u}\|_{1,D}$ is small. Using the interpolation operator $R_h$ for $\mathbb{P}_1$,
we get
$$
\| \mathbf{u} - R_h \mathbf{u}\|_{1,D} \leq C\, h \vert \mathbf{u} \vert_{2,D}.
$$
If we set $\mathbf{w}_h= R_h \mathbf{u}$, then $\|\mathbf{w}- \mathbf{w}_h\|_{1,D}$ is small enough.
Similarly, for $Q_h$, using density property, there exits $u\in Q\cap H^2(D)$ such that
$\| q- u\|_{0,D}$ is small. Since
$$
\| u - R_h u\|_{0,D} \leq C\, h^2 \vert u \vert_{2,D}
$$
we get that $\| q- q_h\|_{0,D}$ is small enough for $q_h=R_hu$. The approximation properties
 (\ref{density}) 
are satisfied in the above examples.

The space $V_h$ depends on $W_h$ and $Q_h$ by (\ref{V_h}). Let $\mathbf{v}$ be an element of $V$.
We solve the Stokes problem: find $\mathbf{v}_h\in W_h$ and $p_h\in Q_h$ such that
\begin{eqnarray*}
a(\mathbf{v}_h,\mathbf{w}_h) + b(\mathbf{w}_h,p_h)
&=& a(\mathbf{v},\mathbf{w}_h),\quad \forall \mathbf{w}_h\in W_h \\
b(\mathbf{v}_h,q_h)&=&0, \quad \forall q_h\in Q_h
\end{eqnarray*} 
with $b$ given by (\ref{b_D}) and $a:W\times W\rightarrow \mathbb{R}$,
$$
a(\mathbf{u},\mathbf{w})=\nu\int_D \nabla \mathbf{u} : \nabla\mathbf{w} \,d\mathbf{x}.
$$
From the second equation of Stokes problem, we have that $\mathbf{v}_h\in V_h$.
The continuous Stokes problem has $\mathbf{v}$ and $p=0$ as solution.
From \cite{Brezzi2013}, \cite{Girault1986}, we have
$$
\|\mathbf{v}- \mathbf{v}_h\|_{1,D} \leq C
\left(
\inf_{\mathbf{w}_h\in W_h} \|\mathbf{v}- \mathbf{w}_h\|_{1,D}
+\inf_{ q_h \in Q_h} \|p- q_h\|_{0,D}
\right)
$$
which shows that the last limit in (\ref{density}) holds.

We also consider $g \in \mathcal{C}(\overline{D})$, $g < 0$ on $\partial D$ and take $\omega$ of the form
\begin{equation}\label{3.omega}
  \omega=\omega_g=int\left\{ \mathbf{x}\in D;\ g(\mathbf{x})\leq 0\right\},
\end{equation}
that satisfies $\partial D  \subset \overline{\omega}$ as assumed in Section 1.
Since this definition may yield $\omega$ not connected, we choose the connected component
that contains $\partial D$ to be $\omega$ and we may assume $g$ to be positive outside it
by adding to $g$ the distance function to $\omega$, to the square and multiplied by some positive constant. We also assume $\omega$ to be Lipschitz (such details will be clarified in the next section via condition (\ref{4.1})).
  
We define (with the same $h$ as for the discretization) $H^h$ and $\tilde{H}^h:\mathbb{R}\rightarrow \mathbb{R}$ by 
\begin{equation}
\label{3.1}
H^h (r) =
\left\{
\begin{array}{ll}
1, & r \geq h,\\
\frac{(-2r+3h)r^2}{h^3} , & 0 <r < h,\\
0 & r \leq 0 ,
\end{array}
\right.
\quad
\tilde{H}^h (r) =
\left\{
\begin{array}{ll}
1, & r \geq 0,\\
\frac{(-2r+h)(r+h)^2}{h^3} , & -h < r <0,\\
0 &  r \leq -h .
\end{array}
\right.
\end{equation}
We have $0\leq H^h \leq 1$ and $0\leq \tilde{H}^h \leq 1$, they are differentiable
with respect
to their argument and they approximate $H(r)=1$ if $r\geq 0$ and $H(r)=0$ if $r < 0$ (the Heaviside function)
as $h \rightarrow 0$.

We introduce 
$a^\epsilon_h : W \times W \rightarrow \mathbb{R}$,
$c^\epsilon_h : W \times W \times W \rightarrow \mathbb{R}$,
$\tilde{c}^\epsilon_h : W \times W \times W \rightarrow \mathbb{R}$
by 
\begin{eqnarray*} 
a^\epsilon_h (\mathbf{v},\mathbf{w}) &=&
\nu \int_D [1-H^h (g)] \nabla \mathbf{v} : \nabla\mathbf{w} \,d\mathbf{x}\\
&& +\epsilon
\int_D \tilde{H}^h (g) \nabla \mathbf{v} : \nabla\mathbf{w} \,d\mathbf{x}
 +\epsilon
\int_D \tilde{H}^h (g) \mathbf{v} \cdot \mathbf{w} \,d\mathbf{x},
\\
c^\epsilon_h (\mathbf{u},\mathbf{v},\mathbf{w}) & =&
\int_D [1-\tilde{H}^h (g)]\left[ 
(\mathbf{u}\cdot\nabla ) \mathbf{v}
  \right] \cdot \mathbf{w} \,d\mathbf{x}\\
&& +\epsilon
\int_D \tilde{H}^h (g)\left[ 
(\mathbf{u}\cdot\nabla ) \mathbf{v}
  \right] \cdot \mathbf{w} \,d\mathbf{x},
\\
\tilde{c}^\epsilon_h(\mathbf{u},\mathbf{v},\mathbf{w})&=&
\frac{1}{2}c^\epsilon_h (\mathbf{u},\mathbf{v},\mathbf{w})
-\frac{1}{2}c^\epsilon_h (\mathbf{u},\mathbf{w},\mathbf{v}).
\end{eqnarray*}

We define $F_h\in W^\prime$ by
$$
\langle F_h, \mathbf{v}\rangle_{*,D} =
\int_D [1-\tilde{H}^h (g)] \mathbf{f}\cdot \mathbf{v}\,d\mathbf{x}
+\int_{\Gamma_N} \boldsymbol{\psi}\cdot \mathbf{w}\, ds
,\quad \forall \mathbf{v}\in W.
$$

Obviously, $1-H^h (g)$ is a regularization of the characteristic function of $\omega_{g}$,
that is the above defined functionals have a smooth dependence on the geometry induced by $g$ via (\ref{3.omega}).
In this way, we avoid the jumps in their coefficients at the border of $\omega_{g}$.

Let $\{ \boldsymbol{\varphi}_i \}_{i=1,\dots,n} $ be
a basis of $V_h$. We introduce the finite element approximation and regularization of (\ref{2.ns3}):
find $\mathbf{y}^\epsilon_h\in V_h$ such that
\begin{equation}
\label{3.ns3}
a^\epsilon_h(\mathbf{y}^\epsilon_h,\mathbf{v}_h)
+\tilde{c}^\epsilon_h(\mathbf{y}^\epsilon_h,\mathbf{y}^\epsilon_h,\mathbf{v}_h)
=\langle F_h, \mathbf{v}_h\rangle_{*,D},\quad \forall \mathbf{v}_h\in V_h.
\end{equation}

From the trace theorem, there exists $C_T>0$, such that
\begin{equation}
\label{trace}
\forall \mathbf{w}\in W,\ \| \mathbf{w}\|_{0,\Gamma_N}\leq C_T\| \mathbf{w}\|_{1,D}.
\end{equation}

\begin{proposition}\label{prop:3.1}
For each $0 < \epsilon$, the problem (\ref{3.ns3}) has at least one solution
$\mathbf{y}^\epsilon_h \in V_h$ and
$$
\| \mathbf{y}^\epsilon_h \|_{1,D}
\leq \frac{\| \mathbf{f} \|_{0,D} + C_T \| \boldsymbol{\psi}\|_{0,\Gamma_N}}{\alpha^\epsilon_\omega}
$$
where $\alpha^\epsilon_\omega = \min\left(\frac{\nu}{C_P^2+1},\epsilon\right)$.
\end{proposition}

\noindent
\textit{Proof.} 
Since $[1-H^h (g)]=1$ in $\omega$ and $\tilde{H}^h(g)=1$ in $\Omega$ and they are
positive everywhere, we have
for all $\mathbf{w}\in W$
$$
a^\epsilon_h(\mathbf{w},\mathbf{w}) \geq \nu \vert \mathbf{w} \vert_{1,\omega}^2
+\epsilon \| \mathbf{w} \|_{1,\Omega}^2 \geq
\frac{\nu}{C_P^2+1} \| \mathbf{w} \|_{1,\omega}^2 + \epsilon  \| \mathbf{w} \|_{1,\Omega}^2
\geq \alpha^\epsilon_\omega \| \mathbf{w} \|_{1,D}^2.
$$

If $w\in L^p(D)$, $p\geq 1$, then
$\|[1-H^h (g)]\, w \|_{L^p(D)}\leq \| w \|_{L^p(D)}$. It follows
$$
\int_D [1-H^h (g)] \nabla \mathbf{v} : \nabla\mathbf{w} \,d\mathbf{x}
\leq
\| [1-H^h (g)] \nabla \mathbf{v} \|_{0,D}
\| \nabla\mathbf{w} \|_{0,D}
\leq
\| \nabla \mathbf{v} \|_{0,D}
\| \nabla\mathbf{w} \|_{0,D}
$$
and we get
$$
\vert a^\epsilon_h (\mathbf{v},\mathbf{w}) \vert
\leq (\nu +\epsilon) \| \mathbf{v} \|_{1,D} \| \mathbf{w} \|_{1,D},
\quad \forall\mathbf{v},\mathbf{w}\in W.
$$

Similarly, we have
\begin{eqnarray*}
&&\left\vert \int_D [1-\tilde{H}^h (g)]\left[ 
(\mathbf{u}\cdot\nabla ) \mathbf{v}
  \right] \cdot \mathbf{w} \,d\mathbf{x} \right\vert
=
\left\vert \int_D \left[ 
(\mathbf{u}\cdot\nabla ) \mathbf{v}
  \right] \cdot \left( \mathbf{w}[1-\tilde{H}^h (g)]\right) \,d\mathbf{x} \right\vert\\
&&
\leq
C_d \, \| \mathbf{u} \|_{L^4(D)} \, \vert \mathbf{v} \vert_{1,D} \,
\| \mathbf{w} [1-\tilde{H}^h (g)]\|_{L^4(D)}
\leq
C_d \, \| \mathbf{u} \|_{L^4(D)} \, \vert \mathbf{v} \vert_{1,D} \,
\| \mathbf{w} \|_{L^4(D)}
\\
&&
\leq
C_D \| \mathbf{u} \|_{1,D}\| \mathbf{v} \|_{1,D} \| \mathbf{w} \|_{1,D}
\end{eqnarray*}
where $C_d$ depends just on the dimension ($d=2$ in this subsection)
and $C_D>0$ is independent of $\epsilon$ and $g$.
We get that
$$
\vert \tilde{c}^\epsilon_h (\mathbf{u},\mathbf{v},\mathbf{w}) \vert
\leq (1 +\epsilon)C_D \| \mathbf{u} \|_{1,D} \| \mathbf{v} \|_{1,D} \| \mathbf{w} \|_{1,D},
\quad \forall\mathbf{u},\mathbf{v},\mathbf{w}\in W.
$$

Rem. \ref{1.1} is also valid in $D$, that is
if $\mathbf{u}^n$ converges weakly to $\mathbf{u}$ in $H^1(D)^2$,  then
$$
\lim_{n\rightarrow\infty}
\int_D \left[ 
(\mathbf{u}^n\cdot\nabla ) \mathbf{u}^n
  \right] \cdot \mathbf{v} \,d\mathbf{x}
=\int_D \left[ 
(\mathbf{u}\cdot\nabla ) \mathbf{u}
  \right] \cdot \mathbf{v} \,d\mathbf{x},\quad \forall \mathbf{v}\in H^1(D)^2
$$
and
$$
\lim_{n\rightarrow\infty}
\int_D \left[ 
(\mathbf{u}^n\cdot\nabla ) \mathbf{v}
  \right] \cdot \mathbf{u}^n \,d\mathbf{x}
=\int_D \left[ 
(\mathbf{u}\cdot\nabla ) \mathbf{v}
  \right] \cdot \mathbf{u} \,d\mathbf{x},\quad \forall \mathbf{v}\in H^1(D)^2.
$$
But
$$
\int_D [1-\tilde{H}^h (g)]\left[ 
(\mathbf{u}^n\cdot\nabla ) \mathbf{u}^n
  \right] \cdot \mathbf{v} \,d\mathbf{x} 
=
 \int_D \left[ 
(\mathbf{u}^n\cdot\nabla ) \mathbf{u}^n
  \right] \cdot \left( \mathbf{v}[1-\tilde{H}^h (g)]\right) \,d\mathbf{x} 
$$
and
$$
\int_D [1-\tilde{H}^h (g)]\left[ 
(\mathbf{u}^n\cdot\nabla ) \mathbf{v}
  \right] \cdot \mathbf{u}^n \,d\mathbf{x} 
=
 \int_D \left[ 
(\mathbf{u}^n\cdot\nabla ) \mathbf{v}
  \right] \cdot \left( \mathbf{u}^n[1-\tilde{H}^h (g)]\right) \,d\mathbf{x}. 
$$ 
Since $\mathbf{v}[1-\tilde{H}^h (g)]\in L^2(D)^2$ and
$u_i^n u_j^n[1-\tilde{H}^h (g)]$ converges strongly to $u_iu_j[1-\tilde{H}^h (g)]$ in $L^{2}(D)$, for $i,j \in \{1,2\}$,
then
$$
\lim_{n\rightarrow\infty}
\tilde{c}^\epsilon_h(\mathbf{u}^n,\mathbf{u}^n,\mathbf{v})
=\tilde{c}^\epsilon_h(\mathbf{u},\mathbf{u},\mathbf{v}).
$$

We also obtain
\begin{eqnarray*}
  \left\vert \langle F_h, \mathbf{v}\rangle_{*,D}\right\vert
  &\leq& \| [1-\tilde{H}^h (g)]\mathbf{f} \|_{0,D} \| \mathbf{v} \|_{0,D}
+ \| \boldsymbol{\psi}\|_{0,\Gamma_N} \| \mathbf{v} \|_{0,\Gamma_N}\\
  &\leq & \| \mathbf{f} \|_{0,D}   \| \mathbf{v} \|_{1,D}
+ C_T\| \boldsymbol{\psi}\|_{0,\Gamma_N} \| \mathbf{v} \|_{1,D} ,\quad \forall\mathbf{v}\in W.
\end{eqnarray*}  

Let
$\Phi: V_h \rightarrow V_h$ be
the continuous mapping (due to the above arguments), defined by
\begin{eqnarray}
\left( \Phi(\mathbf{v}_h), \boldsymbol{\varphi}_i \right)_{H^1(D)^2}
= a^\epsilon_h(\mathbf{v}_h,\boldsymbol{\varphi}_i )
+\tilde{c}^\epsilon_h(\mathbf{v}_h,\mathbf{v}_h,\boldsymbol{\varphi}_i)
-\langle F_h, \boldsymbol{\varphi}_i\rangle_{*,D},
\label{3.4}
\end{eqnarray}  
for $i=1,\dots,n$, where $\left( \cdot, \cdot \right)_{H^1(D)^2}$ is the scalar product of $H^1(D)^2$.

If $\mathbf{y}^\epsilon_h$ is a solution of (\ref{3.ns3}), then
$$
\left( \Phi(\mathbf{y}^\epsilon_h), \boldsymbol{\varphi}_i \right)_{H^1(D)^2}=0,\ i=1,\dots,n.
$$
Since $\Phi(\mathbf{y}^\epsilon_h)\in V_h$, then $\Phi(\mathbf{y}^\epsilon_h)=0$.
Reciprocally, if $\Phi(\mathbf{y}^\epsilon_h)=0$ for $\mathbf{y}^\epsilon_h\in V_h$, we get
 (\ref{3.ns3}).

We have 
\begin{eqnarray*}
\left( \Phi(\mathbf{v}_h), \mathbf{v}_h \right)_{H^1(D)^2}
& = &a^\epsilon_h(\mathbf{v}_h, \mathbf{v}_h)
+\tilde{c}^\epsilon_h(\mathbf{v}_h,\mathbf{v}_h,\mathbf{v}_h)
-\langle F_h, \mathbf{v}_h\rangle_{*,D}\\
& = &a^\epsilon_h(\mathbf{v}_h, \mathbf{v}_h)
-\langle F_h, \mathbf{v}_h\rangle_{*,D}
\end{eqnarray*}  
since $\tilde{c}^\epsilon_h(\mathbf{v}_h,\mathbf{v}_h,\mathbf{v}_h)=0$.

As in the beginning of the proof, since $V_h\subset W_h\subset W$, we have
$$
\forall \mathbf{v}_h\in V_h,\quad a^\epsilon_h(\mathbf{v}_h,\mathbf{v}_h)\geq
\alpha^\epsilon_\omega  \| \mathbf{v}_h \|_{1,D}^2.
$$
Using
\begin{eqnarray*}
\vert \langle F_h, \mathbf{v}_h\rangle_{*,D} \vert & \leq &
\left( \| \mathbf{f} \|_{0,D} + C_T \| \boldsymbol{\psi}\|_{0,\Gamma_N}\right) \|\mathbf{v}_h \|_{1,D}
\end{eqnarray*}  
we get
\begin{eqnarray}
\left( \Phi(\mathbf{v}_h), \mathbf{v}_h \right)_{H^1(D)^2}
\geq \alpha^\epsilon_\omega  \| \mathbf{v}_h \|_{1,D}^2
-  \left( \| \mathbf{f} \|_{0,D}  + C_T \| \boldsymbol{\psi}\|_{0,\Gamma_N}\right) \|\mathbf{v}_h \|_{1,D} .
\label{3.5}
\end{eqnarray}  
We put $R^\epsilon=\frac{\| \mathbf{f} \|_{0,D} + C_T \| \boldsymbol{\psi}\|_{0,\Gamma_N}}{\alpha^\epsilon_\omega}$ and from (\ref{3.5})
we obtain $\left( \Phi(\mathbf{v}_h), \mathbf{v}_h \right)_{H^1(D)^2}\geq 0$ for all
$\mathbf{v}_h\in V_h$,
$\| \mathbf{v}_h \|_{1,D}=R^\epsilon$.
From \cite{Girault1986}, Cor. 1.1, p. 279, we get that the problem (\ref{3.ns3})
has at least one solution
$\mathbf{y}^\epsilon_h\in V_h$ and
\begin{equation}\label{3.6}
\| \mathbf{y}^\epsilon_h \|_{1,D}\leq R^\epsilon
=\frac{\| \mathbf{f} \|_{0,D}+ C_T \| \boldsymbol{\psi}\|_{0,\Gamma_N}}{\alpha^\epsilon_\omega}.
\end{equation}

This ends the proof.
 \quad $\Box$

\begin{proposition}\label{prop:3.2}
Let $\{ \mathcal{T}_h\}_{h>0}$ be a family of regular triangulation of $D$.  
Let $0 <\theta <1$ be fixed. We assume that $\epsilon$ is fixed and satisfies
(\ref{2.H5}). 
Under the hypotheses (\ref{2.H4}) and (\ref{2.H3}), then
\begin{equation}\label{3.7}
  \lim_{h\rightarrow 0} \mathbf{y}^\epsilon_h=\mathbf{y}^\epsilon, 
\end{equation}
weakly in $H^1(D)^2$ and strongly in $L^2(D)^2$, where $\mathbf{y}^\epsilon$ is the
solution of (\ref{2.ns3}).
\end{proposition}

\noindent
\textit{Proof.} 
From Prop. \ref{prop:2.1} and \ref{prop:2.3iii}, we get that
the problem (\ref{2.ns3}) has a unique solution $\mathbf{y}^\epsilon$.
Since $\epsilon$ is fixed, from (\ref{3.6}), we get that $\mathbf{y}^\epsilon_h$
is bounded in $H^1(D)^2$ and we can extract a subsequence weakly convergent to some
$\widehat{\mathbf{y}}^\epsilon \in W$. We will prove that $\widehat{\mathbf{y}}^\epsilon \in V$.
Let $q\in Q$ be an arbitrary element, then there exists $\{q_h\}_{h>0}$, such that
$q_h\in Q_h$ and $\lim_{h\rightarrow 0} q_h=q$ strongly in $Q$. From the definition of $V_h$ (\ref{V_h}),
we have $b(\mathbf{y}^\epsilon_h, q_h)=0$ and by passing to the limit $h\rightarrow 0$, we get
$b(\widehat{\mathbf{y}}^\epsilon,q)=0$ for all $q\in Q$. Consequently,
$\nabla\cdot\widehat{\mathbf{y}}^\epsilon=0$ in $D$, then $\widehat{\mathbf{y}}^\epsilon \in V$.

Let $\mathbf{v}$ be in $V$. From (\ref{density}),
there exists $\mathbf{v}_h\in V_h$ such that $\lim_{h\rightarrow 0}\mathbf{v}_h=\mathbf{v}$
strongly in $H^1(D)^2$.
We have
$$
\vert \langle F_\omega, \mathbf{v}\rangle_{*,\omega} - \langle F_h, \mathbf{v}_h\rangle_{*,D} \vert \leq
\vert \langle F_\omega, \mathbf{v}\rangle_{*,\omega} - \langle F_h, \mathbf{v} \rangle_{*,D} \vert
+\vert\langle F_h, \mathbf{v} \rangle_{*,D} - \langle F_h, \mathbf{v}_h\rangle_{*,D} \vert
$$
and
$$
\vert \langle F_h, \mathbf{v}-\mathbf{v}_h\rangle_{*,D} \vert  \leq 
\left( \| \mathbf{f} \|_{0,D} + C_T \| \boldsymbol{\psi}\|_{0,\Gamma_N}\right) \|\mathbf{v}-\mathbf{v}_h \|_{1,D}
$$
together with
$$
\langle F_\omega, \mathbf{v}\rangle_{*,\omega} - \langle F_h, \mathbf{v} \rangle_{*,D}=
\int_\omega \mathbf{f}\cdot\mathbf{v}\,d\mathbf{x} 
- \int_D [1-\tilde{H}^h (g)] \mathbf{f}\cdot \mathbf{v}\,d\mathbf{x}.
$$
Since $[1-\tilde{H}^h (g)]\mathbf{f}\cdot\mathbf{v}\in L^1(D)$ converges to
$[1-H(g)]\mathbf{f}\cdot\mathbf{v}\in L^1(D)$,
almost everywhere (a.e.) in $D$ and
$\vert[1-\tilde{H}^h (g)]\mathbf{f}\cdot\mathbf{v}\vert\leq \vert\mathbf{f}\cdot\mathbf{v}\vert\in L^1(D)$,
By the Lebesgue dominated convergence theorem, see \cite{Brezis2005}, p. 54,
we get $\lim_{h\rightarrow 0} \langle F_h, \mathbf{v}\rangle_{*,D}= \langle F_\omega, \mathbf{v}\rangle_{*,\omega}$,
then
$$
\lim_{h\rightarrow 0} \langle F_h, \mathbf{v}_h\rangle_{*,D}= \langle F_\omega, \mathbf{v}\rangle_{*,\omega}.
$$

We also have
\begin{eqnarray*}
a^\epsilon_h(\mathbf{y}^\epsilon_h,\mathbf{v}_h)-a^\epsilon(\widehat{\mathbf{y}}^\epsilon,\mathbf{v}) & = &
a^\epsilon_h(\mathbf{y}^\epsilon_h,\mathbf{v}_h-\mathbf{v}) 
+ a^\epsilon_h(\mathbf{y}^\epsilon_h,\mathbf{v})-a^\epsilon(\widehat{\mathbf{y}}^\epsilon,\mathbf{v})
\end{eqnarray*}
and
\begin{eqnarray*}
  \vert a^\epsilon_h(\mathbf{y}^\epsilon_h,\mathbf{v}_h-\mathbf{v})\vert
  &\leq & (\nu +\epsilon) \|\mathbf{y}^\epsilon_h\|_{1,D}\|\mathbf{v}_h-\mathbf{v} \|_{1,D}
  \leq (\nu +\epsilon)R^\epsilon\|\mathbf{v}_h-\mathbf{v} \|_{1,D}.
\end{eqnarray*} 
We have $\lim_{h\rightarrow 0} \tilde{H}^h (g)=\lim_{h\rightarrow 0} H^h(g)=H(g)$ a.e. in $D$.
We employ the Lemma 6.1 from \cite{HMT2016}:
if $a,a_n\in L^\infty(D)$, $\left\| a_n\right\|_{0,\infty,D} \leq M$, $a_n\rightarrow a$ a.e.
in $D$, $b_n\rightarrow b$ weakly in $L^2(D)$ and $v\in L^2(D)$, then
$$
\lim_{n\rightarrow \infty}\int_D a_n b_n v \,dx=\int_D abv \,dx.
$$
Here $\left\| \cdot \right\|_{0,\infty,D} $ is the norm of $L^\infty(D)$.
We can apply this Lemma for $a_n=\tilde{H}^h (g)$ or $a_n=H^h (g)$,
$b_n$ components of $\widehat{\mathbf{y}}^\epsilon_h$
or $\nabla\widehat{\mathbf{y}}^\epsilon_h$,
and $v$ components of $\mathbf{v}$, then
we get
$$
\lim_{h\rightarrow 0} a^\epsilon_h(\widehat{\mathbf{y}}^\epsilon_h,\mathbf{v}) 
=\lim_{h\rightarrow 0} a^\epsilon(\widehat{\mathbf{y}}^\epsilon,\mathbf{v})
$$
on a subsequence. Finally, we get
$$
\lim_{h\rightarrow 0} a^\epsilon_h(\mathbf{y}^\epsilon_h,\mathbf{v}_h)
=a^\epsilon(\widehat{\mathbf{y}}^\epsilon,\mathbf{v}).
$$

We have proved in Prop. \ref{prop:3.1} that
$$
\vert \tilde{c}^\epsilon_h (\mathbf{u},\mathbf{v},\mathbf{w}) \vert
\leq (1 +\epsilon)C_D \| \mathbf{u} \|_{1,D} \| \mathbf{v} \|_{1,D} \| \mathbf{w} \|_{1,D},
\quad \forall \mathbf{u},\mathbf{v},\mathbf{w}\in W.
$$
We obtain here that
$$
\vert\tilde{c}^\epsilon_h(\mathbf{y}^\epsilon_h,\mathbf{y}^\epsilon_h,\mathbf{v}_h-\mathbf{v}) \vert
\leq (1 +\epsilon)C_D
(R^\epsilon)^2 \|\mathbf{v}_h-\mathbf{v} \|_{1,D}
$$
and
\begin{eqnarray*}
  \tilde{c}^\epsilon_h(\mathbf{y}^\epsilon_h,\mathbf{y}^\epsilon_h,\mathbf{v}_h)
  &=& \tilde{c}^\epsilon_h(\mathbf{y}^\epsilon_h,\mathbf{y}^\epsilon_h,\mathbf{v}_h-\mathbf{v})
  +\tilde{c}^\epsilon_h(\mathbf{y}^\epsilon_h,\mathbf{y}^\epsilon_h,\mathbf{v}).
\end{eqnarray*}  
But
\begin{eqnarray*}
\tilde{c}^\epsilon_h(\mathbf{y}^\epsilon_h,\mathbf{y}^\epsilon_h,\mathbf{v}) & = &
\frac{1}{2} c_D([1-\tilde{H}^h (g)]\mathbf{y}^\epsilon_h,\mathbf{y}^\epsilon_h,\mathbf{v})
+\epsilon \frac{1}{2} c_D(\tilde{H}^h (g)\mathbf{y}^\epsilon_h,\mathbf{y}^\epsilon_h,\mathbf{v})\\
&& 
-\frac{1}{2} c_D([1-\tilde{H}^h (g)]\mathbf{y}^\epsilon_h,\mathbf{v},\mathbf{y}^\epsilon_h)
-\epsilon \frac{1}{2} c_D( \tilde{H}^h (g)\mathbf{y}^\epsilon_h,\mathbf{v},\mathbf{y}^\epsilon_h)
\end{eqnarray*}
where $c_D$ is defined as $c_\omega$ by replacing $\omega$ by $D$.

As in the Rem. \ref{rem:1.1}, using that the injection $H^1(D)\subset L^4(D)$ is compact,
we obtain $\mathbf{y}^\epsilon_h \rightarrow \widehat{\mathbf{y}}^\epsilon$
strongly in $L^4(D)^2$ and $\nabla\mathbf{y}^\epsilon_h \rightarrow \nabla\widehat{\mathbf{y}}^\epsilon$
weakly in $L^2(D)^4$.
We obtain that
$$
(\mathbf{y}^\epsilon_h)_i(\mathbf{y}^\epsilon_h)_j
\rightarrow (\widehat{\mathbf{y}}^\epsilon)_i(\widehat{\mathbf{y}}^\epsilon)_j\quad
i,j=1,2
$$
strongly in $L^2(D)$ and also weakly in $L^2(D)$.
We employ the above Lemma with $a_n=1-\tilde{H}^h (g)$,
$b_n=(\mathbf{y}^\epsilon_h)_i(\mathbf{y}^\epsilon_h)_j$ and 
$v=(\nabla \mathbf{v})_{ji}=\frac{\partial v_j}{\partial x_i}$, then
$$
\lim_{h\rightarrow 0} c_D([1-\tilde{H}^h (g)]\mathbf{y}^\epsilon_h,\mathbf{v},\mathbf{y}^\epsilon_h)
= c_D([1-H(g)]\widehat{\mathbf{y}}^\epsilon,\mathbf{v},\widehat{\mathbf{y}}^\epsilon)
=c_\omega(\widehat{\mathbf{y}}^\epsilon,\mathbf{v},\widehat{\mathbf{y}}^\epsilon).
$$

Using Egorov's Theorem, see \cite{Brezis2005} p. 75, for all $\delta>0$, there exists
$D_\delta\subset D$ such that $meas(D\setminus D_\delta)<\delta$ and $1-H(g)$ converges
uniformly on $D_\delta$. We have
\begin{eqnarray*}
\int_D[1-\tilde{H}^h (g)](\mathbf{y}^\epsilon_h)_i(\nabla \mathbf{y}^\epsilon_h)_{ji}
(\mathbf{v})_j d\mathbf{x}
&=&\int_{D\setminus D_\delta} [1-\tilde{H}^h (g)](\mathbf{y}^\epsilon_h)_i(\nabla \mathbf{y}^\epsilon_h)_{ji}
(\mathbf{v})_j d\mathbf{x}\\
&&+\int_{D_\delta} [1-\tilde{H}^h (g)](\mathbf{y}^\epsilon_h)_i(\nabla \mathbf{y}^\epsilon_h)_{ji}
(\mathbf{v})_j d\mathbf{x}.
\end{eqnarray*}
In $D_\delta$, $[1-\tilde{H}^h (g)](\mathbf{v})_j$ converges to $[1-H(g)](\mathbf{v})_j$
strongly in $L^4(D_\delta)$ and $(\mathbf{y}^\epsilon_h)_i(\nabla \mathbf{y}^\epsilon_h)_{ji}$
converges to $(\widehat{\mathbf{y}}^\epsilon)_i(\widehat{\mathbf{y}}^\epsilon\nabla )_{ji}$
weakly in $L^{4/3}(D_\delta)$. It follows that, on a subsequence,
$$
\lim_{h\rightarrow 0}\int_{D_\delta} [1-\tilde{H}^h (g)](\mathbf{y}^\epsilon_h)_i
(\nabla \mathbf{y}^\epsilon_h)_{ji}(\mathbf{v})_j d\mathbf{x}
=\int_{D_\delta}[1-H(g)](\widehat{\mathbf{y}}^\epsilon)_i
(\nabla \widehat{\mathbf{y}}^\epsilon)_{ji}(\mathbf{v})_j d\mathbf{x}.
$$

In $D\setminus D_\delta$, we have
\begin{eqnarray*}
&&\left\vert \int_{D\setminus D_\delta} [1-\tilde{H}^h (g)](\mathbf{y}^\epsilon_h)_i
(\nabla \mathbf{y}^\epsilon_h)_{ji}(\mathbf{v})_j d\mathbf{x} \right\vert \\
&\leq &\| (\mathbf{y}^\epsilon_h)_i\|_{L^4(D\setminus D_\delta)}
\| (\nabla \mathbf{y}^\epsilon_h)_{ji} \|_{L^2(D\setminus D_\delta)}
\| (\mathbf{v})_j\|_{L^4(D\setminus D_\delta)}\\
&\leq& \| \mathbf{y}^\epsilon_h\|_{L^4(D\setminus D_\delta)}
\| \nabla \mathbf{y}^\epsilon_h \|_{L^2(D\setminus D_\delta)}
\| (\mathbf{v})_j\|_{L^4(D\setminus D_\delta)}
\end{eqnarray*}
From the continuity of the injection $H^1(D\setminus D_\delta)\subset L^4(D\setminus D_\delta)$ 
$$
\|\mathbf{y}^\epsilon_h\|_{L^4(D\setminus D_\delta)} \leq C_1 \|\mathbf{y}^\epsilon_h\|_{1,D\setminus D_\delta}
\leq C_1 \|\mathbf{y}^\epsilon_h\|_{1,D}
$$
and we have that $\|\mathbf{y}^\epsilon_h\|_{1,D}\leq R^\epsilon$.
Also we have, see \cite{Adams1975}, Th. 28, p. 25,
$$
\| (\mathbf{v})_j\|_{L^4(D\setminus D_\delta)} \leq meas(D\setminus D_\delta)^{1/4}
\| (\mathbf{v})_j\|_{L^2(D\setminus D_\delta)} 
$$
then
$$\left\vert \int_{D\setminus D_\delta} [1-\tilde{H}^h (g)](\mathbf{y}^\epsilon_h)_i
(\nabla \mathbf{y}^\epsilon_h)_{ji}(\mathbf{v})_j d\mathbf{x} \right\vert
\leq \delta^{1/4} C_1 (R^\epsilon)^2 \| (\mathbf{v})_j\|_{L^2(D\setminus D_\delta)}
$$
and consequently, on a subsequence,
$$
\lim_{h\rightarrow 0} c_D([1-\tilde{H}^h (g)]\mathbf{y}^\epsilon_h,\mathbf{y}^\epsilon_h,\mathbf{v})
= c_D([1-H(g)]\widehat{\mathbf{y}}^\epsilon,\widehat{\mathbf{y}}^\epsilon,\mathbf{v})
=c_\omega(\widehat{\mathbf{y}}^\epsilon,\widehat{\mathbf{y}}^\epsilon,\mathbf{v}).
$$

We have similar results for the terms
$c_D(\tilde{H}^h (g)\mathbf{y}^\epsilon_h,\mathbf{y}^\epsilon_h,\mathbf{v})$ and
$c_D( \tilde{H}^h (g)\mathbf{y}^\epsilon_h,\mathbf{v},\mathbf{y}^\epsilon_h)$, then
$$
\lim_{h\rightarrow 0} \tilde{c}^\epsilon_h(\mathbf{y}^\epsilon_h,\mathbf{y}^\epsilon_h,\mathbf{v})
= \tilde{c}^\epsilon(\widehat{\mathbf{y}}^\epsilon,\widehat{\mathbf{y}}^\epsilon,\mathbf{v}).
$$
By passing to the limit on a subsequence $h\rightarrow 0$ as before, we get
$$
\lim_{h\rightarrow 0} \tilde{c}^\epsilon_h(\mathbf{y}^\epsilon_h,\mathbf{y}^\epsilon_h,\mathbf{v}_h)
= \tilde{c}^\epsilon(\widehat{\mathbf{y}}^\epsilon,\widehat{\mathbf{y}}^\epsilon,\mathbf{v}).
$$

Consequently, by passing to the limit on a subsequence $h\rightarrow 0$
in (\ref{3.ns3})
we obtain that $\widehat{\mathbf{y}}^\epsilon=\mathbf{y}^\epsilon$ the solution of (\ref{2.ns3}).
 As this is unique, the convergence is valid without taking subsequences.
 \quad $\Box$

\section{Applications}
\setcounter{equation}{0}

\subsection{A numerical example\label{sec:3.1}}

For the numerical tests, we have used the finite element software FreeFem++, \cite{freefem++}.
We consider the bounded domain $D\subset\mathbb{R}^2$ with the boundary
$\partial D=\Gamma_1 \cup \Gamma_2 \cup \Gamma_3 \cup \Gamma_4$, defined by:
left side $\Gamma_1=\{-0.5\}\times ]-0.5,0.5[$,
bottom side $\Gamma_2=]-0.5,0.5[ \times \{-0.5\}$,
top side $\Gamma_4=]-0.5,0.5[ \times \{0.5\}$ and the right side
\begin{eqnarray*}
\Gamma_3&=&\left\{ \left( x_1(t),x_2(t)\right),
\ t\in \left]-\frac{\pi}{2},\frac{\pi}{2}\right[;\right.\\
&& \left. x_1(t)=0.5+ 0.5\cos(t),\ x_2(t)=0.5\sin(t)\right\}.
\end{eqnarray*}

For $\Omega \subset \subset D$ not necessarily connected, we set
$\omega=D\setminus\overline{\Omega}$.
Using the notations from the first section, we solve the Navier-Stokes system in $\omega$
for $\Gamma_D=\Gamma_2 \cup \Gamma_3 \cup \Gamma_4$, $\Gamma_N=\Gamma_1$,
the imposed traction $\boldsymbol{\psi}=\left( 100y, 0 \right)$ 
  on $\Gamma_N$,
homogeneous Dirichlet boundary condition on $\Gamma_D$
and homogeneous Neumann boundary condition on $\partial\Omega$,
viscosity $\nu=1$ and $\mathbf{f}=(0,0)$.

\begin{figure}[ht]
\centering
  \includegraphics[width=6.5cm]{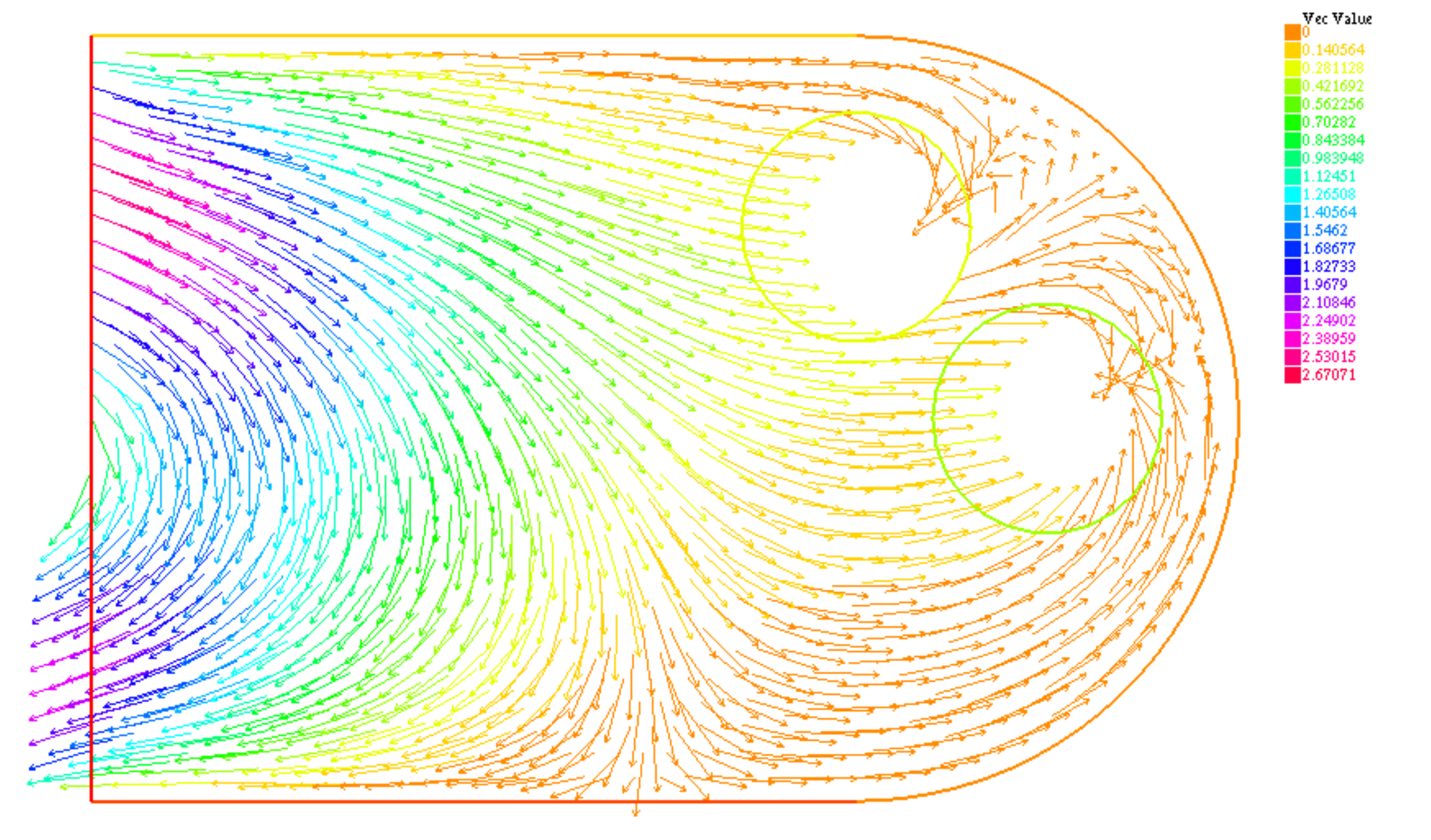}\quad
  \includegraphics[width=6.5cm]{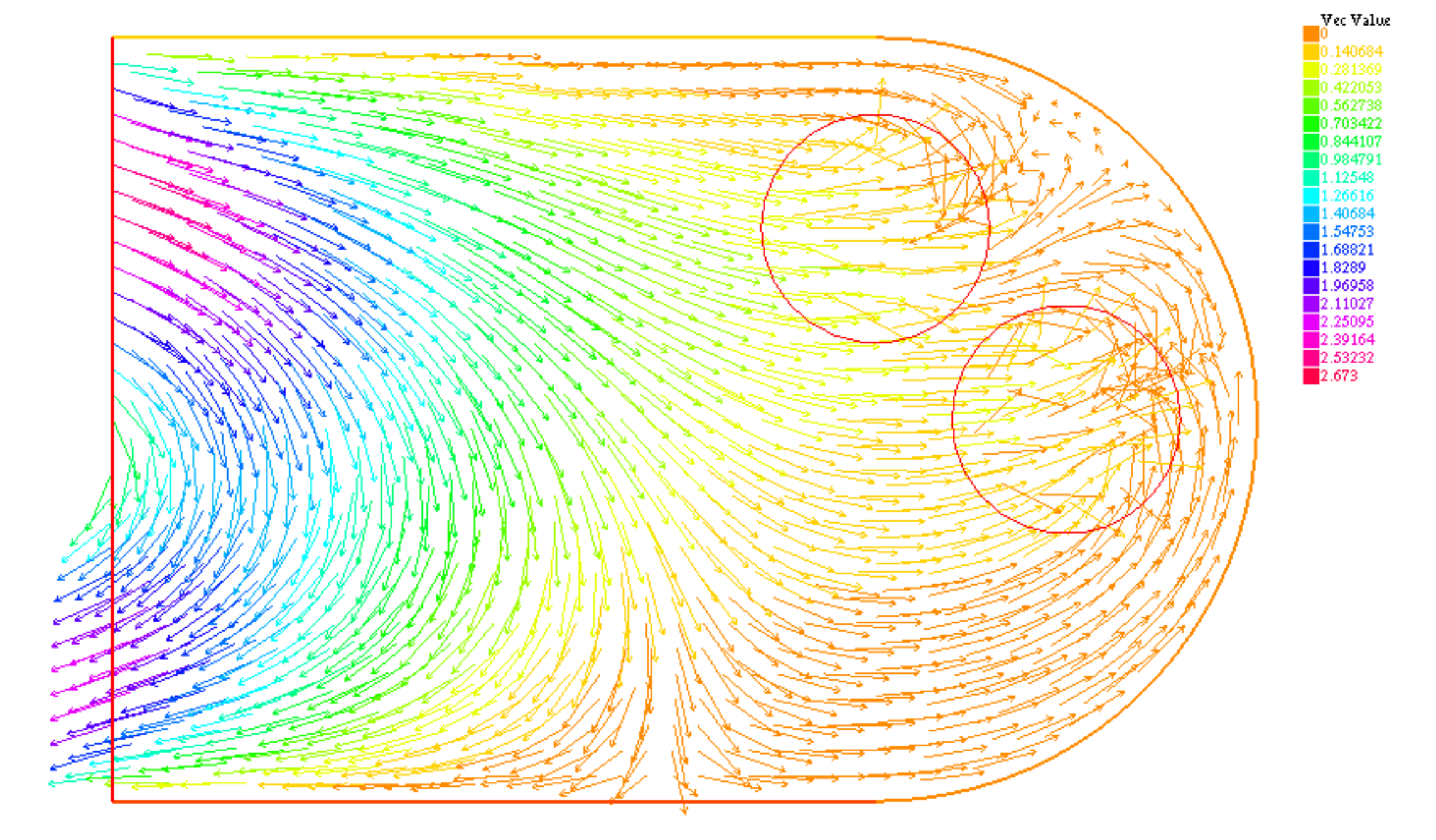}
\caption{The velocities of the Navier-Stokes equations in $\omega$ (left).
  The velocities of the Navier-Stokes equations in $D$; we have plotted $\partial\Omega$, too (right).}
\label{fig:test3-vit}
\end{figure}

Here, $\Omega$ is composed by two disks $\Omega_1$ and $\Omega_2$ of radius $0.15$ and
centers $(0.5,0.25)$ and $(0.75,0)$, respectively. See Figure \ref{fig:test3-vit}.
We use the parametrization 
$$
g=\max \{ -(x_1-0.5)^2-(x_2-0.25)^2+0.15^2,-(x_1-0.75)^2-(x_2)^2+0.15^2 \}.
$$

In order to compute numerically $\widehat{\mathbf{y}}_\omega\in V_1$, we start
  from (\ref{1.ns1_v1})
and we treat the constraints concerning the normal flux on each internal boundary
$\partial \Omega_i$ by a Lagrangian multiplier: find $\widehat{\mathbf{y}}^\omega_h\in W^\omega_h$,
$\widehat{p}^\omega_h\in Q^\omega_h$,
$\ell_i\in \mathbb{R}$, $i=1,2$, such that
\begin{eqnarray}
a^\omega_h(\widehat{\mathbf{y}}^\omega_h,\mathbf{w}^\omega_h)
+b^\omega_h(\mathbf{w}^\omega_h,\widehat{p}^\omega_h)
&&
\nonumber\\
+\tilde{c}^\omega_h(\widehat{\mathbf{y}}^\omega_h,\widehat{\mathbf{y}}^\omega_h,\mathbf{w}^\omega_h)
+\sum_{i=1}^2 \ell_i\int_{\partial \Omega_i}\mathbf{w}^\omega_h\cdot \mathbf{n}\,ds
&=&\langle F_\omega, \mathbf{w}^\omega_h\rangle_{*,\omega},\forall \mathbf{w}^\omega_h\in W^\omega_h
\label{3.1.ns_mult_1}\\
b^\omega_h(\widehat{\mathbf{y}}^\omega_h,q^\omega_h) &= & 0,\forall q^\omega_h\in Q^\omega_h
\label{3.1.ns_mult_2}\\
\int_{\partial \Omega_i}\widehat{\mathbf{y}}^\omega_h\cdot \mathbf{n}\,ds &=&0,\quad i=1,2.
\label{3.1.ns_mult_3}
\end{eqnarray}

We have employed the mixed finite elements $\mathbb{P}_1+bubble$ for the velocity and
$\mathbb{P}_1$ for the pressure.
We solve the nonlinear Navier-Stokes system by the Newton method and initializing
the iterations by the Stokes solution and we get
$\widehat{\mathbf{y}}^\omega_h$,
a finite element approximation
of $\widehat{\mathbf{y}}_\omega:\overline{\omega}\rightarrow\mathbb{R}^2$.
We also solve the discrete version of problem (\ref{2.20})-(\ref{2.21}) in $D$,
in order to get $\mathbf{y}^\epsilon_h:\overline{D}\rightarrow\mathbb{R}^2$,
but with $b$ given by (\ref{b_D}) in place of $b^\epsilon$.

To get the convergence for the Newton method, 2-3 iterations are required.
The fluid velocities are plotted in Figure \ref{fig:test3-vit}.
Then, we compute the $L^2$ and $H^1$ relative errors
$$
L^2\ err_{rel}=\frac{\| \widehat{\mathbf{y}}^\omega_h-\mathbf{y}^\epsilon_h\|_{L^2(\omega)}}
{\| \widehat{\mathbf{y}}^\omega_h\|_{L^2(\omega)}},\quad
H^1\ err_{rel}=\frac{\| \widehat{\mathbf{y}}^\omega_h-\mathbf{y}^\epsilon_h\|_{H^1(\omega)}}
{\| \widehat{\mathbf{y}}^\omega_h\|_{H^1(\omega)}}
$$
for various choices of $\epsilon$ and $h$.

We point out that, for the computed solution in $\omega$, we have
$\int_{\Gamma_N} \widehat{\mathbf{y}}^\omega_h\cdot \mathbf{n}\,ds=6.9e-17$,
$\int_{\partial\Omega_1} \widehat{\mathbf{y}}^\omega_h\cdot \mathbf{n}\,ds=2.7e-17$,
$\int_{\partial\Omega_2} \widehat{\mathbf{y}}^\omega_h\cdot \mathbf{n}\,ds=2.9e-17$, where
$\mathbf{n}$ is the unit outward normal to the boundary $\partial\omega$, then
total normal flux
$\int_{\partial\omega} \widehat{\mathbf{y}}^\omega_h\cdot \mathbf{n}\,ds=1.2e-16$.
We have $\vert\vert\nabla\cdot \widehat{\mathbf{y}}^\omega_h \vert\vert_{0,\omega} =0.139678$.

When we solve (\ref{2.20})-(\ref{2.21}) with $b$ given by (\ref{b_D}) in place of $b^\epsilon$,
$\nabla\mathbf{y}^\epsilon_h$ as well as the pressure $p^\epsilon_h$
are discontinuous on $\partial\Omega$.
We have obtained
$\int_{\Gamma_N} \mathbf{y}^\epsilon_h\cdot \mathbf{n}\,ds=2.7e-16$ and
$\vert\vert\nabla\cdot y^\epsilon_h \vert\vert_{0,\omega} =0.142572$.

\begin{figure}[ht]
\centering
  \includegraphics[width=6cm]{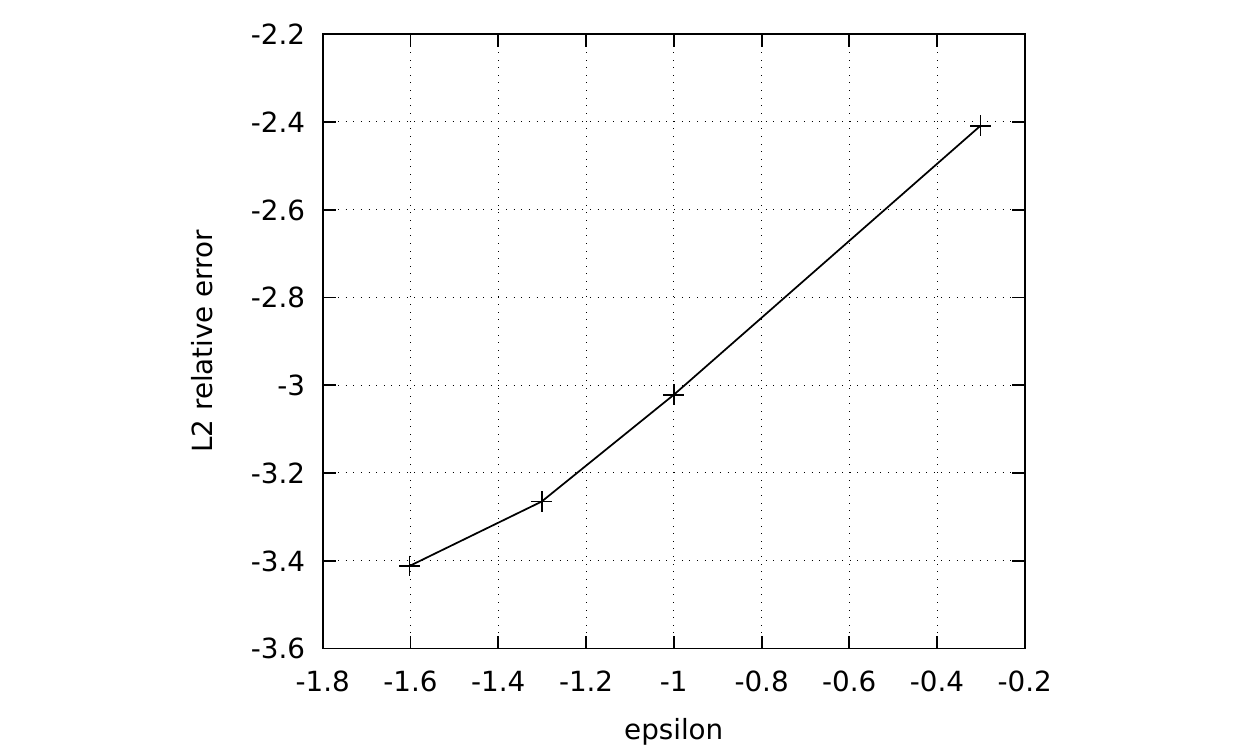}\quad
  \includegraphics[width=6cm]{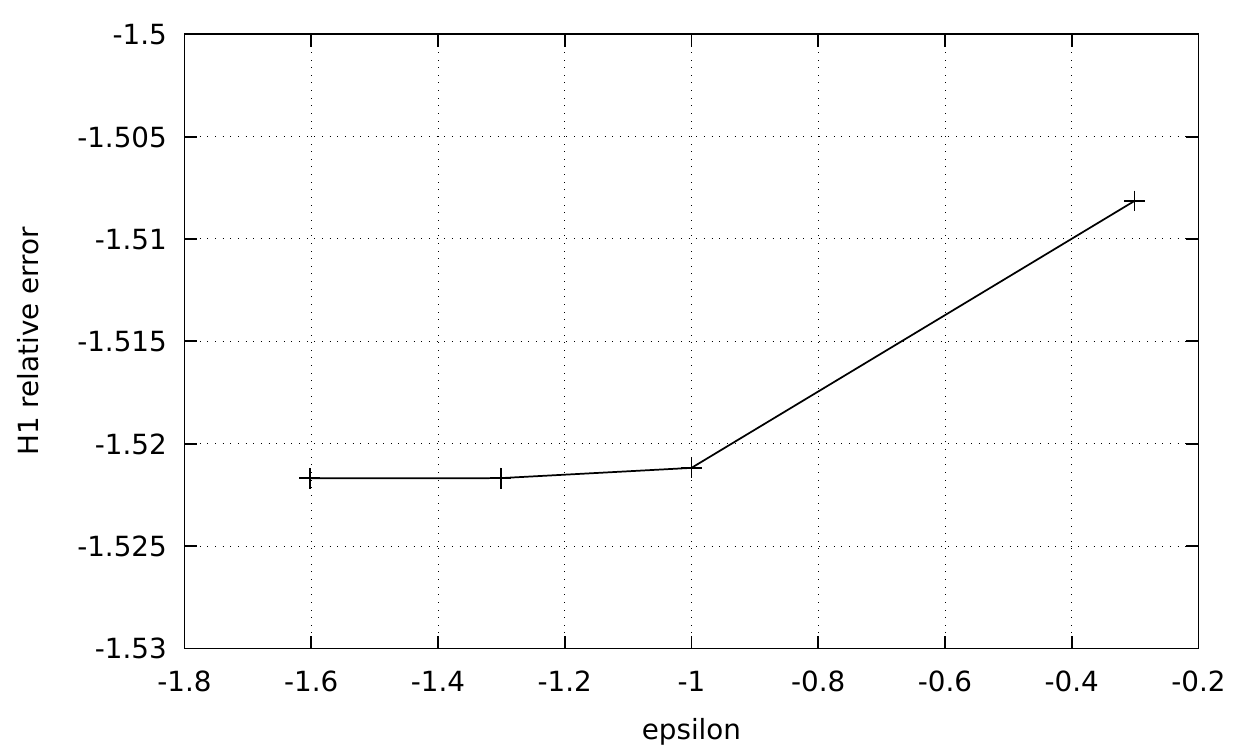}
\caption{The $L^2$ (left) and $H^1$ (right) relative errors for different
  $\epsilon$ in $\log_{10}$ scale.}
\label{fig:test3_c_epsilon}
\end{figure}

We have solved the problem in $D$ for different
$\epsilon$: 0.5, 0.1, 0.05, 0.025, on a mesh with 30270 triangles and 15361 vertices.
The $L^2$ and $H^1$ relative errors are plotted in Figure \ref{fig:test3_c_epsilon} and
the slope of the regression line (least squares method, see \cite{Theodor1982}, Ch. 4)
are $0.79$ and $0.01$, respectively. 

\begin{figure}[ht]
\centering
  \includegraphics[width=4cm]{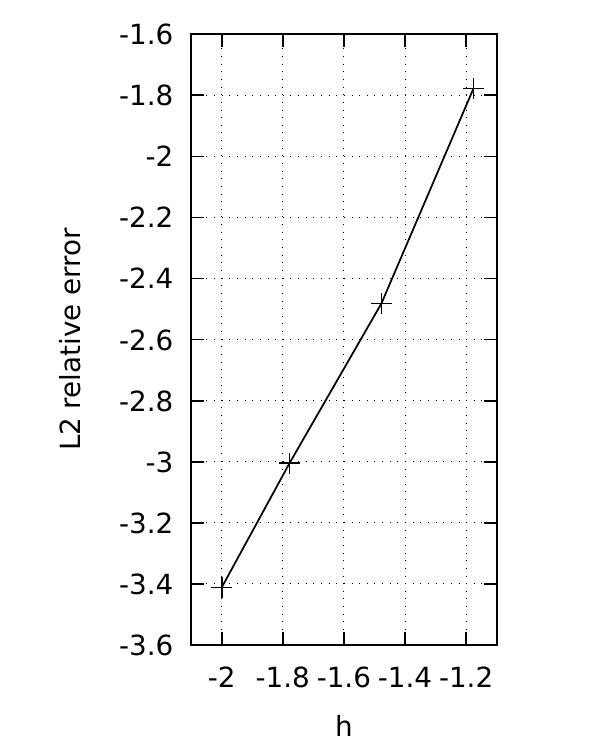}\quad
  \includegraphics[width=7cm]{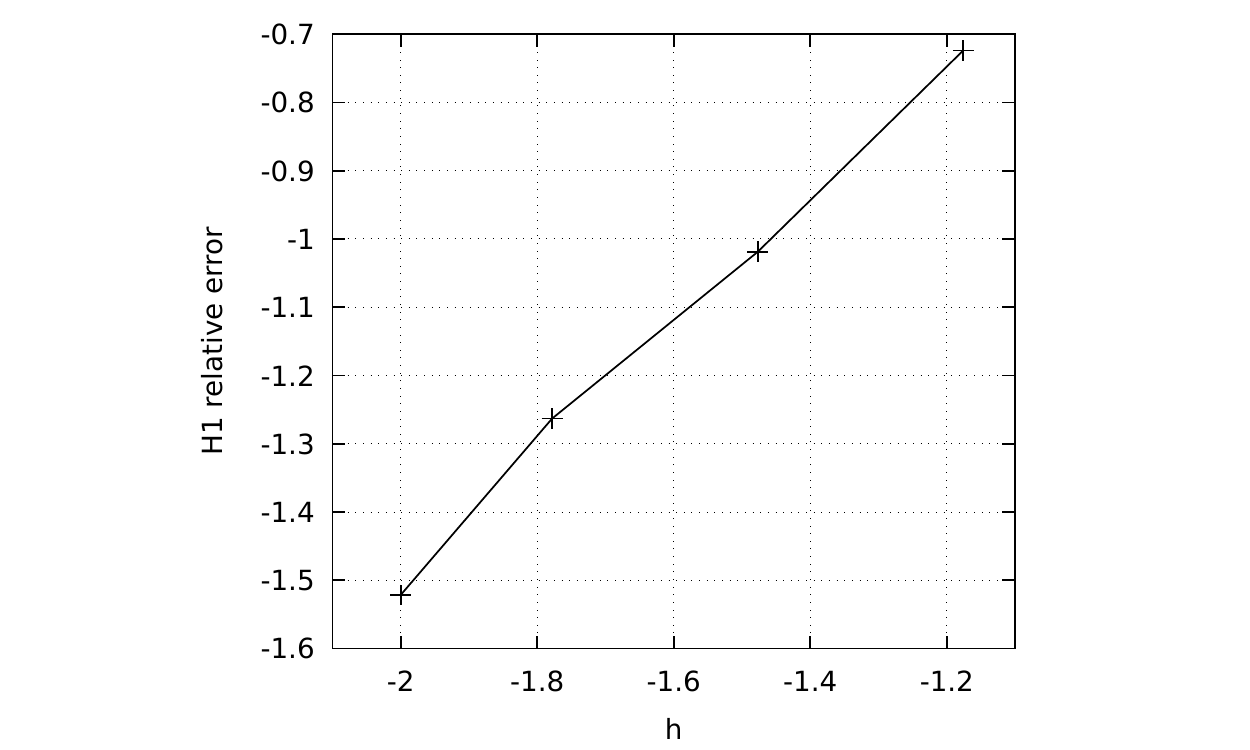}
\caption{The $L^2$ (left) and $H^1$ (right) relative errors for different
  mesh size in $\log_{10}$ scale.}
\label{fig:test3_c_h}
\end{figure}

We have also solved the problem in $D$ for $\epsilon=0.025$ and different meshes.
We show the $L^2$ and $H^1$ relative errors in Figure \ref{fig:test3_c_h} and
the slope of the regression line are $1.96$ and $0.94$, respectively.

\subsection{Shape and topology optimization}

We associate to the system (\ref{1.1})-(\ref{1.5}) (or (\ref{1.ns3}) in the weak formulation),
a geometric optimization problem defined over the class of all admissible domains
$\omega\in \mathcal{O}$, $\omega\subset D\subset \mathbb{R}^2$, (a given family of domains).
Equivalently, the family of admissible obstacles $\Omega=D\setminus\overline{\omega}$ may be indicated.

We consider a family $\mathcal{F}\subset \mathcal{C}(\overline{D})$ of level functions and
the admissible open sets $\omega=\omega_g$ are defined by (\ref{3.omega}).
Here , we impose more regularity on the admissible level functions:
$\mathcal{F}\subset \mathcal{C}^1(\overline{D})$ and $D$ is a polygonal domain and
\begin{eqnarray}
&& \vert\nabla g\vert\neq 0\hbox{ on }
  G_g=\{ \mathbf{x}\in D;\ g(\mathbf{x})=0\},\ \forall g\in \mathcal{F};
  \label{4.1}\\
&& g(\mathbf{x}) < 0,\ \forall\mathbf{x}\in \partial D,\ \forall g\in \mathcal{F}.
  \label{4.2} 
\end{eqnarray}
Under hypotheses (\ref{4.1}), (\ref{4.2}), we have $\partial D\subset \overline{\omega}_g$
(given by (\ref{3.omega})) and this fixes the connected component of the open set $\omega_g$ to be taken into
account in (\ref{3.omega}) as the domain $\omega_g$. Moreover, $\omega_g$ may not be simply connected  and
$\partial \omega_g \setminus \partial D $ is
of class $\mathcal{C}^1$ by the implicit function theorem and (\ref{4.1}).
The relation (\ref{3.omega}) may be written as
\begin{equation}\label{4.3}
\omega_g=\left\{ \mathbf{x}\in D;\ g(\mathbf{x}) < 0\right\}, \ \forall g\in \mathcal{F}.
\end{equation}

By (\ref{4.1}), (\ref{4.3}), one can obtain a global parametrization of the boundary $\partial \omega_g$
via a simple ordinary differential Hamiltonian system. This is useful in boundary observation problems
and other questions, see \cite{MT_2022}.
In (\ref{4.3}), we also use the convention that $g>0$ outside the chosen component of $\omega_g$.
This is possible since in (\ref{4.3}), we may add to $g$ the square of the distance
function, multiplied by a positive constant, to the component of $\omega_g$ with $\partial D\subset \partial \omega_g$ and
this preserves $g\in\mathcal{C}^1(\overline{D})$ and (\ref{4.1}), (\ref{4.2}),
(\ref{4.3}). Moreover, (\ref{4.1}) ensures that the number of connected components of the open set defined in (\ref{3.omega}) or the number of holes  of $\omega_g$ are finite, \cite{MT_2022}.

To (\ref{1.ns3}) and its solution $\mathbf{y}_g\in V_{\omega_g}$, general cost
functionals (to be minimized) may be associated
\begin{equation}\label{4.4}
  \min_{g\in \mathcal{F}}\int_\Lambda j\left(\mathbf{x}, \mathbf{y}_g(\mathbf{x})\right)d\mathbf{x}
  (ds),
\end{equation}
where one may choose, for instance, $\Lambda=\omega_g$, $\Lambda=E\subset \omega_g$
some prescribed subdomain, $\Lambda=\partial \omega_g\setminus \partial D$ or
some prescribed open part $\gamma$ of it, etc.
In order that such choices are possible, one should impose simple algebraic conditions on
$\mathcal{F}$:
\begin{eqnarray}
&& 
 g(\mathbf{x})<0,\ \forall \mathbf{x}\in E,\ \forall g\in \mathcal{F},
  \label{4.5}\\
&& g(\mathbf{x}) = 0,\ \forall\mathbf{x}\in \gamma,\ \forall g\in \mathcal{F}.
  \label{4.6} 
\end{eqnarray}
Moreover, $j:\overline{D}\times \mathbb{R}\rightarrow\mathbb{R}$ is a Carath\'eodory
function, bounded from below by a constant and $j(\mathbf{x},\cdot)\in \mathcal{C}^1(\mathbb{R})$.

The optimal design problem that we discuss here is given by (\ref{1.ns3}), (\ref{4.4})
and a family of admissible controls $g \in \mathcal{F}\subset \mathcal{C}^1(\overline{D})$ that satisfy
(\ref{4.1}), (\ref{4.2}), (\ref{4.5}) or (\ref{4.6}).
More constraints, for instance on the state $\mathbf{y}_g\in V_{\omega_g}$, may be added.
The setting that we consider here was developed for general linear elliptic equations in
\cite{N_Tiba2012}, \cite{MT_2022}. In \cite{MT_2021} we have studied a simpler case associated
to Stokes flows.

Notice that $\omega_g$ defined in (\ref{4.3}) is not necessarily simply connected and
$\Omega_g=D\setminus \overline{\omega}_g$ is not necessarily connected.
That's why the problem (\ref{1.ns3}), (\ref{4.4}) combines shape optimization with
topology optimization and the following numerical experiments confirm this.
We quote as well the monograph of Plotnikov and Sokolowski \cite{PS} and
the recent papers \cite{DFOP}, \cite{SN2022}, p. 1875 that stress the difficulty of the
topological aspects in geometric
optimization problems, especially when the state system represents models from
fluid mechanics.

At the numerical level, solving shape optimization or topology optimization
problems meets a high computational cost due to the necessity to update the discrete
mesh in each iteration and to recompute the mass matrix, that are very time
consuming. Such difficulties stimulated the development of the so-called fixed domain
methods and the penalization approach is an important example in this direction.
Initially proposed for the solution of boundary value problems in domains with
complicated geometry, Mignot \cite{mi}, Astrakmantsev \cite{Ak}, the fixed domain methods may be very
efficient in unknown or variable domain problems like geometric optimization or free
boundary problems, see the survey \cite{N_Tiba2012} .

We follow here the approach ``first discretize, then optimize'' and we analyze the
problem (\ref{1.ns3}), (\ref{4.4}) when $\Lambda=\omega_g$, under hypotheses (\ref{4.1}), (\ref{4.2}).

For the numerical computations, we employ the mixed finite elements
$\mathbb{P}_1+bubble$ for the velocity
$\mathbf{y}^\epsilon_h \in W_h$ and
$\mathbb{P}_1$ elements for the pressure $p^\epsilon_h \in Q_h$.
For $g \in \mathcal{F}\subset \mathcal{C}^1(\overline{D})$, we define the discretization $g_h \in M_h$, where $M_h$ is obtained using 
$\mathbb{P}_1$ finite elements.
We denote $dim(W_h)=2N_1$, $dim(Q_h)=N_2$ and $dim(M_h)=N_3$, $M=2N_1+N_2$ and $N=M+N_3$.

In Section \ref{sec:2.2}, we have used the basis
$\{ \boldsymbol{\varphi}_i \}_{i=1,\dots,n} $ of $V_h$. But numerically this is disadvantageous and
we use here the basis $\{ \boldsymbol{\phi}_i \}_{i=1,\dots,2N_1} $ of $W_h$ defined
in the previous sentences and similarly for $Q_h$ and $M_h$, that is we employ usual finite
element basis composed by hat functions. 
The vectors associated to $\mathbf{y}^\epsilon_h$, $p^\epsilon_h$, $g_h$ are
$Y_h\in \mathbb{R}^{2N_1}$, $P_h\in \mathbb{R}^{N_2}$, $G_h\in \mathbb{R}^{N_3}$, respectively.

We use the definition of $a^\epsilon_h$, $\tilde{c}^\epsilon_h$, $F_h$ from Section \ref{sec:2.2},
but with $g$ replaced by its discretization $g_h$.
Also, we set $b^\epsilon_h : W \times Q \rightarrow \mathbb{R}$, by
$$
b^\epsilon_h(\mathbf{w},q)=
-\int_D [1-\tilde{H}^h(g_h)] (\nabla \cdot \mathbf{w})q\,d\mathbf{x}
-\epsilon \int_D \tilde{H}^h(g_h) (\nabla \cdot \mathbf{w})q\,d\mathbf{x}.
$$

We introduce:
\begin{eqnarray*}
  \mathcal{A}(G_h) &=&
  \left( a^\epsilon_h(\boldsymbol{\phi}_j,\boldsymbol{\phi}_i) \right)_{1\leq i,j\leq 2N_1}
\in \mathbb{R}^{2N_1\times 2N_1}\\
\tilde{\mathcal{C}}_1(G_h,Y_h)&=&\left( \tilde{c}^\epsilon_h\left(\sum_{k=1}^{2N_1} Y_k\boldsymbol{\phi}_k,
\boldsymbol{\phi}_j,\boldsymbol{\phi}_i\right) \right)_{1\leq i,j\leq 2N_1}
\in \mathbb{R}^{2N_1\times 2N_1}\\
\mathcal{B}(G_h) &=&
\left( b^\epsilon_h(\boldsymbol{\phi}_j,\pi_i) \right)_{1\leq i\leq N_2,1\leq j\leq 2N_1}
\in \mathbb{R}^{N_2\times 2N_1}\\
\mathcal{L}(G_h)&=&\left(\langle F_h, \boldsymbol{\phi}_i\rangle_{*,D}\right)_{1\leq i\leq 2N_1}
\in \mathbb{R}^{2N_1}.
\end{eqnarray*}

We consider a discretization of the minimization problem
\begin{equation}\label{4.7bis}
\inf_{G_h\in \mathbb{R}^{N_3}} J_h(G_h)
\end{equation}
subject to the algebraic formulation obtained from the mixed version of formulation (\ref{3.ns3})
(with $g_h$ in the place of $g$):
find $(Y_h,P_h) \in \mathbb{R}^{2N_1}\times \mathbb{R}^{N_2}$
such that
\begin{equation}\label{3.8bis}
\left(
\begin{array}{c}
  \mathcal{A}(G_h)Y_h + \tilde{\mathcal{C}}_1(G_h,Y_h) Y_h + \mathcal{B}^T(G_h)P_h\\ 
  \mathcal{B}(G_h)Y_h
\end{array}  
\right)
=
\left(
\begin{array}{c}
  \mathcal{L}(G_h)\\
  0
\end{array}  
\right)  
\end{equation}
and the constraints
\begin{eqnarray}
&&  g_h < 0\hbox{ on }\partial D\label{4.8}\\
&& \forall g_h,\ \nexists T\in \mathcal{T}_h,\ T \hbox{ triangle},
  \ g_h(\mathbf{x})=0,\forall \mathbf{x}\in T.
  \label{4.10}
\end{eqnarray}

Notice that $\mathcal{A}(G_h),\ \tilde{\mathcal{C}}_1(G_h,Y_h),\ \mathcal{B}(G_h),\ \mathcal{L}(G_h),
\ Y_h,\ P_h$
depend on $G_h$ 
according to the definitions before.
Moreover, we shall not impose (\ref{4.10}) in the sequel since, by (\ref{4.3}) and the
convention $g > 0$ in $D \setminus \omega_g$ and $g < 0$ in  $\omega_g$, (\ref{4.10}) is
automatically fulfilled by $g_h$, for $\mathbb{P}_1$  elements. 

Usually, the approximating minimization problem (\ref{4.7bis}), (\ref{3.8bis}) is interpreted as a discretized constrained
optimal control problem with the control $G_h\in \mathbb{R}^{N_3}$ satisfying 
(\ref{4.8}), describing the unknown geometry and
acting in the coefficients.
However, (\ref{3.8bis}) has no known
uniqueness property. The uniqueness argument from Prop. \ref{prop:2.3iii} is  sharp and seems not
possible to be extended to the discretized and regularized formulation (\ref{3.8bis}).
Such control problems are known in the literature as singular control problems (they may
lack uniqueness or even
existence for the state system, for certain data), see \cite{Li}, \cite[Ch. 3.1.3]{NST}. In the singular control approach the minimization is performed with respect to all unknowns $(Y_h,P_h,G_h)$ as independent variables and (\ref{3.8bis}) is interpreted as a restriction.
It is also to be noticed that the discrete optimization problem 
 (\ref{4.7bis})-(\ref{4.8}) may have no global optimal
solution due to the strict inequality on $\partial D$.

First, we minimize the dissipated energy given by
$j\left(\mathbf{x}, \mathbf{y}(\mathbf{x})\right)
=\mathbf{e}(\mathbf{y}(\mathbf{x})):\mathbf{e}(\mathbf{y}(\mathbf{x}))$
where $\mathbf{e}(\mathbf{y})=\frac{1}{2}\left( \nabla \mathbf{y} + (\nabla \mathbf{y})^T \right)$
and we set
\begin{equation}\label{3.18}
  \mathcal{J}_h(Y_h,P_h,G_h)=
  \int_D [1-\tilde{H}^h(g_h)] j\left(\mathbf{x}, \mathbf{y}^\epsilon_h(\mathbf{x})\right)d\mathbf{x},
\end{equation}

\noindent
that also achieves an approximating extension of the cost functional to $D$.

We define $C:\mathbb{R}^N\rightarrow \mathbb{R}^M$ by
\begin{equation}\label{3.19}
C(Y_h,P_h,G_h)=
\left(
\begin{array}{c}
  \mathcal{A}(G_h)Y_h + \tilde{\mathcal{C}}_1(G_h,Y_h) Y_h + \mathcal{B}^T(G_h)P_h -\mathcal{L}(G_h)\\
  \mathcal{B}(G_h)Y_h
\end{array}  
\right).
\end{equation}
We include the Dirichlet boundary condition for the velocity by
choosing 1 on the diagonal and 0 otherwise
in the lines of $\mathcal{A}(G_h)+\tilde{\mathcal{C}}_1(G_h,Y_h)$
corresponding to the components $Y_i$, in the Dirichlet boundary nodes.
We also put 0 on the corresponding lines of
the vector $\mathcal{L}(G_h)\in \mathbb{R}^{2N_1}$, considered to be a matrix $2N_1 \times 1$.

To simplify, we denote $X=(Y_h,P_h,G_h)\in \mathbb{R}^N$, then the constrained optimization problem
to be solved is
\begin{eqnarray}
  \inf_{X\in \mathbb{R}^N}  \mathcal{J}_h(X)\label{3.20}\\
  C(X)=0\label{3.21}
\end{eqnarray}
where $C(X)=\left( c_1(X),\dots , c_M(X)\right)^T$.

We can treat the above constrained optimization problem by the classical penalization of the
constraints in the cost functional, \cite{bert},
\begin{equation}
\inf_{X\in \mathbb{R}^N}  \mathcal{J}_\rho (X) = \mathcal{J}_h(X) +\frac{\rho}{2} C^T(X)C(X)\label{A}
\end{equation}
with $\rho >0$. The gradient of the penalized function is
$$
\nabla\mathcal{J}_\rho (X) = \nabla\mathcal{J}_h(X) +\rho \left(jac\, C(X)\right)^T C(X)
$$
and we can employ the steepest descent method.

We have to provide
$\nabla\mathcal{J}_h$, $jac\, C$, $X^0$,
where $\nabla\mathcal{J}_h(X)\in \mathbb{R}^N$ is the gradient of $\mathcal{J}_h(X)$,
$jac\, C(X)\in \mathbb{R}^{M \times N}$ is the jacobian matrix of $C(X)$ and $X^0$ is the
initial iteration of $X$.
It is also possible to impose (\ref{4.8}), but for the numerical tests we have
removed this condition, for simplicity. And
since $g_h$ is not modified near $\partial D$ in certain examples below and can be checked a posteriori
 or it plays no essential role even when violated.

The gradient of the objective function is computed by
$$
\nabla\mathcal{J}_h(X)=
\left(
\begin{array}{c}
  \int_D [1-\tilde{H}^h(g_h)] 2\mathbf{e}(\mathbf{y}^\epsilon_h):
  \mathbf{e}(\boldsymbol{\phi}_i) d\mathbf{x},\quad 1\leq i\leq 2N_1\\
  0\in \mathbb{R}^{N_2}\\
  \int_D [-\tilde{H}^h]^\prime (g_h)\gamma_i\, \mathbf{e}(\mathbf{y}^\epsilon_h):
  \mathbf{e}(\mathbf{y}^\epsilon_h) d\mathbf{x},\quad 1\leq i\leq N_3  
\end{array}
\right).
$$

Before introducing the $jac\, C$, we define
\begin{eqnarray*}
\tilde{\mathcal{C}}_2(G_h,Y_h)&=&\left( \tilde{c}^\epsilon_h\left(\boldsymbol{\phi}_j,
\sum_{k=1}^{2N_1} Y_k\boldsymbol{\phi}_k, \boldsymbol{\phi}_i\right) \right)_{1\leq i,j\leq 2N_1}
\in \mathbb{R}^{2N_1\times 2N_1}
\end{eqnarray*}
\begin{eqnarray*}
\mathcal{A}^\prime(G_h,Y_h)&=&
\left(
\nu \int_D [-{H}^h]^\prime (g_h)\gamma_j
\nabla \mathbf{y}^\epsilon_h : \nabla \boldsymbol{\phi}_i \,d\mathbf{x}
\right.\\
&&\left.
+\epsilon
\int_D  [\tilde{H}^h]^\prime (g_h)\gamma_j
\nabla \mathbf{y}^\epsilon_h : \nabla \boldsymbol{\phi}_i \,d\mathbf{x}
 +\epsilon
\int_D  [\tilde{H}^h]^\prime (g_h)\gamma_j \mathbf{y}^\epsilon_h \cdot \boldsymbol{\phi}_i \,d\mathbf{x}
\right)
\end{eqnarray*}
for $1\leq i\leq 2N_1,1\leq j\leq N_3$,
\begin{eqnarray*}
  \tilde{\mathcal{C}}^\prime (G_h,Y_h)&=&
  \frac{1}{2}\mathcal{C}_2^\prime (G_h,Y_h)-\frac{1}{2}\mathcal{C}_3^\prime (G_h,Y_h)
\end{eqnarray*}
where
\begin{eqnarray*}
  \mathcal{C}_2^\prime (G_h,Y_h) &=&
\int_D  [-\tilde{H}^h]^\prime (g_h)\gamma_j
\left[ (\mathbf{y}^\epsilon_h\cdot\nabla) \mathbf{y}^\epsilon_h 
  \right] \cdot \boldsymbol{\phi}_i \,d\mathbf{x}\\
&& +\epsilon
\int_D [\tilde{H}^h]^\prime (g_h)\gamma_j
\left[ (\mathbf{y}^\epsilon_h\cdot\nabla )\mathbf{y}^\epsilon_h 
  \right] \cdot \boldsymbol{\phi}_i \,d\mathbf{x} \\
\mathcal{C}_3^\prime (G_h,Y_h) &=&
\int_D  [-\tilde{H}^h]^\prime (g_h)\gamma_j
\left[ (\mathbf{y}^\epsilon_h\cdot\nabla ) \boldsymbol{\phi}_i
  \right] \cdot \mathbf{y}^\epsilon_h  \,d\mathbf{x}\\
&& +\epsilon
\int_D [\tilde{H}^h]^\prime (g_h)\gamma_j
\left[ (\mathbf{y}^\epsilon_h\cdot\nabla ) \boldsymbol{\phi}_i
  \right] \cdot  \mathbf{y}^\epsilon_h  \,d\mathbf{x} 
\end{eqnarray*}
for $1\leq i\leq 2N_1,1\leq j\leq N_3$ and moreover
\begin{eqnarray*}
&&\mathcal{B}^\prime(G_h,Y_h)=\\
&&\left(
-\int_D [-\tilde{H}^h]^\prime (g_h)\gamma_j  (\nabla \cdot \mathbf{y}^\epsilon_h) \pi_i\,d\mathbf{x}
-\epsilon \int_D [\tilde{H}^h]^\prime (g_h)\gamma_j (\nabla \cdot \mathbf{y}^\epsilon_h) \pi_i\,d\mathbf{x}
\right)
\end{eqnarray*}
for $1\leq i\leq N_2,1\leq j\leq N_3$,
\begin{eqnarray*}
&&\mathcal{BT}^\prime(G_h,P_h)=\\
&&\left(
-\int_D [-\tilde{H}^h]^\prime (g_h)\gamma_j  (\nabla \cdot \boldsymbol{\phi}_i) p_h\,d\mathbf{x}
-\epsilon \int_D [\tilde{H}^h]^\prime (g_h)\gamma_j (\nabla \cdot \boldsymbol{\phi}_i) p_h\,d\mathbf{x}
\right)
\end{eqnarray*}
for $1\leq i\leq 2N_1,1\leq j\leq N_3$,
\begin{eqnarray*}
\mathcal{L}^\prime(G_h)&=&
\left(
\int_D [-\tilde{H}^h]^\prime (g_h)\gamma_j  \mathbf{f}\cdot \boldsymbol{\phi}_i
\,d\mathbf{x}
\right)  
\end{eqnarray*}  
for $1\leq i\leq 2N_1,1\leq j\leq N_3$.

Then, the jacobian matrix is
$$
jac\, C(X)=
\left(
\begin{array}{ccc}
  \mathcal{A}(G_h) +\tilde{\mathcal{C}}_1(G_h,Y_h)+\tilde{\mathcal{C}}_2(G_h,Y_h)
  & \mathcal{B}^T(G_h)
  &  jac_{13}\,C(X) \\
\mathcal{B}(G_h) & 0 & \mathcal{B}^\prime(G_h,Y_h)
\end{array}
\right)  
$$
where
$$
jac_{13}\,C(X)=\mathcal{A}^\prime(G_h,Y_h) + \tilde{\mathcal{C}}^\prime (G_h,Y_h)
  + \mathcal{BT}^\prime(G_h, P_h) -\mathcal{L}^\prime(G_h).
$$
We have to modify the lines of $jac\, C(X)$ corresponding to the Dirichlet boundary condition
for the velocity, similarly to the modifications of $C(X)$: 1 on the diagonal, 0 otherwise on the lines
of $\mathcal{A}(G_h) +\tilde{\mathcal{C}}_1(G_h,Y_h)+\tilde{\mathcal{C}}_2(G_h,Y_h)$,
and 0 on the lines of $jac_{13}\,C(X)$ corresponding to the components $Y_i$, where $Y_i=0$ is imposed.

\newpage
\textbf{Test 1.}
\medskip

We use the same domain $D$ as in Section \ref{sec:3.1}, but
$\Omega$ is composed by two disks $\Omega_1$ and $\Omega_2$ of radius $0.1$ and
centers $(-0.2,0.2)$ and $(-0.2,-0.2)$.
As in Section \ref{sec:3.1}, we set the viscosity $\nu=1$, the body force $\mathbf{f}=(0,0)$,
the homogeneous Dirichlet boundary condition on $\Gamma_D$,
the homogeneous Neumann boundary condition on $\partial\Omega$, but
the imposed traction is now $\boldsymbol{\psi}=\left( 100\,x_2, 0 \right)$ on $\Gamma_N$.

The penalization parameter is $\rho=0.8$.
The mesh has 44722 triangles, 22632 vertices. The optimization
problem (\ref{3.20})-(\ref{3.21}) has $N=179972$ and $M=157340$.
We have used $\epsilon=0.01$ and we start with the solution of the Navier-Stokes
equations in the fixed domain $D$ obtained for the initial parametrization of
$\Omega$ and in this case we have $\| C(X^0) \|_\infty=1.9e-12$.
For given $G_h$, the nonlinear system in $(Y_h,P_h)$ was solved by the Newton method.
It is the only place where the Navier-Stokes equations are solved in this way.

\begin{figure}[ht]
\centering
  \includegraphics[width=10cm]{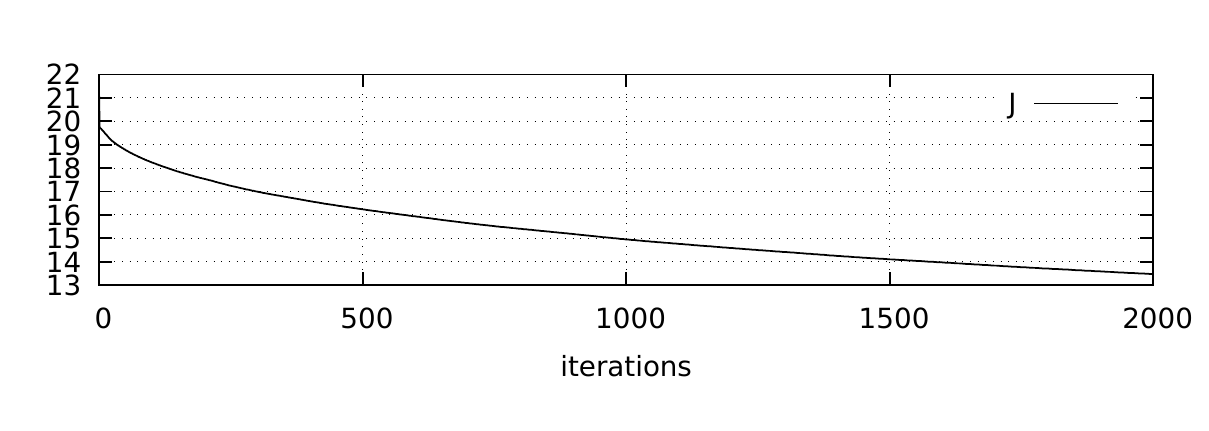}
\caption{Test 1. The history of $\mathcal{J}_h(X)$.}
\label{fig:test1_J}
\end{figure}

\begin{figure}[ht]
\centering
  \includegraphics[width=10cm]{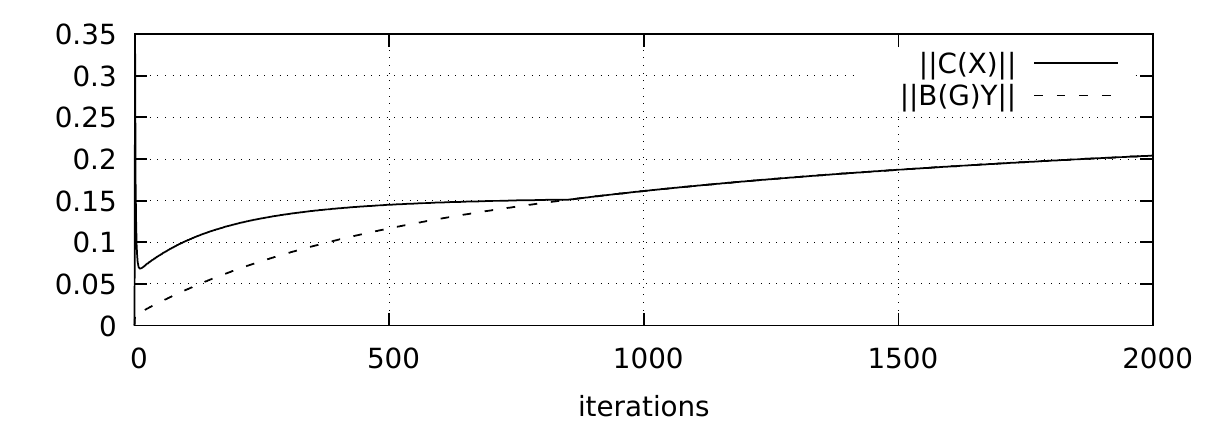}
\caption{Test 1. The history of $\| C(X) \|_\infty$ and $\| \mathcal{B}(G_h)Y_h \|_\infty$.}
\label{fig:test1_CX}
\end{figure}

In Figure \ref{fig:test1_J}, we can see the evolution of the objective function
$\mathcal{J}_h(X)$. At the initial iteration we have $\mathcal{J}_h(X^0)=21.083848$.
A second quantity to be observed is $\| C(X) \|_\infty$ in Figure \ref{fig:test1_CX},
that may be interpreted
as an error related indicator. After 2000 iteration, its value is
$\| C(X^{2000}) \|_\infty=0.204184$. 
In addition, we can observe $\| \mathcal{B}(G_h)Y_h \|_\infty$ in Figure \ref{fig:test1_CX}, which appears in the
second line of (\ref{3.19})
and represents the discrete equation involving the divergence of the
fluid velocity. The components of the vector
$\mathcal{B}(G_h)Y_h$ are $b^\epsilon_h(\mathbf{y}^\epsilon_h,\pi_i)$.
Of course, we have $\| \mathcal{B}(G_h)Y_h \|_\infty \leq \| C(X) \|_\infty$.
We have also computed $\vert\vert\nabla\cdot y^\epsilon_h \vert\vert_{0,\omega} = 0.223471$. 

\begin{figure}[ht]
\centering
  \includegraphics[width=5cm]{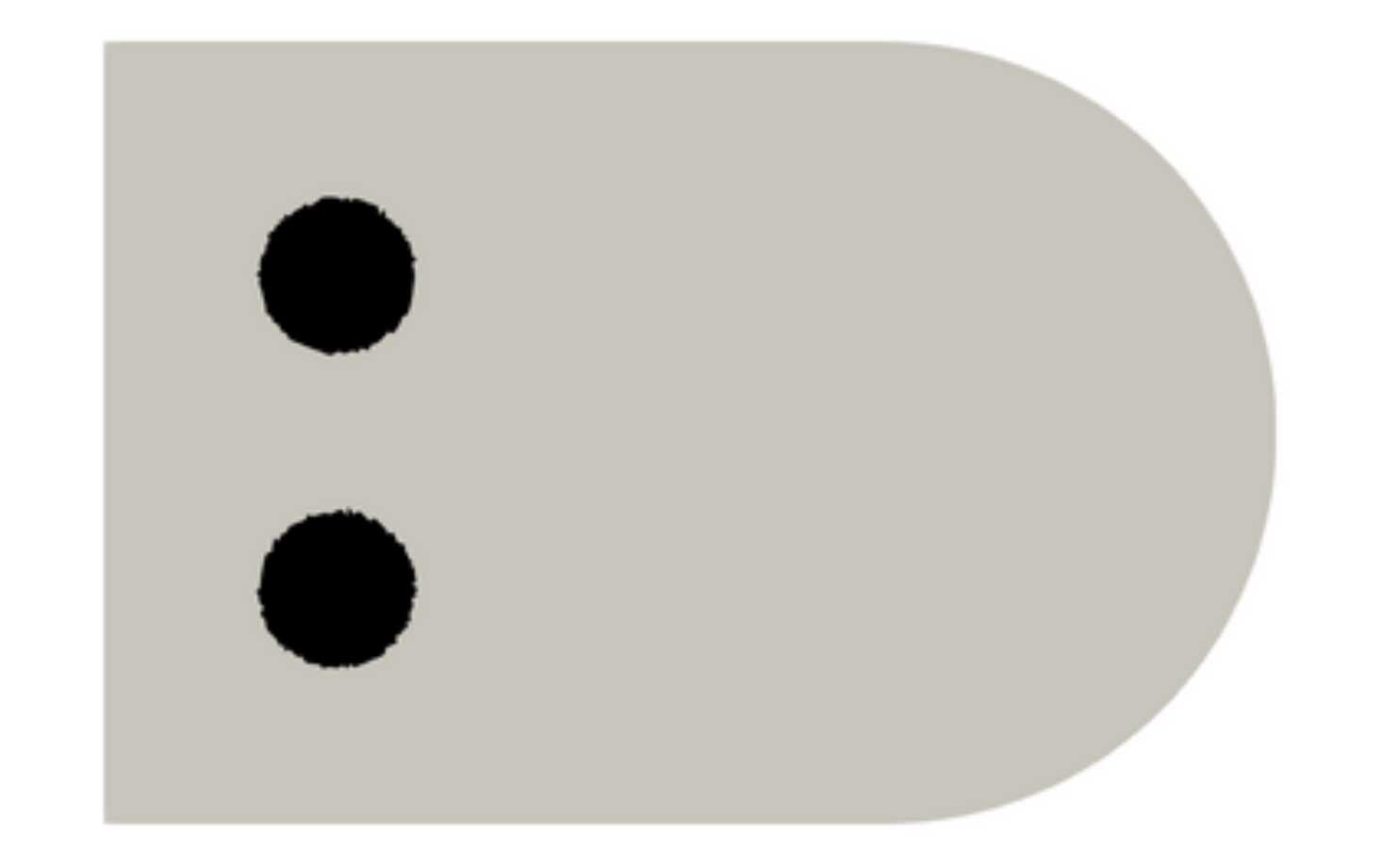}\quad
  \includegraphics[width=5cm]{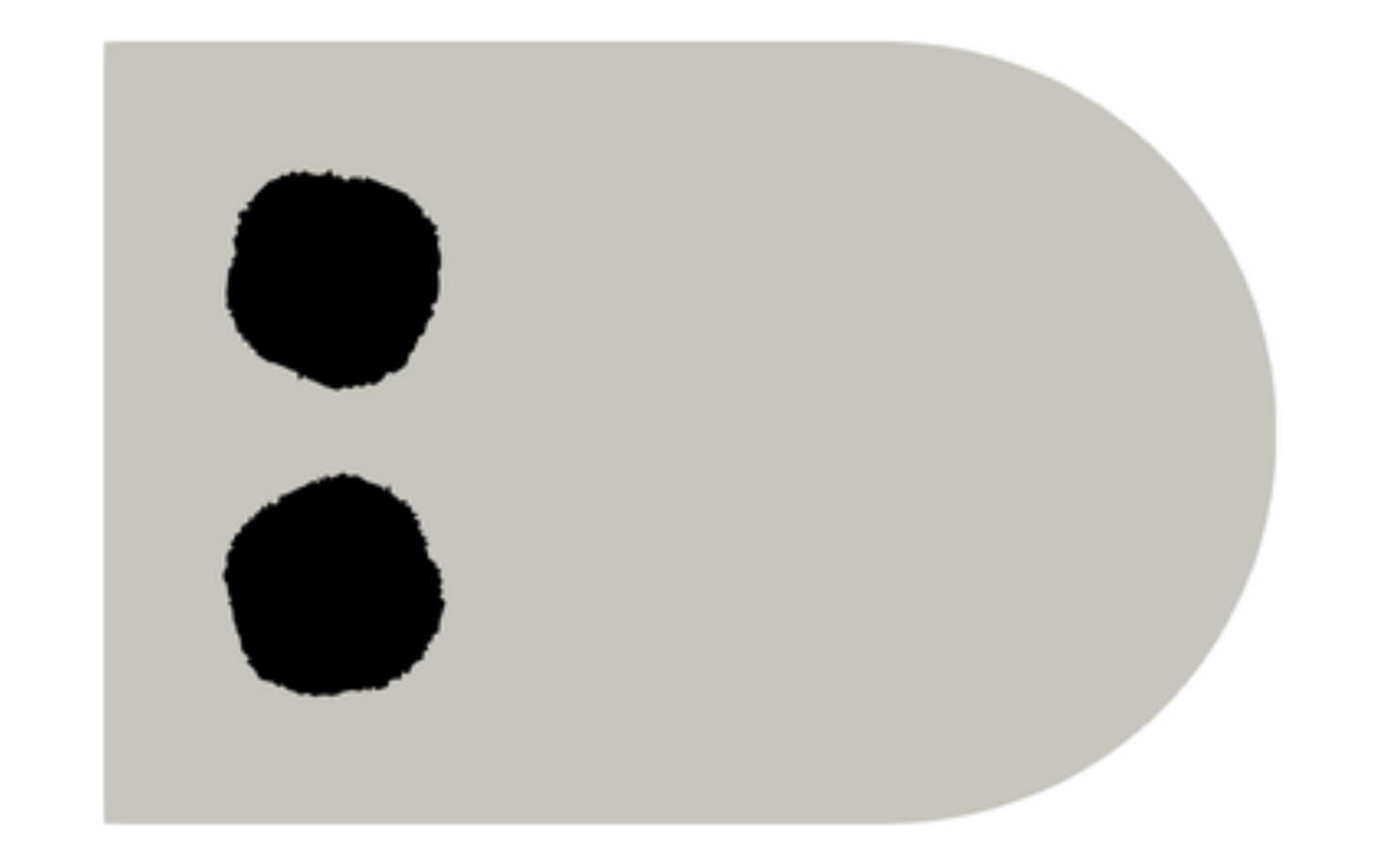}\\
  \includegraphics[width=5cm]{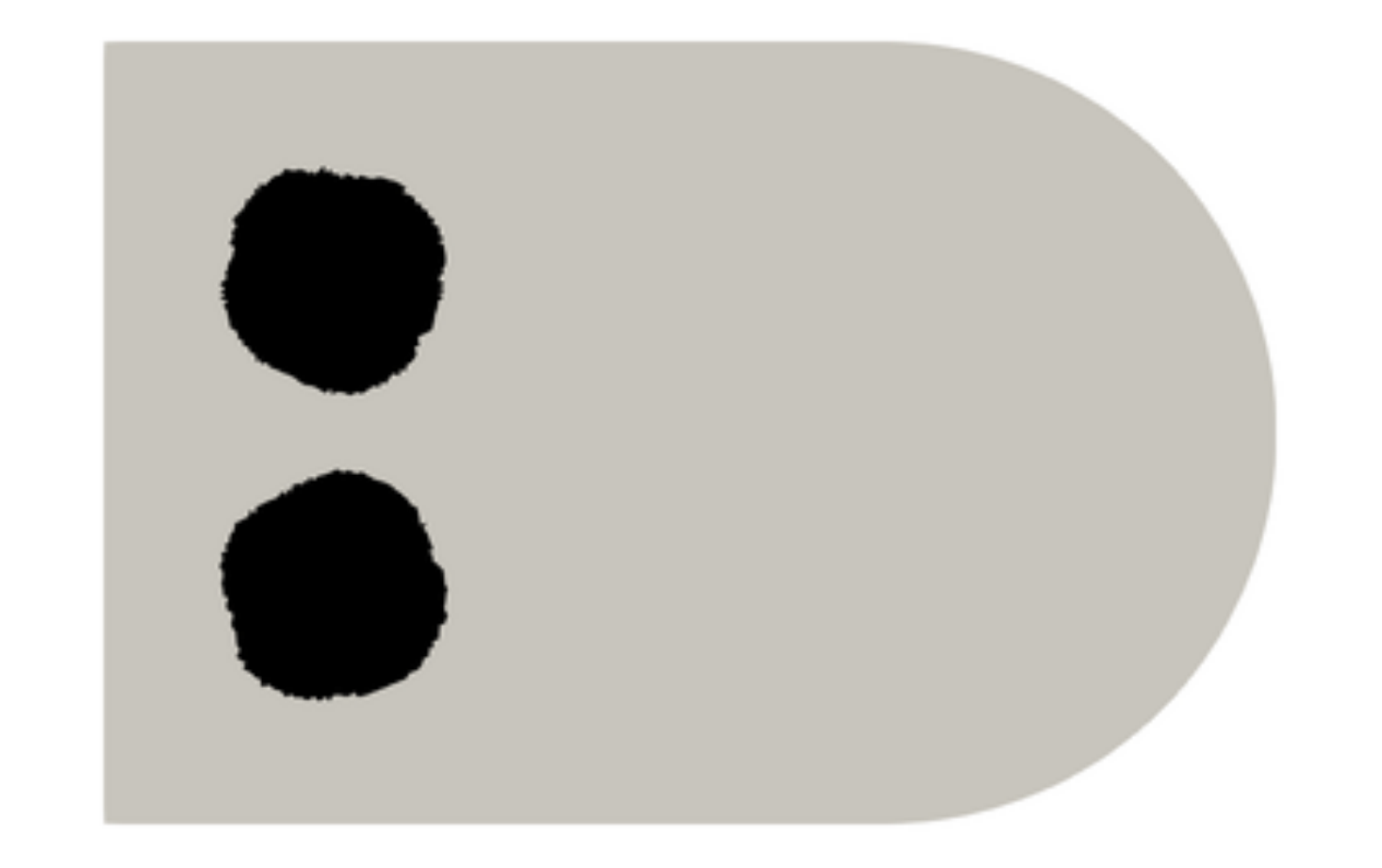}\quad
  \includegraphics[width=5cm]{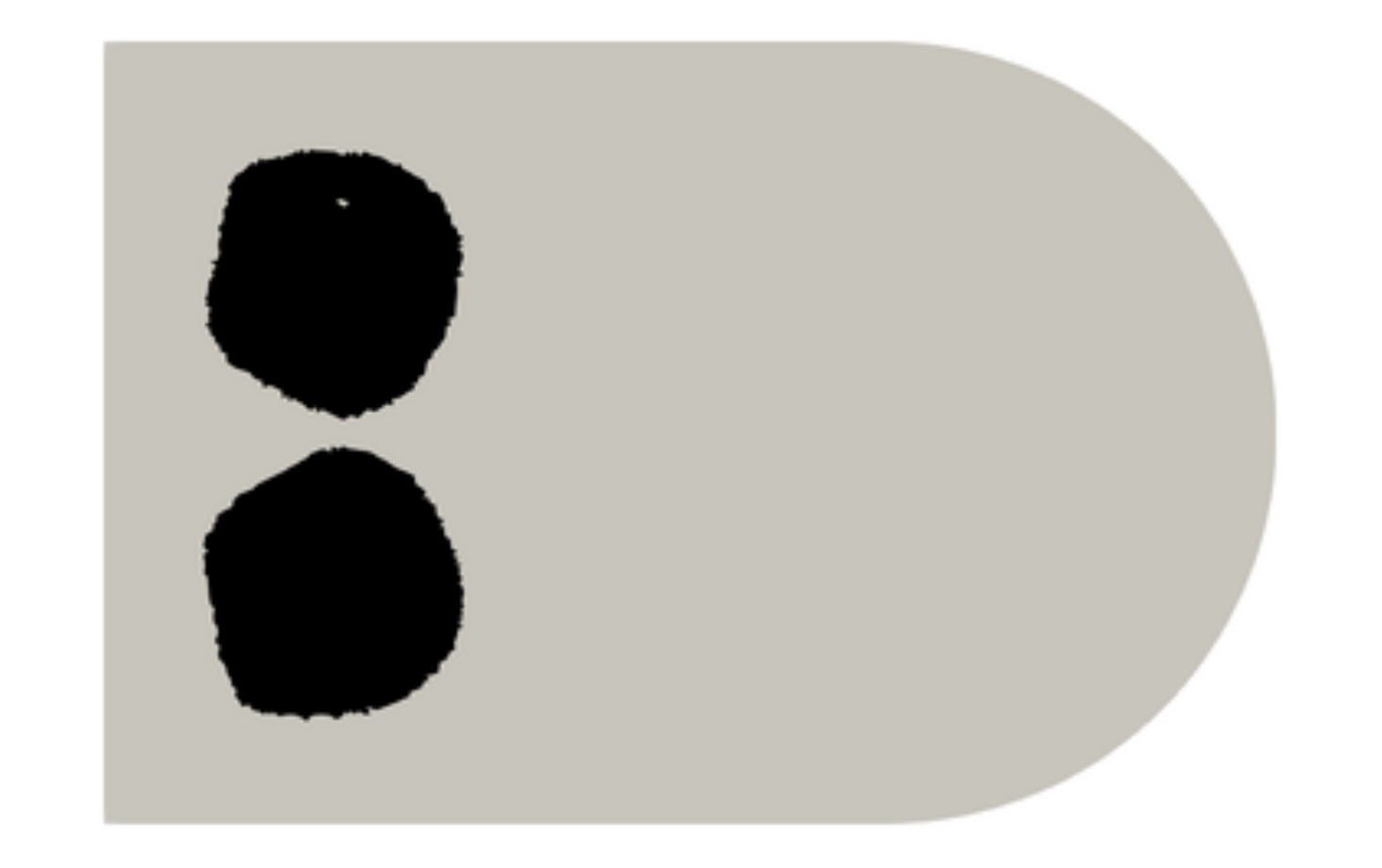}\\
  \includegraphics[width=5cm]{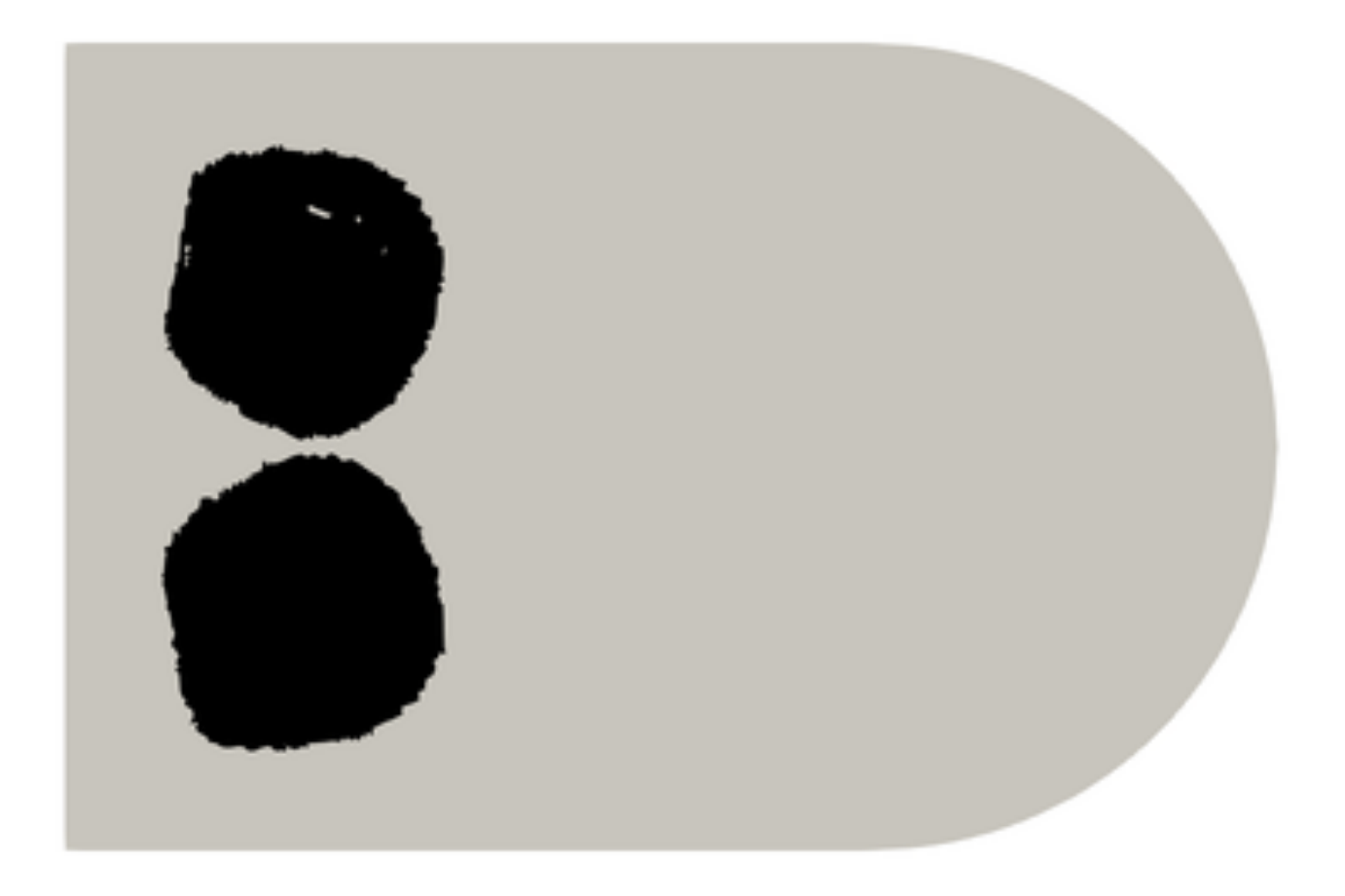}\quad
  \includegraphics[width=5cm]{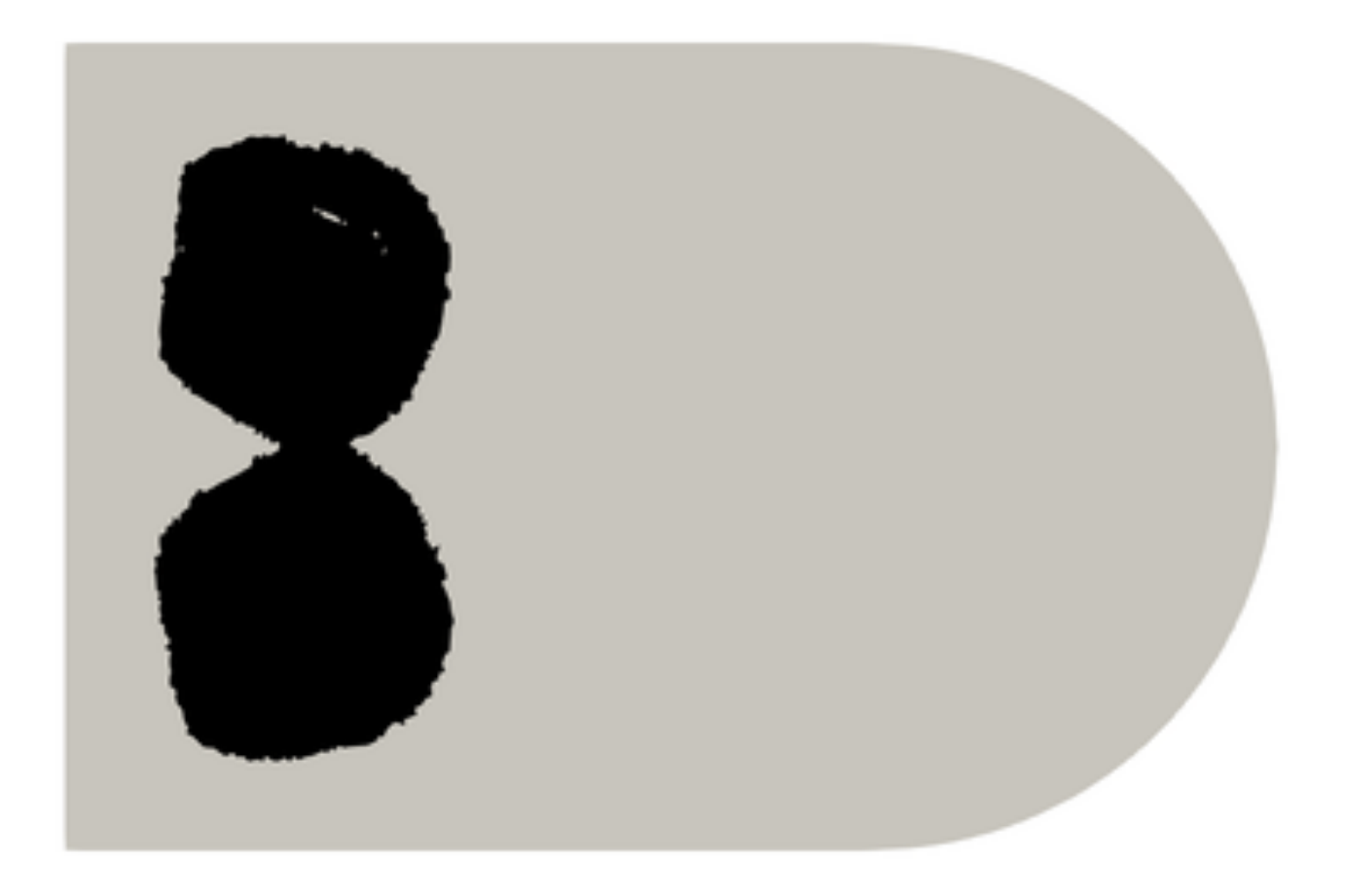}
\caption{Test 1. $\Omega$ at iterations: 0, 25, 50, 500, 1000, 2000.}
\label{fig:Omega}
\end{figure}

In Figure \ref{fig:Omega}, we can observe the evolution of $\Omega$.
It does not touch the boundary of $D$ and $g_h$ is not updated on $\partial D$.
 Then the constraint (\ref{4.8}), valid at the initial iteration, is still observed
at the following iterations. This numerical experiment was continued up to the iteration 2000 and
the numerical convergence was achieved in the sense that the difference
$\vert\mathcal{J}_h(X^k)-\mathcal{J}_h(X^{k-1})\vert$ is of
order $10^{-3}$ for $k \geq 1900$. 
The value $\mathcal{J}_h(X^{2000}) = 13.475690$. 
  In Figure \ref{fig:Omega}, one can see that
the two obstacles merge (starting with the iteration 1282, in fact). Notice the small holes on
one side of the obstacle which we consider as a computational disturbance. They are not admissible
and we neglect them. Moreover, this experiment and the next one show that our approach can
close/merge/generate holes combined with shape optimization (see \cite{MT_2022} as well,
where the Hamiltonian approach is employed). We also remark that, in general, 
$\| \mathcal{B}(G_h)Y_h \|_\infty$ 
varies between $0.1$ and $0.204184$ (mainly due to the violation of the discrete divergence condition), which is a weak point. 
In the next experiment, this aspect is much improved.

\bigskip
\textbf{Test 2.}
\medskip

We have also computed the optimization of the tracking type objective function 
$j\left(\mathbf{x}, \mathbf{y}(\mathbf{x})\right)
=\left(\mathbf{y}(\mathbf{x})-\mathbf{y}_d(\mathbf{x})\right)
\cdot\left(\mathbf{y}(\mathbf{x})-\mathbf{y}_d(\mathbf{x})\right)$
where $\mathbf{y}_d\in W$ is the solution obtained in Section \ref{sec:3.1}
but for
$$
g=-\frac{(x_1+0.2)^2}{0.2^2} - \frac{(x_2)^2}{0.4^2} + 1.
$$

The gradient of the objective function is now
$$
\nabla\mathcal{J}_h(X)=
\left(
\begin{array}{c}
  \int_D [1-\tilde{H}^h(g_h)] 2\left(\mathbf{y}^\epsilon_h-\mathbf{y}_d\right)\cdot
  \boldsymbol{\phi}_i d\mathbf{x},\quad 1\leq i\leq 2N_1\\
  0\in \mathbb{R}^{N_2}\\
  \int_D [-\tilde{H}^h]^\prime (g_h)\gamma_i\,
  \left(\mathbf{y}^\epsilon_h-\mathbf{y}_d\right)\cdot
  \left(\mathbf{y}^\epsilon_h-\mathbf{y}_d\right) d\mathbf{x},\quad 1\leq i\leq N_3  
\end{array}
\right).
$$

\begin{figure}[ht]
\centering
  \includegraphics[width=10cm]{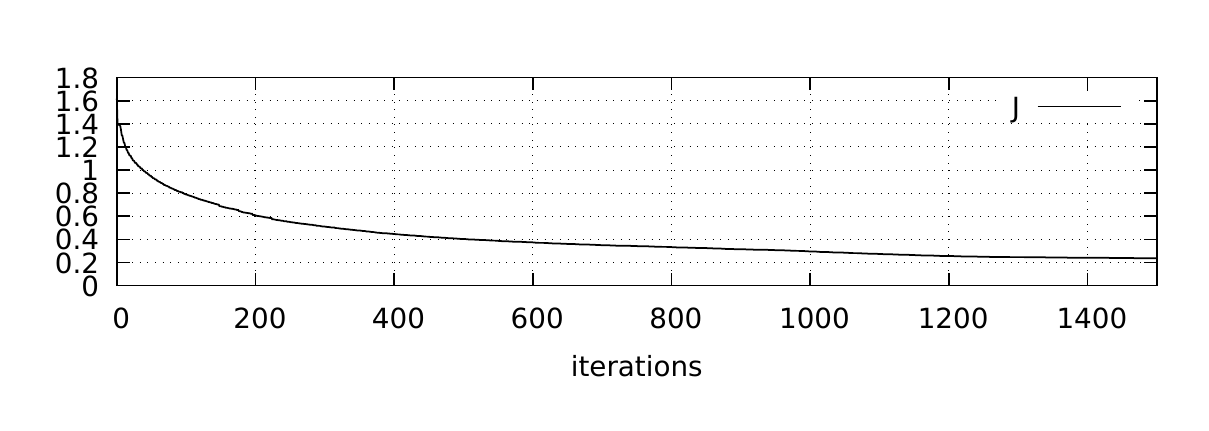}
\caption{Test 2. The history of $\mathcal{J}_h(X)$.}
\label{fig:grad_J}
\end{figure}

\begin{figure}[ht]
\centering
  \includegraphics[width=10cm]{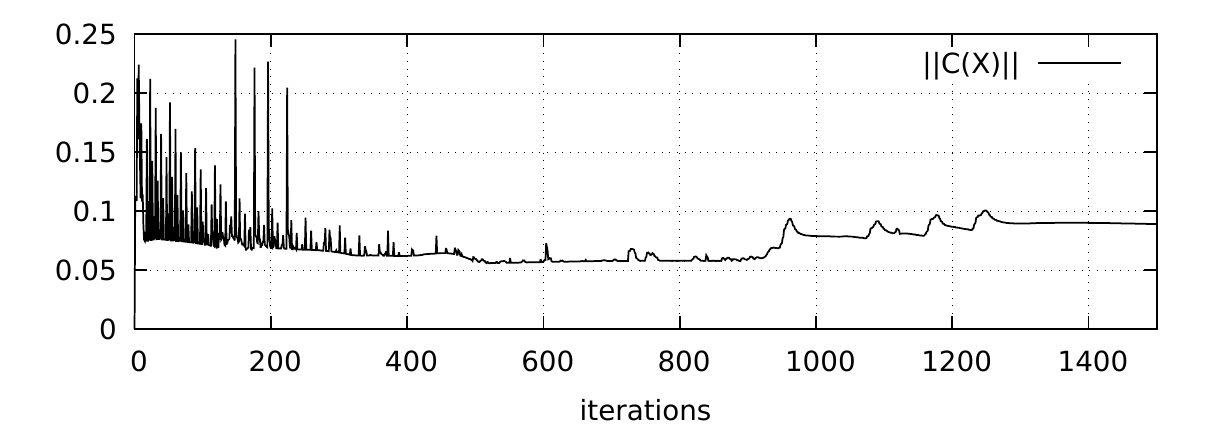}
\caption{Test 2. The history of $\| C(X) \|_\infty$.}
\label{fig:grad_CX}
\end{figure}

\begin{figure}[ht]
\centering
  \includegraphics[width=10cm]{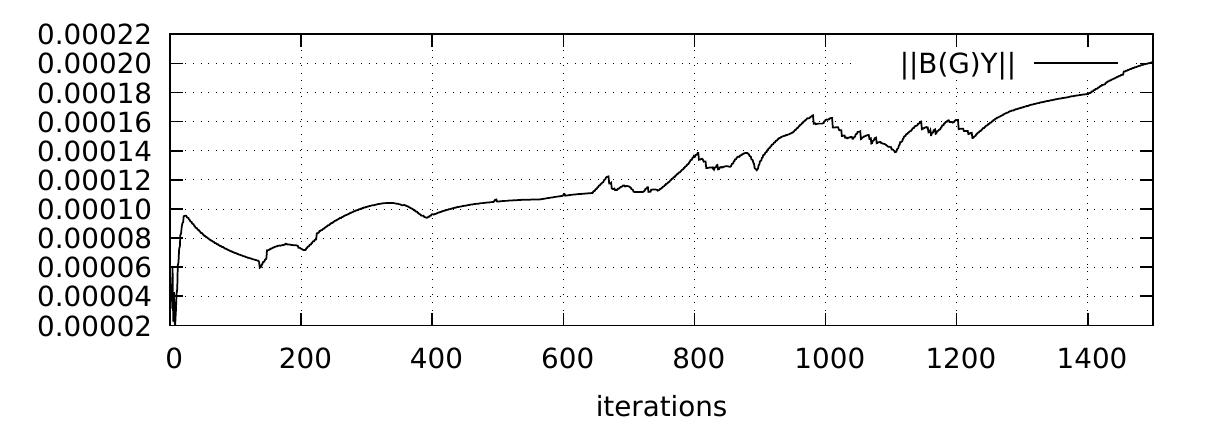}
\caption{Test 2. The history of $\| \mathcal{B}(G_h)Y_h \|_\infty$.}
\label{fig:grad_div}
\end{figure}

The numerical parameters are as before, but the initial domain $\Omega$
is composed by two disks $\Omega_1$ and $\Omega_2$ of radius $0.15$ and
centers $(-0.2,0.2)$ and $(-0.2,-0.2)$. The penalized parameter is $\rho=0.02$.
We can see the evolution of $\mathcal{J}_h$ and $C$ in Figure \ref{fig:grad_J}
and \ref{fig:grad_CX}.
At the initial iteration we have $\mathcal{J}_h=1.618213$ and
$\| C(X) \|_\infty=7.16e-10$. In this example, we know that the global optimal value is $0$ and we have performed more iterations as in Test 1 to obtain convergence, but the descent property is slow.
After 1500 iterations, we have $\mathcal{J}_h=0.237144$ and $\| C(X) \|_\infty=8.9e-02$.
In Figure \ref{fig:grad_div}, we can observe $\| \mathcal{B}(G_h)Y_h \|_\infty$ (the discrete divergence) that is close to $0$ in all iterations.

\begin{figure}[ht]
\centering
  \includegraphics[width=5cm]{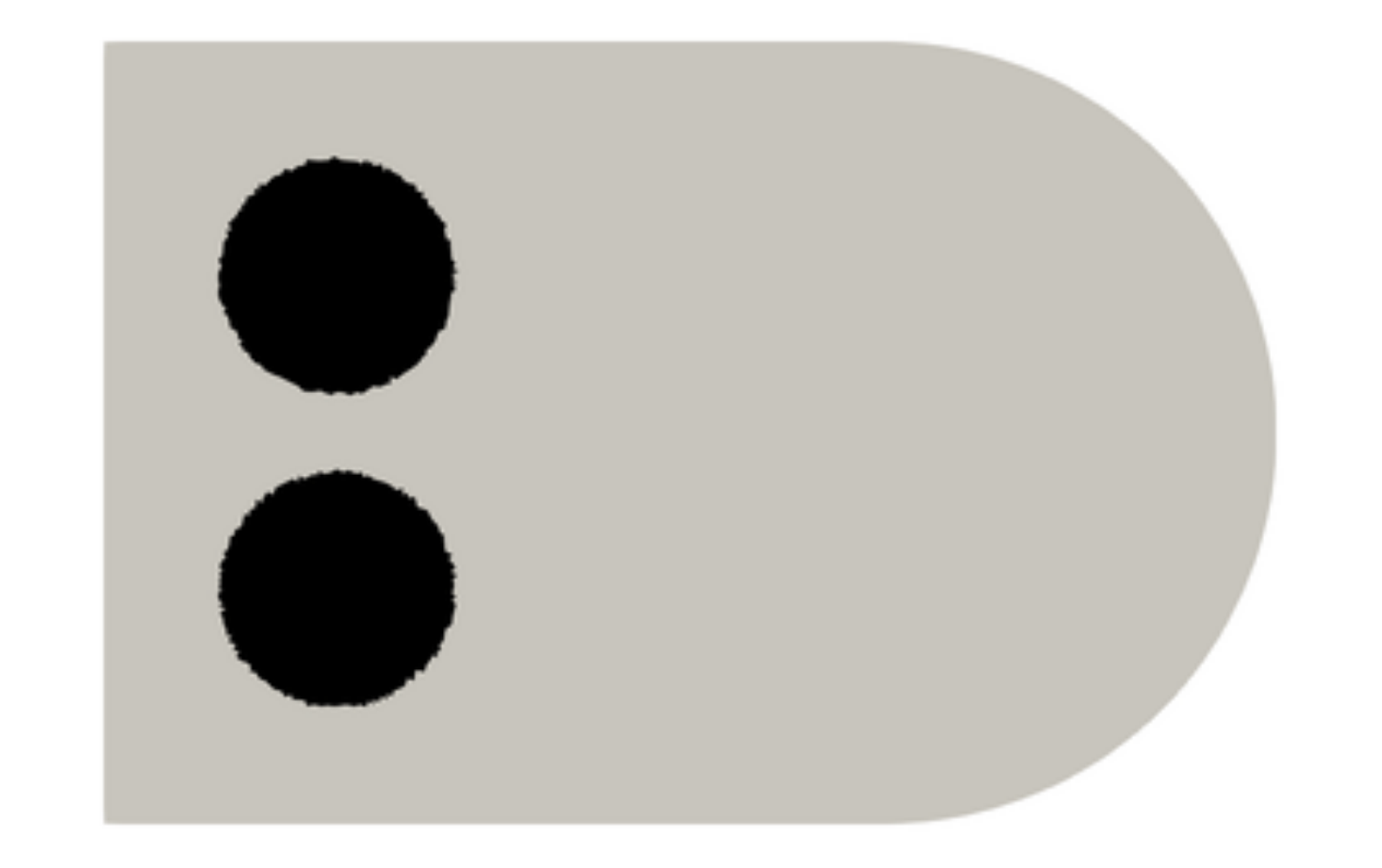}\quad
  \includegraphics[width=5cm]{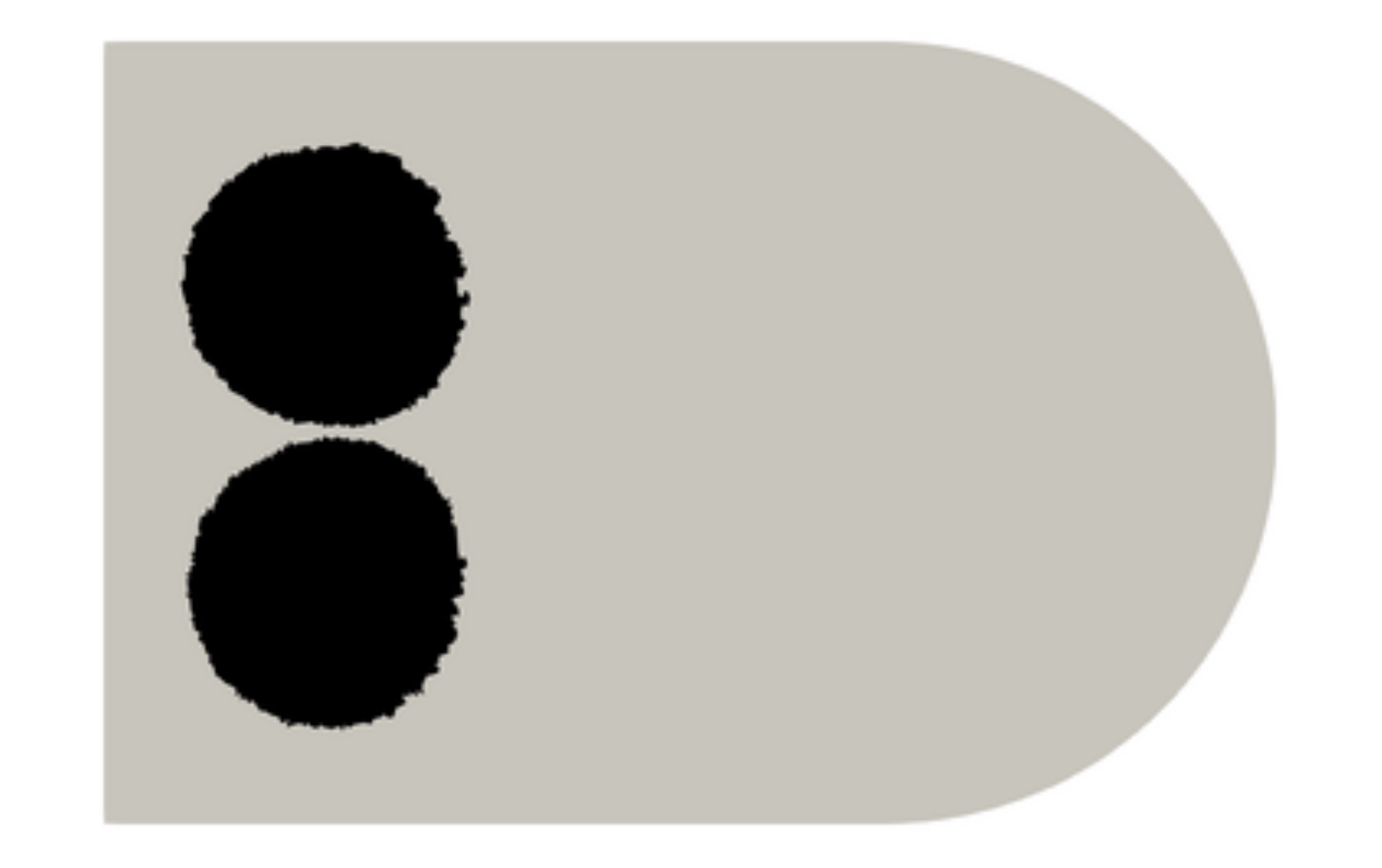}\\
  \includegraphics[width=5cm]{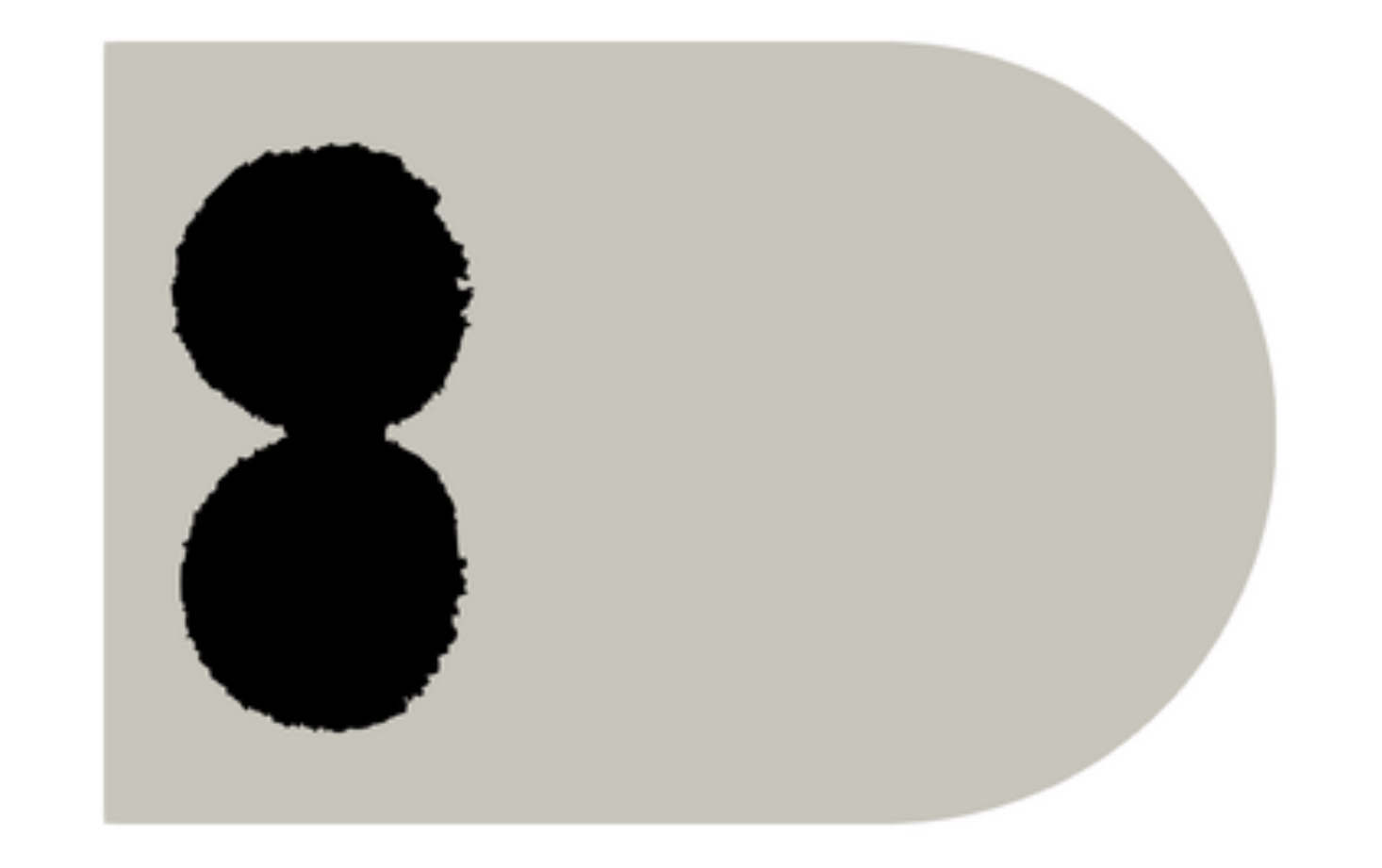}\quad
  \includegraphics[width=5cm]{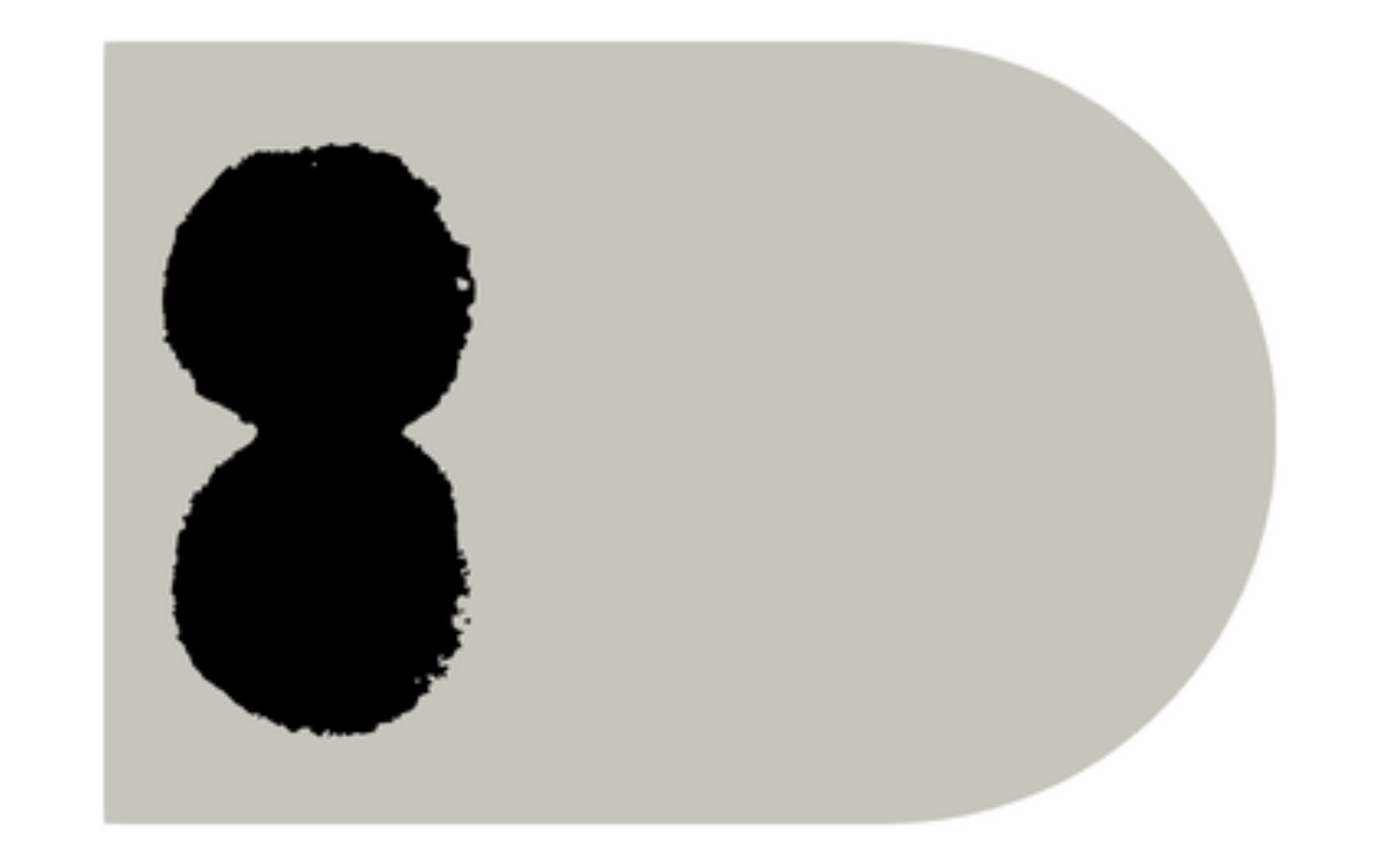}\\
  \includegraphics[width=5cm]{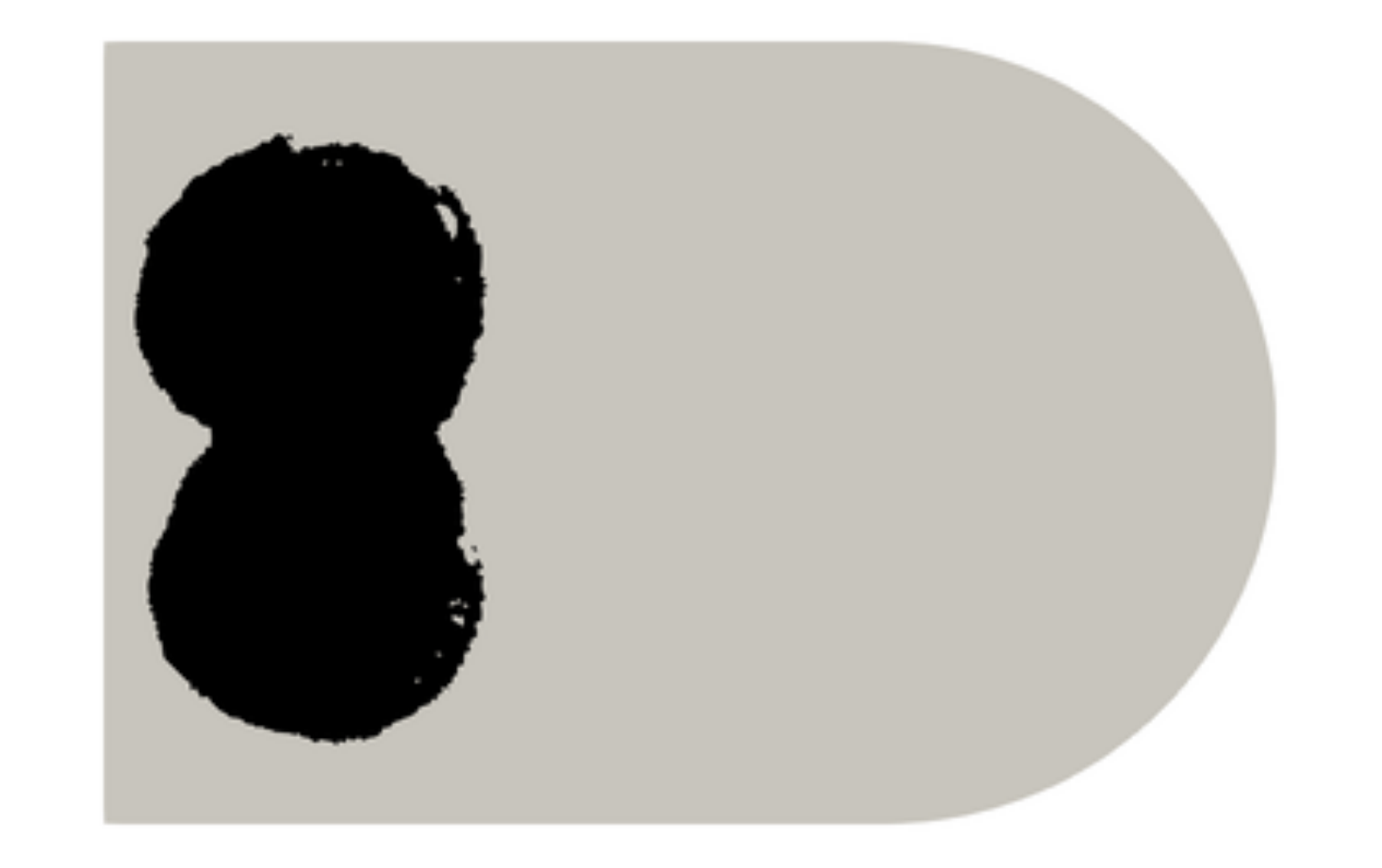}\quad
  \includegraphics[width=5cm]{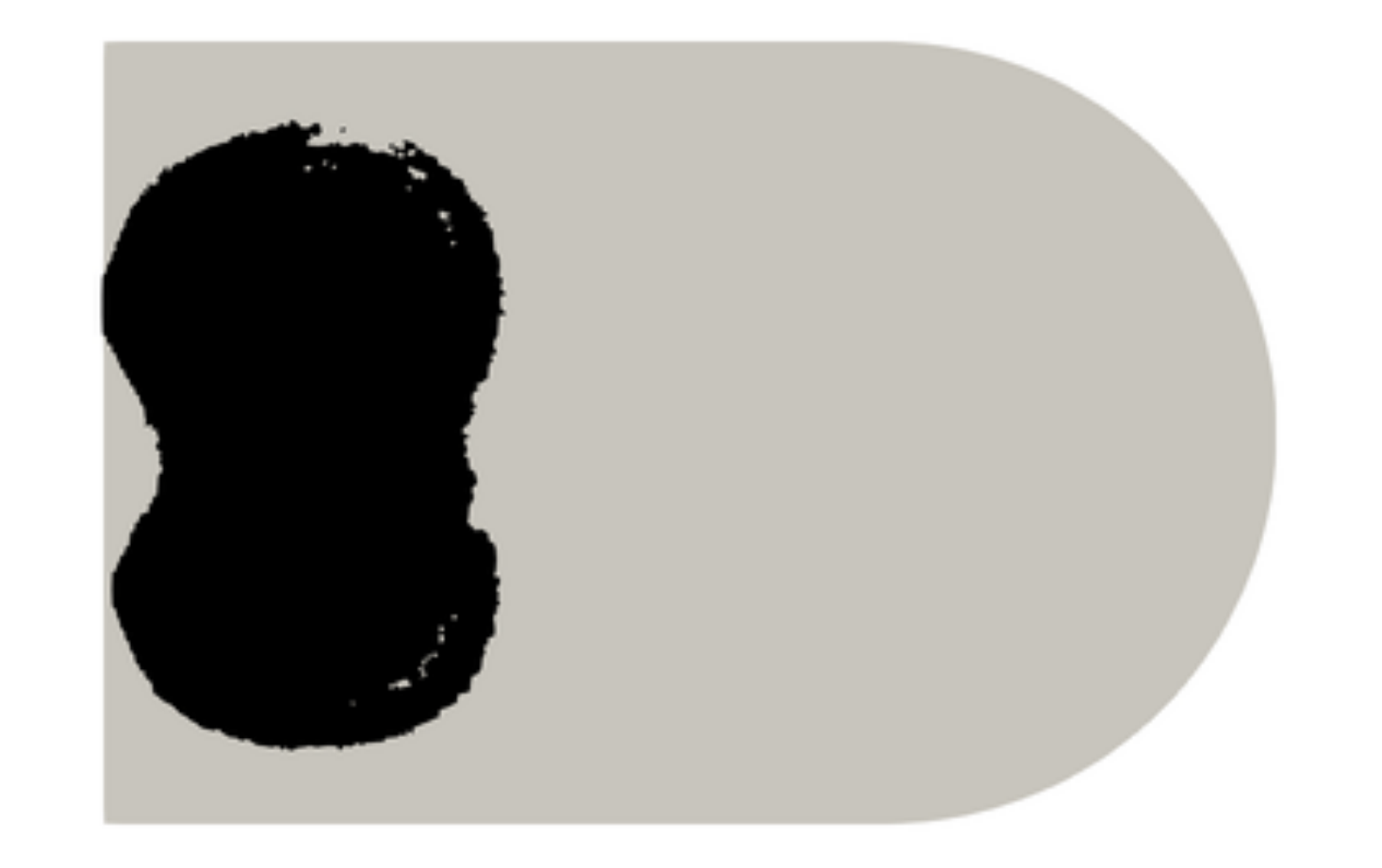}
\caption{Test 2. $\Omega$ at iterations: 0, 10, 25, 50, 200, 1500.}
\label{fig:grad_Omega}
\end{figure}

The evolution of $\Omega$ is presented in Figure \ref{fig:grad_Omega} (notice that the topology changes in this test too). After iteration 730 the obstacle touches $\partial D$, but this does not affect our numerical approach. The shape/topology optimization problems are strongly nonconvex and the global solution may be not unique.

In this test, if $\rho$ is 0.2 or bigger, the equality constraint is strongly enforced
(``stringent'', \cite{SN2022}) and the cost functional $\mathcal{J}_h(X)$ or the domain
$\Omega$ cannot change significantly  since the algorithm cannot find significant
admissible variations, see Remark 6 in \cite{SN2022}.

We end this subsection by presenting a convergence result.
As in Section \ref{sec:2.2}, we use the basis
$\{ \boldsymbol{\varphi}_i \}_{i=1,\dots,n} $ of $V_h$ and  we put
$m=N_3$, that is $G_h\in \mathbb{R}^m$. We introduce
\begin{eqnarray*}
A_h&=&\left( a^\epsilon_h(\boldsymbol{\varphi}_j,\boldsymbol{\varphi}_i) \right)_{1\leq i,j\leq n}
\in \mathbb{R}^{n\times n}\\
\tilde{C}_h(Y)&=&\left( \tilde{c}^\epsilon_h\left(\sum_{k=1}^n Y_k\boldsymbol{\varphi}_k,
\boldsymbol{\varphi}_j,\boldsymbol{\varphi}_i\right) \right)_{1\leq i,j\leq n}
\in \mathbb{R}^{n\times n},\quad Y= \left(Y_k\right)_{1\leq k\leq n}\in \mathbb{R}^n\\
\mathbf{f}_h&=&\left(\langle F_h, \boldsymbol{\varphi}_i\rangle_{*,D}\right)_{1\leq i\leq n} \in \mathbb{R}^n.
\end{eqnarray*}

The discrete minimization problem (\ref{3.20}), (\ref{3.21}) takes the form (the pressure doesn't appear explicitly anymore):
\begin{equation}\label{4.7}
\inf_{(G_h,Y_h)\in \mathbb{R}^m \times \mathbb{R}^n} J_h(G_h, Y_h)
\end{equation}
subject to the algebraic formulation (with $g$  replaced by $g_h$) of (\ref{3.ns3}):
find $Y_h \in \mathbb{R}^n$ such that
\begin{equation}\label{3.8}
A_h Y_h+ \tilde{C}_h(Y_h)Y_h=\mathbf{f}_h 
\end{equation}
(and the constraint (\ref{4.8}), neglected in the numerical tests).
Notice that $A_h$, $\tilde{C}_h$, $\mathbf{f}_h$, $Y_h$, $\mathbf{y}^\epsilon_h$ depend on $g_h$ (or $G_h$).

The standard penalization technique eliminates the equality constraint (\ref{3.8}) ($\delta > 0$ is
a ``small'' parameter):
\begin{equation}\label{3.12}
\inf_{(G_h,Y_h)\in \mathbb{R}^m \times \mathbb{R}^n} J_h(G_h,Y_h)
+ \frac{1}{\delta} \vert A_h Y_h+ \tilde{C}_h(Y_h)Y_h-\mathbf{f}_h \vert_{\mathbb{R}^n} ^2.
\end{equation}

The problem (\ref{3.12}), (\ref{4.8}) may again have no global solution
and we shall work with minimizing sequences. We assume admissibility, that is the existence
of a pair $[\hat{G}_h, \hat{Y}_h] \in \mathbb{R}^m \times \mathbb{R}^n$ with
$J_h (\hat{G}_h, \hat{Y}_h)$ finite and such that (\ref{3.8}), (\ref{4.8}) are satisfied.
We have the inequality:
\begin{equation}\label{3.13}
\inf \{(\ref{3.12}), (\ref{4.8})\} \leq \inf \{(\ref{4.7}), (\ref{3.8}), (\ref{4.8})\},
\end{equation}
for any $\delta > 0$. This follows since for any admissible $[G_h, Y_h]$ for
(\ref{4.7}), (\ref{3.8}), (\ref{4.8}), we have
$$
\inf\{(\ref{3.12}), (\ref{4.8})\} \leq J_h(G_h, Y_h)
$$
and (\ref{3.13}) is obtained by taking infimum in the right-hand side.
Consequently, for any $\delta \in (0,1]$,
there is $[G_{\delta}, Y_{\delta}] \in \mathbb{R}^m \times \mathbb{R}^n$ satisfying (\ref{4.8}),
such that we have:
\begin{eqnarray}
&& J_h(G_\delta,Y_\delta)
  + \frac{1}{\delta} \vert A_h Y_\delta+ \tilde{C}_h(Y_\delta)Y_\delta -\mathbf{f}_h \vert_{\mathbb{R}^n} ^2
  \nonumber\\
  &\leq &  \inf \{(\ref{3.12}), (\ref{4.8})\} + \delta
  \leq  \inf \{(\ref{4.7}), (\ref{3.8}), (\ref{4.8})\} + \delta
  \leq C
\label{3.14}
\end{eqnarray}
where $C > 0$ is an absolute constant, independent of $\delta \in (0,1]$.

The set of admissible $G_h \in \mathbb{R}^m$, i.e. satisfying (\ref{4.8}), may be restricted
to be bounded since in (\ref{4.3}) or (\ref{3.omega}) one may scale $g$ by positive
constants without modifying the corresponding geometry $\omega_g$. Moreover, we assume
coercivity of the cost functional, with respect to $Y_h$:
\begin{equation}\label{3.15}
J_h(G_h,Y_h) \rightarrow \infty \ \hbox{ if }\ \vert Y_h\vert_{\mathbb{R}^n} \rightarrow \infty,
\end{equation}
uniformly with respect to admissible $G_h$ (bounded). This hypothesis (\ref{3.15}) is clearly
satisfied by many cost functionals of interest. By (\ref{3.14}), (\ref{3.15}) we get
that $[G_{\delta}, Y_{\delta}]$ are bounded in $\mathbb{R}^m \times \mathbb{R}^n$ with respect
to $\delta \in (0,1]$. We also impose that $J_h(\cdot, \cdot)$ is bounded from below by some constant.
Then, (\ref{3.14}) yields 
\begin{equation}\label{3.16}
\vert A_h Y_\delta+ \tilde{C}_h(Y_\delta)Y_\delta -\mathbf{f}_h \vert_{\mathbb{R}^n} ^2 \leq M\delta,
\end{equation}
with $M > 0$ another absolute constant independent of $\delta \in (0,1]$. By (\ref{3.15}),
we get $\{Y_{\delta}\}$ bounded in $\mathbb{R}^n$ with respect to $\delta \in (0,1]$ and, on a
subsequence, we may assume $G_{\delta} \rightarrow \tilde{G}$ in
$\mathbb{R}^m$, $Y_{\delta} \rightarrow \tilde{Y}$ in $\mathbb{R}^n$.
We recall that $A_h=A_h(G_h)\in \mathbb{R}^{n\times n}$ depends on $G_h$,
$\tilde{C}_h(Y_h)=\tilde{C}_h(G_h,Y_h)\in \mathbb{R}^{n\times n}$ depends on $G_h$ and $Y_h$ and
$\mathbf{f}_h=\mathbf{f}_h(G_h)\in \mathbb{R}^n$ depends on $G_h$.
By (\ref{3.16}), due to the continuity properties of $a^\epsilon_h$ and $\tilde{c}^\epsilon_h$ defined in
Subsection \ref{sec:2.2}, 
we obtain that $[\tilde{G}, \tilde{Y}]$ satisfies (\ref{3.8}).
However, $\tilde{G}$ may satisfy
(\ref{4.8}) just with the $\leq$ sign.

Another standard assumption that we require for $J_h(\cdot, \cdot)$ is lower semicontinuity.
Then, if $\delta \rightarrow 0$, relation (\ref{3.14}) gives
\begin{equation}\label{3.17}
  J_h(\tilde{G},\tilde{Y}) \leq \inf\{(\ref{3.12}),(\ref{4.8})\}
  \leq \inf\{(\ref{4.7}), (\ref{3.8}), (\ref{4.8})\}.
\end{equation}
We have proved

\begin{proposition}\label{prop:3.a}
Assume admissibility for the constrained optimization problem (\ref{4.7}), (\ref{3.8}), (\ref{4.8}) and boundedness for the admissible controls $G_h$. If the cost functional is majorized from below, lower semicontinuous and coercive as in (\ref{3.15}), then the sequence $[G_{\delta}, Y_{\delta}]$ satisfies (\ref{4.8}), has the
minimizing property from (\ref{3.17}) and satisfies (\ref{3.8}) in the approximating
variant (\ref{3.16}).
\end{proposition}

Notice that the assumptions from Prop. \ref{prop:3.a} are quite standard and ``realistic''. Their use via penalization, in shape/topology optimization, seems new. In case the usual FEM basis of hat functions $\{ \boldsymbol{\phi}_i \}_{i=1,\dots,2N_1} $ of $W_h$ is employed, the constrained minimization problem has the form (\ref{3.20}), (\ref{3.21}) with unknowns $Y_h$, $P_h$, $G_h$, and both components of $C(X)$ have to be penalized in the cost as in (\ref{A}).

\section*{Acknowledgement}
This work was partially supported by the French - Romanian cooperation program ``ECO Math'', 2022.

\end{document}